\documentclass[12pt]{article}
\usepackage{amsmath,amsthm,amssymb,amscd}
\usepackage{latexsym}
\newtheorem{thm}{Theorem}[section]
 \newtheorem{cor}[thm]{Corollary}
 \newtheorem{lem}[thm]{Lemma}
 \newtheorem{prop}[thm]{Proposition}
 \theoremstyle{definition}
 
 \theoremstyle{remark}

  \numberwithin{equation}{section}

\makeatletter
\newcommand*\owedge{\mathpalette\@owedge\relax}
\newcommand*\@owedge[1]{%
  \mathbin{%
    \ooalign{%
      $#1\m@th\bigcirc$\cr
      \hidewidth$#1\m@th\wedge$\hidewidth\cr
    }%
  }%
}
\makeatother

\begin{document}

\title
{Compact manifolds of dimension $n\geq 12$ with positive isotropic curvature}

\author{ Hong Huang}
\date{}
\maketitle

\begin{abstract}
 We prove the following result: Let $(M,g_0)$ be a compact
manifold of dimension $n\geq 12$ with  positive isotropic curvature. Then $M$ is diffeomorphic  to a  spherical space form, or the total space of an orbifiber bundle over $\mathbb{S}^1$ or $\mathcal{I}$ with generic fiber diffeomorphic to $\mathbb{S}^{n-1}/\Gamma$ such that the total space  admits a metric with positive isotropic curvature, where $\Gamma$ is a finite subgroup of $O(n)$ acting freely on $\mathbb{S}^{n-1}$, and $\mathcal{I}$ is the one dimensional closed orbifold with two singular points both with  local group $\mathbb{Z}_2$ and with $|\mathcal{I}|$  a closed interval,  or  a  connected sum of a finite number of such manifolds.
This extends a recent work of Brendle, and implies a conjecture of Schoen and a conjecture of Gromov in dimensions $n\geq 12$.  The proof uses Ricci flow with surgery on compact orbifolds with isolated singularities.


\end{abstract}


\section {Introduction}

The notion of positive isotropic curvature was introduced by Micallef and Moore \cite{MM}. Recall  that a Riemannian manifold $M$ (of dimension $\geq 4$) is said to have
 positive isotropic curvature ((strictly) PIC) if  for all points $p \in M$ and all
orthonormal 4-frames $\{e_1,e_2,e_3,e_4\}\subset T_pM$, the curvature
tensor satisfies
\begin{equation*}R_{1313}+R_{1414}+R_{2323}+R_{2424}-2R_{1234} > 0.
\end{equation*}

The PIC condition is preserved under the Ricci flow, as shown by Hamilton \cite{H97} in dimension 4, and by Brendle-Schoen \cite{BS} and Nguyen \cite{N} in dimensions $n\geq 5$ (see also \cite{Wi}). Hamilton \cite{H97} also derived curvature pinching estimates for Ricci flow with strictly PIC in dimension 4 and initiated the study of Ricci flow with surgery (see also \cite{P2} and \cite{CZ}). Chen-Tang-Zhu \cite{CTZ} used Ricci flow with surgery on orbifolds with isolated singularities to give a topological classification of compact 4-manifolds with strictly PIC. (See also \cite{Hu13} and \cite{Hu15} for some extensions of \cite{H97}, \cite{CZ} and \cite{CTZ}.)

In a recent breakthrough \cite{B19} Brendle obtained  curvature pinching estimates for Ricci flow on compact
manifolds of dimension $n\geq 12$ with  positive isotropic curvature,  constructed Ricci flow with surgery on these manifolds with an additional topological constraint, that is, not containing any nontrivial incompressible $(n-1)$-dimensional space forms, and classified them.    In this paper we  investigate the case without this extra constraint, and via constructing Ricci flow with surgery on compact orbifolds of dimension $n\geq 12$ with isolated singularities and with positive isotropic curvature, and using some topological arguments, get  the following result.

 \begin{thm} \label{thm 1.1} Let $(M,g_0)$ be a compact
manifold of dimension $n\geq 12$ with  positive isotropic curvature. Then $M$ is diffeomorphic  to a  spherical space form, or the total space of an orbifiber bundle over $\mathbb{S}^1$ or $\mathcal{I}$ with generic fiber diffeomorphic to $\mathbb{S}^{n-1}/\Gamma$ such that the total space  admits a metric with positive isotropic curvature, where $\Gamma$ is a finite subgroup of $O(n)$ acting freely on $\mathbb{S}^{n-1}$, and $\mathcal{I}$ is the one dimensional closed orbifold with two singular points both with  local group $\mathbb{Z}_2$ and with $|\mathcal{I}|$  a closed interval,  or  a  connected sum of a finite number of such manifolds.
\end{thm}

 \noindent (By  \cite{MW} the converse of Theorem 1.1 is also true.) In fact, from the proof of Theorem 1.1 we know that a component appearing in the connected sum decomposition of $M$ in Theorem 1.1 which  is diffeomorphic  to  the total space of an orbifiber bundle over $\mathbb{S}^1$ or $\mathcal{I}$ with generic fiber diffeomorphic to $\mathbb{S}^{n-1}$ is diffeomorphic to $\mathbb{R}P^n \sharp \mathbb{R}P^n$, or $\mathbb{S}^{n-1} \times \mathbb{S}^1$, or $\mathbb{S}^{n-1} \tilde{\times} \mathbb{S}^1$, and is diffeomorphic to a  quotient manifold of $\mathbb{S}^{n-1} \times \mathbb{R}$ by standard isometries.  Here $\mathbb{S}^{n-1} \tilde{\times} \mathbb{S}^1$ is the (unique) nonorientable $\mathbb{S}^{n-1}$-bundle over $\mathbb{S}^1$ with structure group $O(n)$. So Theorem 1.1  extends Theorem 1.4 in \cite{B19}. (See also the remark after the proof of Theorem 1.1.)  Compare also the Main Theorem in \cite{CTZ}; here I don't know that if a $\mathbb{S}^{n-1}/\Gamma$-bundle ($|\Gamma|\geq 2$) over $\mathbb{S}^1$ or $\mathcal{I}$ admits a metric with positive isotropic curvature, whether or not the structure group can always be reduced to the isometry group Isom$(\mathbb{S}^{n-1}/\Gamma)$. (I conjecture that if a $\mathbb{S}^{n-1}/\Gamma$-bundle ($n\geq 5, |\Gamma|\geq 2$) over $\mathbb{S}^1$ or $\mathcal{I}$ admits a metric with positive isotropic curvature,  the structure group can always be reduced to Isom$(\mathbb{S}^{n-1}/\Gamma)$.)
 On the other hand, Chen-Tang-Zhu predicted on p. 47 of \cite{CTZ}  that a  weaker version of Theorem 1.1 should hold  once Hamilton's curvature pinching estimates in \cite{H97} can be extended to higher dimensions.  For the definition of  connected sum of smooth manifolds see for example,  \cite{B}, \cite{BJ},  \cite{K}, \cite{M07}, \cite{Mu} and \cite{W}.

 Theorem 1.1 implies  a  conjecture of R. Schoen \cite{S} in dimension $n\geq 12$, which in turn implies a conjecture of M. Gromov \cite{G} (see also p. 346 in \cite{F};  note that the main theorem of \cite{F} may be viewed as a partial result on Gromov's conjecture for $n\geq 5$) in dimension $n\geq 12$; the 4-dimensional case of the Schoen conjecture was solved in Corollary 1 in \cite{CTZ}.

\begin{cor} \label{cor 1.2}  Let $(M,g_0)$ be a compact
manifold of dimension $n\geq 12$ with  positive isotropic curvature. Then there is a finite cover of $M$ which is diffeomorphic  to $\mathbb{S}^n$, or $\mathbb{S}^{n-1} \times \mathbb{S}^1$, or a connected sum of a finite number of copies of $\mathbb{S}^{n-1} \times \mathbb{S}^1$. In particular, the fundamental group of $M$ is virtually free.
\end{cor}

In the process of proof of Theorem 1.1 we get a slightly more general theorem.

\begin{thm} \label{thm 1.3} \ \  Let $(\mathcal{O},g_0)$ be a compact
orbifold of dimension $n\geq 12$ with at most isolated singularities and with  positive isotropic curvature. Then $\mathcal{O}$ is diffeomorphic to  a spherical orbifold with at most isolated singularities, or  a weak connected sum of  at most two spherical orbifolds with at most isolated singularities  which admits a metric with positive isotropic curvature,  or a connected sum of a finite number of such orbifolds.
\end{thm}

\noindent (By an obvious extension of Theorem 1.1 in \cite{MW} to the case of orbifold connected sum (cf. also \cite{CH}), certain version of the converse of Theorem 1.3 is also true.) Compare Theorem 2.1 in \cite{CTZ}. The  weak connected sum of orbifolds is defined as follows. Let $\mathcal{O}_i$ ($i=1,2$) be two connected (effective)  $n$-orbifolds with at most isolated singularities, and let
$D_i \subset \mathcal{O}_i$ be two embedded suborbifolds-with boundary, both
diffeomorphic to a quotient orbifold $D^n//\Gamma$, where $D^n$
is a closed, standard $n$-disk (contained in the Euclidean space $\mathbb{R}^n$), and $\Gamma$ is a finite subgroup of
$O(n)$ acting freely on $\mathbb{S}^{n-1}$. We glue  together
$\mathcal{O}_1\setminus \text{Int} \hspace*{0.5mm} D_1 $ and  $\mathcal{O}_2\setminus \text{Int} \hspace*{0.5mm} D_2$ along their boundaries  using a diffeomorphism $f$ of their boundaries,  where $\text{Int} \hspace*{0.5mm} D_i$ denotes the interior of $D_i$, $i=1,2$. Of course in the case that both $\mathcal{O}_i$ are orientable and oriented, we should take care of the orientations, that is, in this case $f$ should be orientation reversing. The result is called the weak connected sum of
$\mathcal{O}_1$ and $\mathcal{O}_2$ via the gluing map $f$, and is
denoted by $\mathcal{O}_1 \tilde{\sharp}_f  \mathcal{O}_2 $.  (Compare Section 7.5.2 of \cite{B}.) If $ D_i$
($i=1,2$) are disjoint embedded suborbifolds-with boundary (both
diffeomorphic to the quotient orbifold $D^n//\Gamma$) in the same connected
$n$-orbifold $\mathcal{O}$, the result of similar process as above
is called a weak connected sum on the single orbifold
$\mathcal{O}$, and is denoted by $\mathcal{O} \tilde{\sharp}_f$.

By the proof of Theorem 2.3 in Chapter 8 of \cite{Hi}, given $\mathcal{O}_1$, $\mathcal{O}_2 $ (with their orientations if both are orientable), $D_1$  and $D_2$, the diffeomorphism type of $\mathcal{O}_1 \tilde{\sharp}_f  \mathcal{O}_2 $   depends only on the isotopy class of $f$, and given $\mathcal{O}$, $D_1$  and $D_2$, the diffeomorphism type of $\mathcal{O} \tilde{\sharp}_f$   depends only on the isotopy class of $f$.

 The  weak connected sum here is called orbifold connected sum in  \cite{CTZ} (see also \cite{Hu15}). We also call the similar operation defined in p. 48 of \cite{CTZ} involving more orbifolds and more gluing maps a weak connected sum.

Let's compare the weak connected sum and connected sum of smooth orbifolds. Recall the  connected sum of two  smooth orbifolds as defined in for example, p. 135-136 in Bonahon \cite{B02}. Let $\mathcal{O}_i$ ($i=1,2$) be two (connected) $n$-orbifolds with at most isolated singularities, and let  $f_i: D^n // \Gamma \rightarrow  \mathcal{O}_i$ be two embeddings,  where $\Gamma$ is a finite subgroup of
$O(n)$ acting freely on $\mathbb{S}^{n-1}$. Let $f=f_2 \circ f_1^{-1}|_{f_1(\mathbb{S}^{n-1}/\Gamma)}$.  If both $\mathcal{O}_i$ are orientable and oriented, suppose $f_1$ is orientation preserving and $f_2$ is orientation reversing. Let $\mathcal{O}= (\mathcal{O}_1\setminus f_1(\text{Int} (D^n // \Gamma))) \cup_f (\mathcal{O}_2\setminus f_2(\text{Int}(D^n // \Gamma)))$.
   The result is called the  connected sum of
$\mathcal{O}_1$ and $\mathcal{O}_2$ via the gluing map $f$, and is
denoted by $\mathcal{O}_1 \sharp_f  \mathcal{O}_2 $. (Compare Kleiner-Lott \cite{KL}.) If $ D_i$
($i=1,2$) are disjoint embedded suborbifolds-with boundary (both
diffeomorphic to the quotient orbifold $D^n//\Gamma$) in the same (connected)
$n$-orbifold $\mathcal{O}$, the result of similar process as above
is called a connected sum on the single orbifold $\mathcal{O}$  (or self-connected sum of $\mathcal{O}$), and is denoted by $\mathcal{O}\sharp_f$.
   In contrast, in the definition of weak connected sum  the gluing map $f$  is not required to be the restriction of a diffeomorphism from  $D_1$ to $D_2$.  (Of course, both the  connected sum and the weak connected sum can also be defined for  orbifolds with not necessarily isolated singularities.)

   From the existence of exotic spheres  we know that for an arbitrary self diffemorphism $f$ of $\mathbb{S}^{n-1}$ with $n\geq 7$,  $f$ may not  be extended to a self diffemorphism of $D^n$ in general.  A closely related fact is that, in general, a diffeomorphism of a $(n-1)$-dimensional spherical space form ($n\geq 5$) may not be isotopic to an isometry. Indeed this is one of  several difficulties  one has to face in constructing Ricci flows with surgery on compact $n$-manifolds ($n\geq 5$) with positive isotropic curvature and classifying these manifolds (of course, to get suitable curvature pinching estimates is another difficulty, which was overcome by Brendle \cite{B19} in dimensions $n\geq 12$).  Under the extra condition of having no nontrivial incompressible $(n-1)$-dimensional space form, Brendle \cite{B19} overcame this difficulty by using  a clever surgery argument (see the proof of Proposition 6.17 in \cite{B19}).

Without the extra condition of having no nontrivial incompressible $(n-1)$-dimensional space form we have to face some other difficulties compared to the case of \cite{B19}. As in \cite{CTZ}, we have more types of caps compared to that in \cite{B19}. Inspired by Definition 6.16 in \cite{B19} we impose  some extra properties in  the  definition of $\varepsilon$-caps when the $\varepsilon$-caps have ends diffeomorphic to $\mathbb{S}^{n-1}\times [0,1)$; accordingly, we need to establish certain extra properties for the corresponding (orbifold) ancient $\kappa$-solutions  and smooth standard solution, for which we use a surgery argument adapted from the proof of Proposition 6.17 in \cite{B19} with the help of the strong maximum principle and Cerf-Palais disk theorem. (See the proof of Propositions \ref{prop 3.5} and \ref{prop 3.6}.) These extra properties  of the $\varepsilon$-caps are important for us to recognize  the topology of some components which are removed in the Ricci flow surgery; the isotopy extension theorem and Hamilton's  canonical parametrization  are also useful here.

 As already observed in Section D in \cite{H97}, in the general case, the Ricci flow surgeries on  neck-pinch singularities  may introduce  orbifold singularities even if we start with a manifold. In the orbifold case, a priori, an $\varepsilon$-neck or $\varepsilon$-tube may not be universally noncollapsed; this may lead to  difficulty in some contradiction-via-compactness argument. Chen-Tang-Zhu \cite{CTZ} dealt with this difficulty by using the technique of pulling back the Ricci flow solutions restricted to these domains via Hamilton's canonical uniformizations (cf. Section C.2 in \cite{H97} and the appendix of our paper); see the proof of Propositions 4.2 and 4.4 in \cite{CTZ}.  Of course we can also use arbitrary (smooth) universal covering maps to pull back the solutions restricted to these domains.
 Now by using Hamilton's compactness theorem we are in the following typical situation:  We have a strong $\varepsilon'$-neck upstairs, which is located at some covering space and which may not be equivariant under the action of the deck transformation groups  (recalling that the definition of Cheeger-Gromov convergence involves diffeomorphisms), but we want to get some strong $\varepsilon$-neck downstairs. Of course we may use Hamilton's canonical uniformization  in the final time slice, which is usually equivariant with respect to  the action of the deck transformation groups and descends to give an $\tilde{\varepsilon}$-neck downstairs  in the final time slice. But note that a Hamilton's canonical parametrization in the final time slice may not be such one in the previous time slices. Does Hamilton's canonical parametrization in the final time slice give a strong $\varepsilon$-neck downstairs? By using the argument in Section 2 of \cite{H95} one sees that it does. (This little detail was not written down in \cite{CTZ}. By the way, the word ``strong" is missing in the statement and proof of  Propositions 4.2 and 4.4 in \cite{CTZ}; moreover the points $x_j$ there may not be in the same horn $H$.)  More generally, by using the argument in Section 2 of \cite{H95} and Kotschwar's backwards uniqueness theorem for the Ricci flow  \cite{Ko} we show a neck strengthening lemma (see Lemma \ref{lem 2.5}), which extends and improves Lemma
4.3.5 in \cite{BBB+}.  Lemma \ref{lem 2.5} is used in the proof of Propositions \ref{prop 3.9}, \ref{prop 3.10}, \ref{prop 4.5} and \ref{prop 4.10}. We also refer the reader to p. 101 of  Bessi\`{e}res et al. \cite{BBB+} for some details for the application of their Lemma
4.3.5 in the construction of Ricci flow with surgery.

As mentioned above, our main tool is Ricci flow with surgery on compact orbifolds of dimension $n\geq 12$ with isolated singularities and with positive isotropic curvature. One of the key ingredients of the construction of the Ricci flow with surgery is Brendle's curvature pinching estimates, i.e. Corollary 1.3 in \cite{B19}, which also holds true in the orbifold case. As in \cite{H97}, \cite{BBB+}, \cite{BBM}, \cite{Hu13}  and \cite{Hu15}, we do surgery when the maximal scalar curvature $R_{\max}$ reaches certain  threshold before the curvature blows up. (Of course in the compact case one can also do surgery exactly when the curvature blows up as in \cite{P2}, \cite{B19} and \cite{CTZ}, but it seems that the surgery procedure  in \cite{H97} and \cite{BBB+} is more easily adapted to the noncompact case, so we prefer to adopt the surgery procedure  in \cite{H97} and \cite{BBB+} for the future.)

The proof of our theorems  depends crucially on Brendle's work \cite{B19}, and I was also inspired
by Chen-Tang-Zhu's work \cite{CTZ},  the works of Besson and his
collaborators   \cite{BBB+} \cite{BBM}, and the work of Kleiner-Lott \cite{KL}, etc. As mentioned above, we extensively use results and techniques from topology, for example, Cerf-Palais theorem, isotopy extension theorem and covering space technique.  I also make great efforts to add some necessary details involved
in the construction of Ricci flow with surgery on compact orbifolds, which were
neglected in the previous related works. For  examples, see the proof of Propositions  \ref{prop 3.10}, \ref{prop 4.5} and \ref{gluing1}, and Lemmas \ref{gluing1.5} and \ref{lem 2.5}.

 After this manuscript was submitted, we continue to investigate to improve its main results. Note that  the Remark on p. 312 in Bredon \cite{B72} gives an equivariant ambient tubular neighborhood theorem, and more generally one can also derive an orbifold version of the ambient tubular neighborhood theorem; similarly there is an equivariant/orbifold version of the ambient collar theorem. (Using the latter  we can get an orbifold version of Theorem 2.3 in Chapter 8 of \cite{Hi}.) We can  use the ambient isotopy uniqueness of closed tubular neighborhoods of compact suborbifolds  to prove  fine properties of orbifold ancient $\kappa$-solutions and standard solutions when they contain caps with ends of type $\mathbb{S}^{n-1}/\Gamma \times (-1,1)$ with $\Gamma$ nontrivial, extending  Propositions \ref{prop 3.5},   \ref{prop 3.6} and  \ref{prop 3.10},   and strengthen the  definition of $\varepsilon$-caps in Section 2 in these cases.  Using this, we can show that in fact the weak connected sum in Theorem \ref{thm 1.3} can be strengthened to be connected sum,  and furthermore, we can improve Theorem \ref{thm 1.1} to a more precise classification result; see \cite{Hu23}. Moreover, it is possible to extend the main results in this paper and \cite{Hu23} to the case of open manifolds of dimension $n\geq 12$ with uniformly positive isotropic curvature and bounded geometry; see \cite{Hu23b}.  Recently, Zhengnan Chen \cite{Ch} extended Brendle's curvature pinching estimates to the case of $n\geq 9$. So combined with Chen's estimates, our arguments imply that all these results also hold true for $9 \leq n \leq 11$.

  In Section 2 we classify  odd ($\geq 3$)-dimensional spherical orbifolds with nonempty and isolated singularity set (for the even ($\geq 4$)-dimensional case see \cite{CTZ}), define $\varepsilon$-caps,
and describe the canonical neighborhood structure of orbifold ancient
 $\kappa$-solutions and orbifold standard solutions. In Section 3 we prove the so called bounded curvature at bounded distance property, choose the cutoff
parameters for surgical Ricci flow under the canonical neighborhood assumption,  and construct  an $(r, \delta)$-surgical solution to the Ricci flow starting with a compact, connected
Riemannian orbifold of dimension $n\geq 12$ with isolated singularities and with  positive isotropic curvature. In Section 4, we prove
Theorem 1.3 and Theorem 1.1 using the construction in Section 3, and  give a proof of Corollary 1.2 using Theorem 1.1.   In the appendix, we  introduce various notions on necks and (topological) caps, construct orbifold standard solutions, and establish results on the composition of $\varepsilon$-isometries and the gluing and strengthening of $\varepsilon$-necks, which are used in the previous sections.

\vspace *{0.2cm}

\noindent {\bf Acknowledgements}.  I'm very grateful  to  Prof. Simon Brendle for answering my questions on his paper \cite{B19}.  I would also like to thank
Prof. Huai-Dong Cao for helpful communications on Hamilton's definition of the normal neck in \cite{H97}, and Prof. Xi-Ping Zhu for his comments.  I was partially supported by NSFC (12271040) and Beijing Natural Science Foundation (Z190003).

\section{Orbifold ancient $\kappa$-solutions and standard solutions}

Before we  investigate  the structure of orbifold
 ancient $\kappa$-solutions,
we need several lemmas.

\begin{lem} \label{lem 3.1}
Let $k\geq 2$. There does not exist any  isometric involution of $\mathbb{R}P^k$ with fixed point set nonempty and isolated.
\end{lem}

\noindent {\bf Proof}.\ \ Let $\sigma$ be an  isometric involution of $\mathbb{R}P^k$.  We can lift $\sigma$ to an isometry $\tilde{\sigma}$ of $\mathbb{S}^k$ with $\pi\circ \tilde{\sigma}=\sigma \circ \pi$, where $\pi: \mathbb{S}^k \rightarrow \mathbb{R}P^k$ is the natural projection. Since $\sigma$ is an involution, we see that for any $x \in \mathbb{S}^k$, either $\tilde{\sigma}^2(x)=x$ or $\tilde{\sigma}^2(x)=-x$. As $\tilde{\sigma}^2\in O(k+1)$, we see that  either $\tilde{\sigma}^2=I$ or $\tilde{\sigma}^2=-I$, where $I\in O(k+1)$ is the identity matrix. (This also follows from the unique lifting property.) If   $\tilde{\sigma}^2=I$, then  $\tilde{\sigma}$ is an isometric involution of $\mathbb{S}^k$, and the eigenvalues of $\tilde{\sigma} \in O(k+1)$ are either 1 or $-1$. Since $k\geq 2$, at least one eigenvalue of $\tilde{\sigma}$ has multiplicity $\geq 2$, which implies that the fixed point set of $\sigma$  can not be isolated. If $\tilde{\sigma}^2=-I$ (this is possible only when $k$ is odd), one easily sees that $\sigma$ has no fixed points.
\hfill{$\Box$}

\vspace*{0.2cm}

The following lemma  is an analogue of  Lemma 5.2 in \cite{CTZ}.

\begin{lem} \label{lem 3.2}
Let $k\geq 2$ be an integer,  $G\subset O(2k)$  be a finite subgroup such that each nontrivial element in $G$ has at most one eigenvalue equal to 1 and  there is at least one element which has exactly one eigenvalue equal to 1. Then $G\cong \mathbb{Z}_2$.
\end{lem}

\noindent {\bf Proof}.\ \  Compare  Case 2 in the proof of Theorem 3.4 in \cite{CTZ} and the proof of Lemma 5.1 in \cite{CTZ}.  Let $\sigma \in G$ be an element which has exactly one eigenvalue equal to 1. We observe that the order of $\sigma$ is 2. The reason is as follows. Let $\sigma(v)=v$ with $v\in \mathbb{R}^{2k}$ and $|v|=1$. Consider  the orthogonal complement $W$ of $v$ in $\mathbb{R}^{2k}$ and the unit sphere $\mathbb{S}^{2k-2}\subset W$. Note that the subgroup $\langle \sigma \rangle$ generated by $\sigma$ acts freely and nontrivially on $\mathbb{S}^{2k-2}$. But the only nontrivial group that acts freely on $\mathbb{S}^{2k-2}$ is $\mathbb{Z}_2$.  So $\langle \sigma \rangle \cong \mathbb{Z}_2$, which acts on  $\mathbb{S}^{2k-2}$ antipodally.

Now we see that $\sigma$ is the only element in  $G$ which  has exactly one eigenvalue equal to 1. Otherwise suppose there is another element  $\sigma' \in G$ which also has exactly one eigenvalue equal to 1. Note that $\sigma'$ also has order 2. Let $E$ and $E'$ be the $(-1)$-eigenspaces of $\sigma$ and $\sigma'$ respectively. Then $\dim E=\dim E'=2k-1$. The element $\sigma \sigma'$ acts as the identity on $E\cap E'$, but  $\dim (E\cap E')\geq 2k-1+2k-1-2k=2k-2\geq 2$. A contradiction to the assumption  that each nontrivial element in $G$ has at most one eigenvalue equal to 1. (This argument is similar to that in Case 1 in the proof of Lemma 5.1 in \cite{CTZ}.)

We claim that $G\cap SO(2k)$ contains only one element, that is the identity. Suppose otherwise. Let $I\neq \gamma \in G\cap SO(2k)$, where $I\in O(2k)$ is the identity matrix. Note that since $\sigma\in O(2k)$ has  one eigenvalue equal to 1 and $2k-1$ eigenvalues equal to $-1$, $\sigma: \mathbb{S}^{2k-1}\rightarrow \mathbb{S}^{2k-1} $ must be orientation-reversing, and $\gamma \sigma: \mathbb{S}^{2k-1}\rightarrow \mathbb{S}^{2k-1}$ must also be orientation-reversing. By Lefschetz fixed point theorem we know that any orientation-reversing homeomorphism of $\mathbb{S}^{2k-1}$ has (at least) a fixed point, so  $\gamma \sigma: \mathbb{S}^{2k-1}\rightarrow \mathbb{S}^{2k-1}$ has a fixed point. (Of course, since $\gamma \sigma \in O(2k) \setminus SO(2k)$,  this also follows from the fact that any element in a real orthogonal group is orthogonally similar to certain  block diagonal  matrix  with $\pm 1 $'s and/or some 2-dimensional rotation matrices  on the diagonal.) It follows that $\gamma \sigma$ has exactly one eigenvalue equal to 1. A contradiction to what we have proved in the last paragraph.  It follows that $G\cong \mathbb{Z}_2$.
\hfill{$\Box$}

\vspace*{0.2cm}

The following lemma is a simple application of the Cerf-Palais disk theorem.

\begin{lem} \label{lem 3.3}  Consider two orbifolds both diffeomorphic to $\mathbb{S}^n//(x,\pm x')$ ($n\geq 2$),  and write them as the form $D_1 \cup  \mathcal{O}_1$ and $D_2 \cup \mathcal{O}_2$, where both  $D_1$ and $D_2$ are diffeomorphic to $D^n$,  and $D_1 \cap \mathcal{O}_1=\partial D_1=\partial \mathcal{O}_1$, $D_2 \cap \mathcal{O}_2=\partial D_2=\partial \mathcal{O}_2$.     Assume   that         $f: D_1 \cup \mathcal{O}_1 \rightarrow D_2 \cup   \mathcal{O}_2 $ is a  diffeomorphism, and $h: \partial D_1 \rightarrow \partial D_2$ is a diffeomorphism which extends to a diffeomorphism $\tilde{h}: D_1 \rightarrow D_2$. Moreover assume that $f|_{D_1}$ and  $\tilde{h}$ are either both orientation-preserving or both orientation-reversing when $n$ is even. Then $h: \partial \mathcal{O}_1 \rightarrow  \partial \mathcal{O}_2$ extends to a diffeomorphism $\tilde{h}': \mathcal{O}_1 \rightarrow \mathcal{O}_2$.
\end{lem}
\noindent {\bf Proof}.\ \ By assumption $f|_{D_1}:   D_1 \rightarrow f(D_1)\subset D_2 \cup \mathcal{O}_2$ and $\tilde{h}:  D_1 \rightarrow D_2 \subset  D_2 \cup \mathcal{O}_2 $  are two embeddings of  $D_1$ in   $ D_2 \cup \mathcal{O}_2$ (either both orientation-preserving or both orientation-reversing if $n$ is even).  Note that $\mathbb{S}^n//(x,\pm x')$ ($n\geq 2$) is orientable if and only if $n$ is even, and if and only if the complement of the singularity set is orientable. Applying the Cerf-Palais disk theorem (see  for example Theorem 3.1 in Chapter 8 of \cite{Hi} or Theorem 2.1 on p. 197 of \cite{M07}) to the manifold $ D_2 \cup \mathcal{O}_2 \setminus \{\text{singularities}\}$, we see that there is a diffeomorphism $\varphi:  D_2 \cup \mathcal{O}_2 \rightarrow  D_2 \cup \mathcal{O}_2 $  such that $\varphi \circ f|_{D_1}=\tilde{h}$.  Then  $\varphi \circ f|_{\mathcal{O}_1}: \mathcal{O}_1 \rightarrow \mathcal{O}_2$   extends $h:  \partial \mathcal{O}_1 \rightarrow  \partial \mathcal{O}_2$.
\hfill{$\Box$}

\vspace *{0.2cm}

Let $(\mathcal{O}, g)$ be a connected Riemannian orbifold. Following \cite{KL}, we call a path in $\mathcal{O}$ as defined in Section 2.2.1 of \cite{BMP} a special curve. We emphasize that for a special curve $\gamma: I\rightarrow \mathcal{O}$, to each $t\in I$ (here $I\subset \mathbb{R}$ is an interval) such that $|\gamma|(t)$ is a singular point, there is associated a smooth local lift of $\gamma$ around $t$ which is given as  part of the data. (See also Definition 3.1 in \cite{BB12}  or Definition 2.1 in \cite{BB13}.)  Given a special curve $\gamma$ in $\mathcal{O}$, we can pullback the tangent bundle $T\mathcal{O}$  and the Levi-Civita connection of $(\mathcal{O}, g)$ via $\gamma$ (cf. \cite{CR} and Section 6 of \cite{BB13}), and define parallel transport along $\gamma$ as in the case of  Riemannian manifolds. (One can also define  piecewise smooth special curves and parallel transport along a piecewise smooth special curve.) We also define  homotopy of special curves as in Section 2.2.1 of \cite{BMP}  and Section 2.2 of \cite{KL}. (We can also define  homotopy of piecewise smooth special curves by using piecewise smooth homotopies in  local models; for the manifold case see for example Appendix 7 in \cite{KN}.)

Now let $p\in |\mathcal{O}|$ be a regular point. For a special loop   $\gamma: [0,1] \rightarrow \mathcal{O}$ based at $p$, let $H_\gamma$ be the holonomy around $\gamma$ as defined in Section 2 of \cite{KL}. In \cite{KL} $H_\gamma$ is given by a horizontal lift of $\gamma$, that is, a lift $\tilde{\gamma}:[0,1]\rightarrow F\mathcal{O}$ of $\gamma$ to the orthonormal frame bundle $F\mathcal{O}$ of $\mathcal{O}$ such that the tangent vector $\tilde{\gamma}'(t) \in T_{\tilde{\gamma}(t)}^HF\mathcal{O}$ for each $t\in [0,1]$, where $T^HF\mathcal{O}$ is the horizontal distribution  on $F\mathcal{O}$ arising from the Levi-Civita connection of $(\mathcal{O}, g)$. (Equivalently, one can also define $H_\gamma$ using parallel transport along $\gamma$ as in the case of  Riemannian manifolds.) Let $\text{Hol}^0(\mathcal{O}, g,p)$ be the subgroup of the special orthogonal group $SO(T_p\mathcal{O})$ generated by those $H_\gamma$'s with $\gamma$ running through the set of null-homotopic special loops based at $p$, and call it  the restricted holonomy group of $(\mathcal{O},g)$ at $p$.

We can extend the proofs of the holonomy theorem of Ambrose-Singer and Nijenhuis in \cite{Ba} and \cite{KN} (see also \cite{Bs} and \cite{Sa}) to the orbifold case.
Furthermore, we can extend the proofs in \cite{Si} (see also \cite{Bs} and \cite{Sa}) and \cite{BCO} of Berger's holonomy classification theorem  to the orbifold case. So both theorems hold in the orbifold case.

\vspace *{0.2cm}

As in \cite{B10}, let $\mathcal{C}_B(\mathbb{R}^n)$ be the vector space of algebraic curvature tensors on $\mathbb{R}^n$. Following the notation in \cite{To}, let $C_{\text{PIC}}$, $C_{\text{PIC}1}$ and $C_{\text{PIC}2}$ be the subsets of  $\mathcal{C}_B(\mathbb{R}^n)$ which consist of algebraic curvature tensors of nonnegative isotropic curvature (weakly PIC, see item (i) of Definition 1.1 in \cite{B19}), of weakly PIC1 (see item (ii) of Definition 1.1 in \cite{B19}), and of weakly PIC2 (see item (iii) of Definition 1.1 in \cite{B19}) respectively. Recall (\cite{BS}\cite{B10}) that the curvature tensor of an $n$-dimensional Riemannian orbifold $(\mathcal{O},g)$ lies in the cone $C_{\text{PIC}1}$ if and only if at each point the Riemannian orbifold $\mathcal{O}\times \mathbb{R}$ with the product metric has nonnegative isotropic curvature, and the curvature tensor of  $(\mathcal{O},g)$ lies in the cone $C_{\text{PIC}2}$  if and only if at each point the Riemannian orbifold $\mathcal{O}\times \mathbb{R}^2$ with the product metric has nonnegative isotropic curvature.  (Here,  at a singular point, we should look at the curvature tensor at the corresponding point in a local model.) Note that the curvature tensor of  $(\mathcal{O},g)$ lies in the cone $C_{\text{PIC}2}$  if and only if $(\mathcal{O},g)$ has nonnegative complex sectional curvature; a proof (attributed  to  N. Wallach) of this equivalence  appeared in  a 2007 arXiv preprint of L. Ni and J. Wolfson, see also for example Proposition 7.18 in \cite{B10}. When (at each point) the curvature tensor of  $(\mathcal{O},g)$ lies in the interior of the cone $C_{\text{PIC}2}$, we also say that  $(\mathcal{O},g)$ is strictly PIC2 (in this case, the inequality in item (iii) of Definition 1.1 in \cite{B19} is strict; see also Section 2 in \cite{BN2}); this condition is equivalent to that $(\mathcal{O},g)$ has positive complex sectional curvature.  (One can also consult the first paragraph on p. 55 and Chapter 7 of \cite{B10}.)

\vspace *{0.2cm}

Note that Proposition 6.6 in \cite{B19} can be extended to the orbifold case, as all the key ingredients of its proof, including the Berger holonomy theorem,  Proposition 9 in \cite{BS08} (which is an application of the Bony type strong maximum principle), Theorem 16 in \cite{B10a},  and  the de Rham decomposition theorem,   can be extended to the orbifold case. (For an extension of the de Rham decomposition theorem to the orbifold case see Lemma 2.19 in \cite{KL}.)  For another application of the combination of the strong maximum principle in \cite{BS08}  and orbifold holonomy see p. 585 of \cite{NW}.

The following result is a simple application of the orbifold version of Proposition 6.6 in \cite{B19} and the strong maximum principle for systems.

\begin{prop} \label{prop 3.4}
Let $(\mathcal{O},g(t)), t\in [0,T]$, be a complete solution to the Ricci flow. Suppose that $(\mathcal{O},g(t)), t\in [0,T]$, is weakly PIC2 and strictly PIC, and there  is a point $x_0$ such that the Ricci curvature  of $g(0)$ at $x_0$ is positive. Then for any $t>0$ and any $x\in \mathcal{O}$, the curvature tensor of $g(t)$ at $x$ lies in the interior of the cone $C_{\text{PIC}2}$.
\end{prop}
\noindent {\bf Proof}.\ \  Otherwise, there is a space-time point $(p,t')$ with $t'>0$ such that the curvature tensor at $(p,t')$ lies on the boundary of the cone $C_{\text{PIC}2}$. By the orbifold version of Proposition 6.6 in \cite{B19}  for all $0<t\leq t'$ the universal cover of $(\mathcal{O}, g(t))$ splits off a line. In particular, the Ricci curvature of $(\mathcal{O}, g(t))$ cannot be strictly positive for  $0<t\leq t'$.  On the other hand, by assumption the Ricci curvature of $g(0)$ is positive at the point $x_0$, so the Ricci curvature is positive everywhere  for any $t>0$ by  the strong maximum principle for systems (Proposition 12.47 in \cite{C+}); here note that the reaction term in the evolution equation for the  Ricci tensor is weakly positive definite under the weakly PIC2 assumption. A contradiction.   \hfill{$\Box$}

\vspace *{0.2cm}

Let $\kappa, r>0$. Recall (cf. \cite{P1} and Definition 4.1 in \cite{BBM}) that  an  evolving  Riemannian $n$-orbifold   $\{(\mathcal{O}(t), g(t))\}_{t\in I}$  is
$\kappa$-noncollapsed  on the scale $r$ at time
$t_0\in I$ if for any point $x\in \mathcal{O}(t_0)$, whenever $|\text{Rm}|\leq r^{-2} \hspace{2mm}$ on the parabolic ball $P(x, t_0, r, -r^2):=\{(x',t') \hspace*{1mm} | \hspace*{1mm} x'\in B(x,t_0,r), t'\in [t_0-r^2,t_0]\}$ we have $
\mbox{vol}_{g(t_0)} \hspace*{1mm} B(x,t_0, r)\geq \kappa r^n$.

\vspace *{0.2cm}

Let $n\geq 5$ and $\kappa >0$.
 An  orbifold ancient  $\kappa$-solution of dimension $n$ is  a
 complete, nonflat  solution  to the orbifold Ricci flow  with at most isolated singularities and with bounded curvature  which is  defined on some time interval $(-\infty, T]$ and is weakly PIC2 and  $\kappa$-noncollapsed on all scales at any time $t\in (-\infty, T]$. (See \cite{BN1}, \cite{BN2}, \cite{CL}, and \cite{LZ} for some recent works on ancient $\kappa$-solutions on high-dimensional manifolds.)

\vspace *{0.2cm}

Now we start to analyze the structure of orbifold ancient  $\kappa$-solutions which satisfy certain pinching assumption.  Recall (cf. for example \cite{Bs} and \cite{B19}) that for two symmetric $(0,2)$-tensors $A$ and $B$, their Kulkarni-Nomizu product is the $(0,4)$-tensor $A \owedge B$ given by $(A \owedge B)_{ijkl}=A_{ik}B_{jl}+A_{jl}B_{ik}-A_{il}B_{jk}-A_{jk}B_{il}$.   Let $r_1:\mathbb{R}^n \rightarrow \mathbb{R}^n$ be the  reflection w.r.t. the hyperplane $x_1=0$, we also use $r_1$ to denote the restriction of $r_1$ to $\mathbb{S}^{n-1}$.

\begin{prop} \label{prop 3.5} (cf. Corollary 6.7 in \cite{B19}, and Theorem 3.4 in \cite{CTZ}) \ \  Let  $n\geq 5$, and  $(\mathcal{O},g(t))$, $t\in (-\infty,T]$, be an
orbifold ancient  $\kappa$-solution of dimension $n$ (with at most isolated singularities). Suppose that $(\mathcal{O},g(t))$  satisfies $Rm-\theta R \hspace*{1mm} \text{id} \owedge \text{id} \in C_{\text{PIC}}$ for some uniform constant $\theta>0$, and  there is a spacetime point $(x_0,t_0)$ such that the curvature tensor at $(x_0,t_0)$ lies on the boundary of the cone $C_{\text{PIC}2}$. Then for each $t$, $(\mathcal{O},g(t))$
 is isometric to a shrinking Ricci soliton
$\mathbb{S}^{n-1}/\Gamma \times \mathbb{R}$ or $\mathbb{S}^{n-1}/\Gamma
\times_{\mathbb{Z}_2} \mathbb{R}$ for some finite subgroup $\Gamma$  of $O(n)$ acting freely on $\mathbb{S}^{n-1}$. In particular, if $\mathcal{O}$ has no singularities and has exactly one end, it must be diffeomorphic to some $\mathbb{S}^n// \langle\Gamma, \hat{\sigma} \rangle \setminus \bar{B}$, where $\sigma$ is an isometric involution of the spherical space form $\mathbb{S}^{n-1}/\Gamma$ with no fixed points, and $\hat{\sigma}$ and $B$ are defined as in the appendix; if $\mathcal{O}$ has nonempty isolated singularities,  it must be diffeomorphic to $\mathbb{S}^n// (x,\pm x') \setminus \bar{B}$ as defined in the appendix. Furthermore, there exists a positive constant $\varepsilon_4=\varepsilon_4(n)\leq \varepsilon_1$, where $\varepsilon_1$ is the constant in Proposition \ref{gluing1}, with the following property.

\noindent When $\mathcal{O}$ is diffeomorphic to $\mathbb{R}P^n \setminus \bar{B}$ (that is, when $\mathcal{O}$ is diffeomorphic to  $\mathbb{S}^n// \langle\Gamma, \hat{\sigma} \rangle \setminus \bar{B}$ with  $\Gamma$  trivial),  any $2\varepsilon_4$-neck in  $(\mathcal{O}, g(t))$ must be diffeomorphic to $\mathbb{S}^{n-1}\times (0,1)$, and  the central   cross section $\Sigma$ of any $2\varepsilon_4$-neck in  $\mathcal{O}$ bounds a   compact domain $\Omega$, moreover, if $f: \mathbb{S}^{n-1}\rightarrow \Sigma$ is a $2\varepsilon_4$-homothety coming from the $2\varepsilon_4$-neck structure, and $f': \partial (\mathbb{R}P^n \setminus B) \rightarrow \mathbb{S}^{n-1}$ is a homothety, either $f\circ f'$ or $f\circ r_1\circ f'$ extends to a  diffeomorphism  $F: \mathbb{R}P^n \setminus B \rightarrow \Omega$;

\noindent When $\mathcal{O}$ is diffeomorphic to $\mathbb{S}^n// (x,\pm x') \setminus \bar{B}$,  any $2\varepsilon_4$-neck in  $(\mathcal{O}, g(t))$ must be diffeomorphic to $\mathbb{S}^{n-1}\times (0,1)$, the central cross section $\Sigma$ of any $2\varepsilon_4$-neck in  $\mathcal{O}$ bounds a compact domain $\Omega$, moreover, if $f: \mathbb{S}^{n-1}\rightarrow \Sigma$ is a $2\varepsilon_4$-homothety coming from the $2\varepsilon_4$-neck structure, and $f': \partial (\mathbb{S}^n// (x,\pm x') \setminus B) \rightarrow \mathbb{S}^{n-1}$ is a homothety, either $f\circ f'$ or $f\circ r_1\circ f'$ extends to a  diffeomorphism  $F: \mathbb{S}^n// (x,\pm x') \setminus B \rightarrow \Omega$.
\end{prop}

 \noindent {\bf Proof}.\ \   Compare the proof of Corollary 6.7 and Proposition 6.17 in \cite{B19}  and Theorem 3.4 in \cite{CTZ}. We pull back our solution  to the
universal cover and get $(\widetilde{\mathcal{O}},\tilde{g}(t))$.  As noted above the  orbifold version of Proposition 6.6 in \cite{B19} holds. It follows that at any time $t$, $(\widetilde{\mathcal{O}},\tilde{g}(t))$ is
isometric to a product $(X, g_X(t))\times \mathbb{R}$, where $X$ is a smooth manifold. Since by assumption $(X, g_X(t))\times \mathbb{R}$ is nonflat, weakly PIC2  and satisfies $\text{Rm}-\theta R \hspace*{1mm} \text{id} \owedge \text{id} \in C_{\text{PIC}}$ for some uniform constant $\theta>0$,  the ancient solution $(X, g_X(t))$ is weakly PIC2 and uniformly PIC1.  By Theorem 6.4 in \cite{B19},  we see that at any time
$t \leq t_0$, $(X,g_X(t))$ is a round $(n-1)$-sphere. Now we see that $(\mathcal{O},g(t))$ is
a metric quotient of the evolving round cylinder $\mathbb{S}^{n-1} \times \mathbb{R}$ by standard isometries (compare Theorem 13.3.10 in \cite{R}).
 Then it follows from the $\kappa$-noncollapsing assumption on all scales that
$\mathcal{O}$ is noncompact (cf. p.212 in \cite{CZ} and p.52 in \cite{CTZ}). So $\mathcal{O}$ has exactly one or
two ends, moreover, if $\mathcal{O}$ has two ends, it must split metrically as $\mathcal{O}' \times \mathbb{R}$, where $\mathcal{O}'$ is a metric quotient of $\mathbb{S}^{n-1}$; see Corollary 3.3 in \cite{KL}.   As $\mathcal{O}$ has at most isolated singularities, if it has two ends,  it must be isometric to
$\mathbb{S}^{n-1}/\Gamma \times \mathbb{R}$ for some finite subgroup
$\Gamma$ of $O(n)$ acting freely on $\mathbb{S}^{n-1}$, in particular, in this case $\mathcal{O}$ is a smooth manifold.

If $\mathcal{O}$ has exactly one  end, it must be isometric
to $\mathbb{S}^{n-1}/\Gamma \times_{\mathbb{Z}_2} \mathbb{R}$ for some finite
subgroup $\Gamma$ of $O(n)$ acting freely on $\mathbb{S}^{n-1}$. The reason  is as
follows.  We can write
$\mathcal{O}=\mathbb{S}^{n-1} \times \mathbb{R}//\widetilde{\Gamma}$ for a subgroup
$\widetilde{\Gamma}$ of isometries of the round cylinder $\mathbb{S}^{n-1} \times
\mathbb{R}$. Since $\mathcal{O}$ has one end, we can write $\widetilde{\Gamma}=\Gamma \cup \Gamma^1$, where the second
component of $\Gamma$ (resp. of  $\Gamma^1$) acts on $\mathbb{R}$ as the
identity (resp.  a reflection). Since $\mathcal{O}$ has  only one end, $\Gamma^1 \neq \emptyset$. Pick $\tilde{\sigma} \in \Gamma^1$.
Then $\tilde{\sigma}^2 \in \Gamma$, and $\tilde{\sigma} \Gamma=\Gamma^1$. It follows that $\tilde{\sigma}$ induces an involution, denoted by $\bar{\sigma}$, acting isometrically on $\mathbb{S}^{n-1}/\Gamma \times \mathbb{R}$.   Now we see that
$\mathcal{O}=(\mathbb{S}^{n-1}/\Gamma \times \mathbb{R})//\langle
\bar{\sigma} \rangle$, which is  of the form $\mathbb{S}^{n-1}/\Gamma\times_{\mathbb{Z}_2} \mathbb{R}$ by definition.

In the case that $\mathcal{O}$ has exactly one  end, we see that the action of $\tilde{\Gamma}$ on $\mathbb{S}^{n-1} \times \mathbb{R}$ will leave exactly one cross section, say $\mathbb{S}^{n-1} \times \{0\}$, invariant. (Compare the beginning of the third paragraph on p.53 in \cite{CTZ}.) Let $\pi: \mathbb{S}^{n-1} \times \mathbb{R} \rightarrow \mathcal{O}=\mathbb{S}^{n-1} \times \mathbb{R}//\tilde{\Gamma}$ be the natural projection. We consider further two subcases.

\vspace *{0.2cm}

\noindent Subcase 1:  $\mathcal{O}$ is a smooth manifold. Then the first component, denoted by $\sigma$,  of the isometric involution $\bar{\sigma}$ above must act on $\mathbb{S}^{n-1}/\Gamma$ without any fixed points. This manifold is just $C_\Gamma^{\sigma}$ as defined in the appendix. Note that in this subcase, if $n$ is odd,  the subgroup $\Gamma$ must be trivial, and the first component of the isometric involution $\bar{\sigma}$ above must act on $\mathbb{S}^{n-1}$  antipodally; because as noted above if $n$ is odd $\Gamma$ must be trivial or $\mathbb{Z}_2$, but if $\Gamma\cong \mathbb{Z}_2$, no such manifold  exists, since by the Lefschetz fixed point theorem, any continuous self-map  of an even dimensional real projective space must have a fixed point. (Here, the use of Lefschetz fixed point theorem can also be avoided  by lifting the map to a self-map of an even dimensional sphere.)

In this subcase we  see that $\mathcal{O}$  must be diffeomorphic to some $\mathbb{S}^n// \langle\Gamma, \hat{\sigma} \rangle \setminus \bar{B}$.  (Compare p. 49 in \cite{CTZ} and the appendix.)

 Now we further assume that $\Gamma$ is trivial, then $\mathcal{O}$ is  the quotient of $\mathbb{S}^{n-1} \times \mathbb{R}$ by the $\mathbb{Z}_2$-action generated by $((x_1,x_2, \cdot\cdot\cdot, x_n), s)\mapsto ((-x_1, -x_2, \cdot\cdot\cdot, -x_n), -s)$, and is diffeomorphic to $\mathbb{R}P^n \setminus \bar{B}$. Note that $\pi(\mathbb{S}^{n-1} \times \{0\})$ is  $\mathbb{R}P^{n-1}$.
Fix any $t$ and write $g=g(t)$.

\vspace *{0.2cm}

\noindent {\bf Claim} 1.  Let $0<\varepsilon_4=\varepsilon_4(n)\leq \varepsilon_1$ be sufficiently small.  Any $2\varepsilon_4$-neck in  $(\mathcal{O},g)$ must be diffeomorphic to $\mathbb{S}^{n-1}\times (0,1)$, and the central cross section $\Sigma$ of any $2\varepsilon_4$-neck   bounds a compact, connected, (smooth) submanifold in  $\mathcal{O}$ which contains $\pi(\mathbb{S}^{n-1} \times \{0\})$.

\vspace *{0.2cm}

\noindent {\bf Proof of Claim 1}.  Given  any $2\varepsilon_4$-neck $N$ in  $(\mathcal{O},g)$, note  that the diameter of $\pi (\mathbb{S}^{n-1} \times \{0\})$ is much smaller than that of the neck $N$ since $\varepsilon_4$ is very small. It follows that at least four ninths  of the neck $N$ is contained in a standard cylindrical neck in $\mathcal{O}$.  By Proposition  \ref{gluing1} any cross section of the neck $N$ is diffeomorphic to $\mathbb{S}^{n-1}$, so the neck $N$  is diffeomorphic to $\mathbb{S}^{n-1} \times (0, 1)$.   It is  well-known that  for $n\geq 5$ (in fact for all $n\geq 3$), $\mathbb{R}P^{n-1}$ cannot be embedded in $\mathbb{R}^n$, in particular,   $\mathbb{R}P^{n-1}$ cannot be embedded in $\mathbb{S}^{n-1} \times (0, 1)$. So the intersection of the central cross section $\Sigma$ of the neck $N$ with $\pi(\mathbb{S}^{n-1} \times \{0\})$ is empty. It follows that  the middle half of the neck $N$ is contained in a standard cylindrical neck, say $N'$, in $\mathcal{O}$. By Proposition  \ref{gluing1} again we see that the conclusion of Claim 1 holds. (In fact, by Proposition  \ref{gluing1} we see  that the compact domain $\Omega$ bounded by $\Sigma$ is diffeomorphic to $\mathbb{R}P^n \setminus B$, as  any cross section of a standard cylindrical neck in $(\mathcal{O}, g)$ bounds such a domain.) By the way, once we know that the neck $N$ is   diffeomorphic to $\mathbb{S}^{n-1} \times (0, 1)$ and  the middle half of the neck $N$ is contained in a standard cylindrical neck, we can also proceed as follows.  Recall the following result: Let $n\geq 5$. If $S$ is a smoothly embedded $(n-1)$-sphere in  $D^n\setminus \partial D^n$, then $S$ bounds a smoothly embedded $n$-disk $D'$; moreover, the complement of $\text{Int}\hspace*{0.5mm} D'$ in  $D^n$ is diffeomorphic to $\mathbb{S}^{n-1} \times [0,1]$.  (The former fact follows from the solution of
Schoenflies conjecture in dimensions $n\geq 5$, and the latter fact follows from the former fact and the Cerf-Palais disk theorem.) Then arguing as in the proof of Corollary A.2.2 and Lemma 3.2.2 in \cite{BBB+}, we see that if $\varepsilon_4$ is sufficiently small, $\Sigma$ and certain cross section of $N'$ cobounds a compact domain which is diffeomorphic to $\mathbb{S}^{n-1}\times [0,1]$,  and  Claim 1 follows.
\hfill{$\Box$

\vspace *{0.2cm}

 With the help of Claim 1, the existence of the diffeomorphism $F$ for the $\mathbb{R}P^n \setminus \bar{B}$ case in the conclusion of this proposition  follows from a surgery argument as in  the proof of Proposition 6.17 in \cite{B19}  and Subcase 2 below (in fact it is slightly easier compared to Subcase 2 below), and is omitted.

\vspace *{0.2cm}

\noindent Subcase 2:  $\mathcal{O}$  has at least one isolated orbifold singularity.  If $n$ is even,  by restricting the action of  $\tilde{\Gamma}$  to $\mathbb{S}^{n-1} \times \{0\}$ and  using Lemma \ref{lem 3.2} and the assumption that $\mathcal{O}$ has at most isolated singularities, we see that $\tilde{\Gamma}\cong \mathbb{Z}_2$ and $\mathcal{O}$ is the quotient of $\mathbb{S}^{n-1} \times \mathbb{R}$ by the $\mathbb{Z}_2$-action generated by $((x_1,x_2, \cdot\cdot\cdot, x_n), s)\mapsto ((x_1, -x_2, \cdot\cdot\cdot, -x_n), -s)$; one can also show this by arguing as in Case 2 in the proof of Theorem 3.4 in \cite{CTZ}.  If $n$ is odd, the subgroup $\Gamma$ above must be trivial or $\mathbb{Z}_2$. If  $\Gamma$  is trivial, clearly  $\mathcal{O}$ is the quotient of $\mathbb{S}^{n-1} \times \mathbb{R}$ by the $\mathbb{Z}_2$-action generated by $((x_1,x_2, \cdot\cdot\cdot, x_n), s)\mapsto ((x_1, -x_2, \cdot\cdot\cdot, -x_n), -s)$.
If $\Gamma \cong \mathbb{Z}_2$, by Lemma \ref{lem 3.1} no such $\mathcal{O}$  exists.

 From the last paragraph we see that in Subcase 2, $\mathcal{O}$ is the quotient of $\mathbb{S}^{n-1} \times \mathbb{R}$ by the $\mathbb{Z}_2$-action generated by $((x_1,x_2, \cdot\cdot\cdot, x_n), s)\mapsto ((x_1, -x_2, \cdot\cdot\cdot, -x_n), -s)$. It has exactly two singular points, and is diffeomorphic to $\mathbb{S}^n// (x,\pm x') \setminus \bar{B}$ as defined in the appendix, which can be shown by using certain $O(n)\times \mathbb{Z}_2$-invariant metric on the (topological) $n$-sphere constructed, for example, in the proof of Proposition 3.1 in \cite{Hu23}. Note that $\pi(\mathbb{S}^{n-1} \times \{0\})$ is  the quotient of $\mathbb{S}^{n-1}$ by the $\mathbb{Z}_2$-action generated by $(x_1,x_2, \cdot\cdot\cdot, x_n)\mapsto (x_1, -x_2, \cdot\cdot\cdot, -x_n)$, and it contains the two singular points of $\mathcal{O}$.  Fix any $t$ and write $g=g(t)$.

\vspace *{0.2cm}

\noindent {\bf Claim} 2.   Let $0<\varepsilon_4=\varepsilon_4(n) \leq \varepsilon_1$ be sufficiently small.  Any $2\varepsilon_4$-neck in  $(\mathcal{O},g)$ must be diffeomorphic to $\mathbb{S}^{n-1}\times (0,1)$, and the central cross section $\Sigma$ of any $2\varepsilon_4$-neck bounds a compact, connected, (smooth) suborbifold in  $\mathcal{O}$ which contains  $\pi(\mathbb{S}^{n-1} \times \{0\})$.

\vspace *{0.2cm}

\noindent {\bf Proof of Claim 2}.  Given   any $2\varepsilon_4$-neck $N$ in  $(\mathcal{O},g)$, as noted above, the two singular points of $\mathcal{O}$ are contained in $\pi (\mathbb{S}^{n-1} \times \{0\})$, whose diameter is much smaller than that of the neck $N$ if $\varepsilon_4$ is sufficiently small. Of course none of the two singular points is contained in $N$.  So the intersection of the central cross section $\Sigma$ of the neck $N$ with $\pi(\mathbb{S}^{n-1} \times \{0\})$ is empty.  It follows that  the middle half of the neck $N$ is contained in a standard cylindrical neck, say $N'$, in $\mathcal{O}$.    Now Claim 2 follows
from Proposition  \ref{gluing1}. \hfill{$\Box$

 \vspace*{0.2cm}

Let $N$ be a  $2\varepsilon_4$-neck in  $(\mathcal{O},g)$ (in Subcase 2), by Claim 2, it must be diffeomorphic to $\mathbb{S}^{n-1}\times (0,1)$. Let  $\psi: \mathbb{S}^{n-1}\times (-(2\varepsilon_4)^{-1}, (2\varepsilon_4)^{-1}) \rightarrow N$ be a diffeomorphism giving the $2\varepsilon_4$-neck structure of $N$, and   $\Sigma=\psi(\mathbb{S}^{n-1}\times \{0\})$.  By Claim 2 again, $\Sigma$ bounds a
compact domain, say $\Omega$,  in  $\mathcal{O}$  such that $\Omega$ contains  $\pi(\mathbb{S}^{n-1} \times \{0\})$. Let   $f: \mathbb{S}^{n-1}\rightarrow \Sigma$ be a $2\varepsilon_4$-homothety coming from the $2\varepsilon_4$-neck structure, that is,   $f$ is the same as   $\psi|_{\mathbb{S}^{n-1} \times \{0\}}: \mathbb{S}^{n-1} \times \{0\} \rightarrow \Sigma $ after the canonical  identification of $\mathbb{S}^{n-1} \times \{0\}$ with $\mathbb{S}^{n-1}$, and $f': \partial (\mathbb{S}^n// (x,\pm x') \setminus B) \rightarrow \mathbb{S}^{n-1}$ be a homothety. We want to find a diffeomorphism  $F: \mathbb{S}^n// (x,\pm x') \setminus B \rightarrow \Omega$ with the desired boundary value as in the conclusion of this proposition. For this aim, as in the proof of Proposition 6.17 in \cite{B19},  we will do metric surgery on $\mathcal{O}$ along the $2\varepsilon_4$-neck $N$ to get a new Riemannian orbifold $(\mathcal{O}', g')$.
Recall that
\begin{equation*}
|Q\psi^*g-g_{\text{cyl}}|_{C^{[(2\varepsilon_4)^{-1}]}, g_{\text{cyl}}}< 2\varepsilon_4,
\end{equation*}
where $Q$ is a positive constant and $g_{\text{cyl}}$ is the standard  metric on the cylinder $\mathbb{S}^{n-1}\times (-(2\varepsilon_4)^{-1}, (2\varepsilon_4)^{-1})$ with scalar curvature 1  and   length ${\varepsilon_4}^{-1}$ for the interval $(-(2\varepsilon_4)^{-1}, (2\varepsilon_4)^{-1})$. Suppose that the image of $\mathbb{S}^{n-1} \times (-\frac{1}{2\varepsilon_4}, 0]$ under the map $\psi$ is contained in  $\Omega$. Following \cite{B19}, \cite{H97}, and \cite{KL08}, let $\varphi: (0, b_0]\rightarrow \mathbb{R}$ be a smooth  function as chosen in the appendix, where $b_0> \frac{1}{10}$; in particular,
  $\varphi(z)= ce^{-1/z}$ for $z \in (0, \frac{1}{10}]$, where  $c < \frac{1}{\sqrt{(n-1)(n-2)}}$ is a small positive constant.  (We require that $\frac{1}{\varepsilon_4} \gg b_0$.)
  We also choose a smooth cutoff function $\chi: (\frac{1}{20}, b_0]\rightarrow [0,1]$  with $\chi(z)=1$ for $\frac{1}{20} < z \leq \frac{1}{18}$ and $\chi(z)=0$ for $\frac{1}{12} \leq z \leq b_0$. For convenience, below for $x\in \psi(\mathbb{S}^{n-1} \times (0, b_0])$ we will also denote  $z(x)= p_2\circ \psi^{-1}(x)$ by $z$,  where $p_2: \mathbb{S}^{n-1}\times (-(2\varepsilon_4)^{-1}, (2\varepsilon_4)^{-1}) \rightarrow (-(2\varepsilon_4)^{-1}, (2\varepsilon_4)^{-1})$ is the projection onto the second factor. Let
 \begin{equation*}
 g'(x)= \begin{cases}
    g(x), & {x\in \Omega;}  \\
    e^{-2\varphi(z)}g(x), & {z\in (0,\frac{1}{20}];} \\
    e^{-2\varphi(z)}(\chi(z)g(x)+(1-\chi(z))Q^{-1}(\psi_*g_{\text{cyl}})(x)), & {z\in (\frac{1}{20},  b_0],}
     \end{cases}
\end{equation*}
which gives a new metric on the orbifold $\Omega \cup \psi(\mathbb{S}^{n-1} \times (0,  b_0])$.
  If $\varepsilon_4$ is sufficiently small,  the new metric $g'$  will remain to be weakly PIC2 and strictly PIC. Moreover, as in Section 72 of \cite{KL08} and the appendix,  $\varphi$ is chosen such that   the metric $g'$  restricted to  the portion  $\psi(\mathbb{S}^{n-1} \times [b', b_0])$ has constant positive sectional curvature, where $\frac{1}{10} < b' < b_0$ is as chosen in the appendix.
 We extend the metric $g'$ by smoothly attaching a closed metric $n$-ball of constant positive sectional curvature  to $(\Omega \cup \psi(\mathbb{S}^{n-1} \times (0, b_0]), g')$ along the boundary such that the extended metric is weakly PIC2 and strictly PIC; moreover,  the Ricci curvature (even sectional curvature) of the extended metric is positive at (at least) one point. We denote  the resulting closed Riemannian orbifold by $(\mathcal{O}', g')$.  For more details see also the proof of Proposition 8.2 in \cite{B19} (but note that in the definition of the function $\varphi$ there, one should add a small positive factor (depending on the dimension) in front of  $e^{-\frac{1}{z}}$), Section 72 of \cite{KL08}, Chapter 7 of  \cite{BBB+}, and Chapter 13 of \cite{MT}. By the way, it seems that the formula for the function $f$ near the point $z=4$ displayed on p. 424 of \cite{CaZ} and p. 228 of \cite{CZ} (and p. 1221 of \cite{Hu13}  and p. 343 of \cite{Z11}) is not quite true; cf. Lemma 2.10 in \cite{CK}  and Chapter 1 in \cite{Pe}.

 We can write $\mathcal{O}'=D_1\cup \Omega$, where $D_1$ is diffeomorphic to $D^n$, and
 $D_1\cap \Omega=\partial D_1=\partial \Omega =\Sigma$.    By  inspecting the surgery procedure we see that $f$ extends to a diffeomorphism   from $ D^n$ to $D_1$, which factors through an $O(n)$-equivariant diffeomorphism from $ D^n$ to a rotationally symmetric surgery cap constructed in the appendix.  Then we run the orbifold Ricci flow with initial data $ (\mathcal{O}', g')$.  By Proposition \ref{prop 3.4}  this orbifold Ricci flow solution  is strictly PIC2 for any $t>0$. By the differentiable sphere theorem for compact orbifolds which are strictly PIC2 we see that $ (\mathcal{O}', g')$ is diffeomorphic to a spherical orbifold with  two singular points both with local group $\mathbb{Z}_2$, which must be of the form $\mathbb{S}^n// (x,\pm x')$ by using Lemma \ref{lem 3.2} above, Lemma 5.2 in \cite{CTZ} and
    Corollary 2.4  in Chapter VI of \cite{B72}. Note that one can extend the differentiable sphere theorem for compact manifolds which is strictly PIC2 in  \cite{BS} to  the orbifold case either by using the estimates in \cite{BS}, Perelman's
noncollapsing theorem \cite{P1} (which can be easily extended to the orbifold case) and the  orbifold Ricci flow compactness theorem \cite{KL}\cite{L}, or
by first running the (compact) orbifold Ricci flow  which is strictly PIC2 up to certain time such that the sectional curvature is globally strictly quarter pinched  (the existence of such time is guaranteed by the estimates in \cite{BS} and arguments in the proofs of Theorem 5.1 in \cite{BW} and Corollary 2.3.4 in \cite{Y}), then using Proposition 5.2 in B\"{o}hm-Wilking \cite{BW}, which will reduce the orbifold case to the manifold case. Now we have a diffeomorphism
$H: \mathbb{S}^n// (x,\pm x')=\bar{B} \cup  (\mathbb{S}^n// (x,\pm x')\setminus B) \rightarrow \mathcal{O}'=D_1 \cup \Omega$.  Note that both $f\circ f'$ and $f\circ r_1\circ f'$ extend to a diffeomorphism  from $\bar{B}$  to $D_1$. We choose one of the extensions such that this extension and $H|_{\bar{B}}$  are  either both orientation-preserving or both orientation-reversing if $n$ is even. Then the desired result in this subcase follows from Lemma \ref{lem 3.3}.
 \hfill{$\Box$

 \vspace*{0.2cm}

\noindent{\bf{Remark}.}  By the way, one can show that if $G$ is a nontrivial, finite subgroup of $\text{Isom}(\mathbb{S}^{2k} \times \mathbb{R}^l)$ acting freely on $\mathbb{S}^{2k} \times \mathbb{R}^l$, where $k$ and $l$ are positive integers, then $G \cong\mathbb{Z}_2$.
We can argue as follows. Recall that $\text{Isom}(\mathbb{S}^{2k} \times \mathbb{R}^l) \cong \text{Isom}(\mathbb{S}^{2k})  \times \text{Isom}(\mathbb{R}^l)$; see for example \cite{So}. Suppose that $g\in G$ is  nontrivial.  As $g$ is of finite order, the  $\text{Isom}(\mathbb{R}^l)$-component of $g$ must have a fixed point in $\mathbb{R}^l$; let $q$ be such a fixed point. By restricting the action of $g$ to the slice  $\mathbb{S}^{2k} \times \{q\}$ we see that the $\text{Isom}(\mathbb{S}^{2k})$-component of $g$ must be the antipodal map on $\mathbb{S}^{2k}$. Since the $\text{Isom}(\mathbb{R}^{l})$-component of $g^2$
must have a fixed point in $\mathbb{R}^l$, we see that
$g^2=1$. Similarly, if $g_1, g_2 \in G$ are both nontrivial, then  the $\text{Isom}(\mathbb{S}^{2k})$-component of $g_i$  must be the antipodal map on $\mathbb{S}^{2k}$, $i=1, 2$, and the $\text{Isom}(\mathbb{S}^{2k})$-component of the product $g_1g_2$ must be trivial, but the $\text{Isom}(\mathbb{R}^{l})$-component of $g_1g_2$
 has a fixed point in $\mathbb{R}^l$, so $g_1g_2=1$. As $g_1^2=1$ we see that $g_1=g_2$.

\vspace*{0.2cm}

The following result slightly extends  and strengthens   Proposition 6.17 in \cite{B19}.

\begin{prop}\label{prop 3.6} (cf. Proposition 6.17 in \cite{B19}) Let  $n\geq 5$ and $\varepsilon_4=\varepsilon_4(n)$ be as in Proposition \ref{prop 3.5}. Let  $(M, g)$ be a complete, noncompact manifold of dimension $n\geq 5$ with strictly PIC and weakly PIC2  everywhere and strictly PIC2 at a point. Assume that $N$ is a $2\varepsilon_4$-neck  in $M$  diffeomorphic to $\mathbb{S}^{n-1} \times (0,1)$, and $\Sigma$ is a central cross section of $N$. Then $\Sigma$ bounds a compact domain $\Omega$ in $M$.  Moreover, if $f: \mathbb{S}^{n-1}\rightarrow \Sigma$ is a $2\varepsilon_4$-homothety coming from the $2\varepsilon_4$-neck structure, either $f$ or $f\circ r_1$ extends to a  diffeomorphism  $F: D^n  \rightarrow \Omega $.
\end{prop}

\noindent {\bf Proof}. \ \ By Perelman's proof of the soul conjecture  \cite{P0} we know that under the assumption of this proposition the manifold $M$ is diffeomorphic to $\mathbb{R}^n$. As $\Sigma$ is a smoothly embedded $(n-1)$-sphere in $M$, we see that $\Sigma$ bounds a compact smooth manifold $\Omega$ which is homeomorphic to $D^n$. (In fact, by the solution of
Schoenflies conjecture in dimensions $n\geq 5$, we know that $\Omega$ is diffeomorphic to $D^n$, but we do not need this fact in our proof of Proposition \ref{prop 3.6}; instead our proof implies this fact in our special situation.) Now we use a surgery  argument as in the proof of Proposition \ref{prop 3.5}, with the help of Proposition \ref{prop 3.4} and the Cerf-Palais disk theorem; compare the proof of Proposition 6.17 in \cite{B19}.
 \hfill{$\Box$

 \vspace*{0.2cm}

\noindent {\bf Remark}.   In Proposition \ref{prop 3.6}, if we assume that $(M,g)$ is strictly PIC2 everywhere instead of weakly PIC2 everywhere  and
strictly PIC2 at a point, any $2\varepsilon_4$-neck in $(M,g)$ must be diffeomorphic to
$\mathbb{S}^{n-1} \times (0,1)$  by Theorem A.2 in \cite{B18}.

\vspace*{0.2cm}

 The following definition of $\varepsilon$-caps  extends and strengthens the corresponding definition in  \cite{B19}, and is inspired by  \cite{B19}.

\vspace *{0.2cm}

\noindent {\bf Definition}. (cf. Definition 6.16 in \cite{B19})  Let $\varepsilon_0=\varepsilon_0(n)$ be a small positive constant and $0 < \varepsilon < \frac{1}{4}\varepsilon_0$.
     Given a point $x_0\in (\mathcal{O},g)$, an open subset
  $U$ of $\mathcal{O}$ is an $\varepsilon$-cap centered at $x_0$ if  $U$ is a topological cap and $U\setminus V$ is an $\varepsilon$-neck, where $V$ is a  compact domain in $U$ diffeomorphic to some $D^n//\Gamma$ or  $\mathbb{S}^n// \langle\Gamma, \hat{\sigma} \rangle \setminus B$ or $\mathbb{S}^n// (x,\pm x') \setminus B$ with $\partial V$ a central cross-section of an $\varepsilon$-neck (that is, the image of $\mathbb{S}^{n-1}/\Gamma \times  \{0\}$ under the diffeomorphism associated to the $\varepsilon$-neck) and with $x_0 \in \text{Int}\hspace*{0.5mm} V$, (we will call $V$ a core of the cap $U$,)  and in addition,

  1. when $U$ is diffeomorphic to $\mathbb{R}^n$, if $\Sigma$ is a central  cross section  of  an $\varepsilon_0$-neck in $U\setminus V$ with an  $\varepsilon_0$-homothety $f: \mathbb{S}^{n-1}\rightarrow \Sigma$  coming from the $\varepsilon_0$-neck structure,
    $\Sigma$ bounds  a compact domain $\Omega $ in  $U$, and either $f$ or $f\circ r_1$ extends to a  diffeomorphism  $F: D^n  \rightarrow \Omega $;

 2. when $U$ is diffeomorphic to  $\mathbb{R}P^n \setminus \bar{B}$,  if  $\Sigma$ is a central  cross section  of  an $\varepsilon_0$-neck in $U\setminus V$ with an  $\varepsilon_0$-homothety $f: \mathbb{S}^{n-1}\rightarrow \Sigma$  coming from the $\varepsilon_0$-neck structure,  $\Sigma$ bounds  a compact domain $\Omega $ in  $U$, and either $f\circ f'$ or $f\circ r_1\circ f'$ extends to a  diffeomorphism  $F: \mathbb{R}P^n \setminus B \rightarrow \Omega$, where  $f': \partial (\mathbb{R}P^n \setminus B) \rightarrow \mathbb{S}^{n-1}$ is a homothety;

 3. when $U$ is  diffeomorphic to $\mathbb{S}^n// (x,\pm x') \setminus \bar{B}$,  if   $\Sigma$ is a central  cross section  of  an $\varepsilon_0$-neck in $U\setminus V$ with an $\varepsilon_0$-homothety $f: \mathbb{S}^{n-1}\rightarrow \Sigma$  coming from the $\varepsilon_0$-neck structure, $\Sigma$ bounds  a compact domain $\Omega $ in  $U$, and either $f\circ f'$ or $f\circ r_1\circ f'$ extends to a  diffeomorphism $F: \mathbb{S}^n// (x,\pm x') \setminus B \rightarrow \Omega$, where $f': \partial (\mathbb{S}^n// (x,\pm x') \setminus B) \rightarrow \mathbb{S}^{n-1}$ is a homothety.

\vspace *{0.2cm}

The idea behind this definition is that an $\varepsilon$-cap $U$ should be, after suitable rescaling,  $\eta$-close to the corresponding subset $W$ of an orbifold ancient $\kappa$-solution or an orbifold standard solution with an $\eta$-homothety $\psi: W\rightarrow U$, where $\eta$ is much smaller than $\varepsilon$ (compare p. 522 of \cite{B19} and Section 2.3 of \cite{D08}); moreover, the definition should help one to identify some components which are discarded in the process of the surgery of the Ricci flow  (see Proposition \ref{prop 4.3}). Assume that the end of the $\varepsilon$-cap $U$ is diffeomorphic to $\mathbb{S}^{n-1} \times (0,1)$.   Let $N_1$ be an $\varepsilon_0$-neck given by $\psi_1: \mathbb{S}^{n-1} \times (-\varepsilon_0^{-1}, \varepsilon_0^{-1})\rightarrow N_1$ which is contained in the $\varepsilon$-neck end in $U$. (Note  that by Proposition \ref{gluing1}, given  $0< \varepsilon < \varepsilon' \leq 2\varepsilon_1$ and  an $\varepsilon$-neck diffeomorphic to $\mathbb{S}^{n-1} \times (0,1)$, any $\varepsilon'$-neck contained in this $\varepsilon$-neck is also  diffeomorphic to $\mathbb{S}^{n-1} \times (0,1)$.)  If $\varepsilon_0\leq \varepsilon_4$, by Lemmas \ref{lem 2.6}   and \ref{lem 2.7}, a suitable restriction of $ \psi^{-1}\circ \psi_1$ will give a $2\varepsilon_4$-neck in $W$ if $\eta$ is sufficiently small. Then  we can use Proposition \ref{prop 3.5} or Proposition \ref{prop 3.6} and the map $\psi$ to get the desired extra property for the $\varepsilon$-cap $U$.

\begin{prop}\label{prop 3.7} (cf.  Theorem 6.18 in \cite{B19}  and Theorem 3.9 in \cite{CTZ}) \ \
Given a small positive constant $\varepsilon$  and a constant $\theta >0$, there exist positive constants $C_1=C_1(n, \theta,\varepsilon)$ and $C_2=C_2(n,\theta,\varepsilon)$,
such that given any noncompact orbifold ancient  $\kappa$-solution
$(\mathcal{O},g(t))$ of dimension $n\geq 5$ (with at most isolated singularities) which satisfies $Rm-\theta R \hspace*{1mm} \text{id} \owedge \text{id} \in C_{\text{PIC}}$ and is not locally isometric to an evolving shrinking round cylinder, for each space-time point $(x_0,t_0)$, there
is an open subset $U$ with
$\overline{B(x_0,t_0,C_1^{-1}R(x_0,t_0)^{-\frac{1}{2}})} \subset U \subset B(x_0,t_0, C_1R(x_0,t_0)^{-\frac{1}{2}})$, which falls into
one of the following categories:

(a) $U$ is an  $\varepsilon$-neck diffeomorphic to $\mathbb{S}^{n-1}/\Gamma \times (-1,1)$ centered at $(x_0,t_0)$, or

(b) $U$ is an $\varepsilon$-cap diffeomorphic to $\mathbb{R}^n//\Gamma$  centered at $(x_0,t_0)$

\noindent for some finite subgroup $\Gamma$ of $O(n)$ acting freely on $\mathbb{S}^{n-1}$;
moreover, the scalar curvature in $U$ at time $t_0$ is (strictly) between $C_2^{-1}R(x_0,t_0)$ and $C_2R(x_0,t_0)$, and $vol_{g(t_0)}(U) > (C_2|\Gamma|)^{-1} R(x_0,t_0)^{-n/2}$.
\end{prop}

\noindent {\bf Proof}. \ \ The arguments are adapted from the proof of Theorem 6.18 in \cite{B19}, Theorems 3.9 in \cite{CTZ},  and
Proposition 3.6 in \cite{Hu15}.
 Note that by  Corollary 6.7 in \cite{B19} (see also the first paragraph in the proof of Proposition \ref{prop 3.5} here) we know that $(\mathcal{O}, g(t))$  is strictly PIC2. By an orbifold version of Gromoll-Meyer's soul theorem (see  Proposition 3.4 in \cite{Hu15}), $\mathcal{O}$ is
diffeomorphic to the quotient orbifold $\mathbb{R}^n// \Gamma$ for some finite subgroup
$\Gamma$ of  $O(n)$ acting freely on $\mathbb{S}^{n-1}$. So we have a $\Gamma$-invariant projection $\pi: \mathbb{R}^n \rightarrow \mathcal{O}$.  We pull back the Ricci flow
$(\mathcal{O}, g(t))$ to $(\mathbb{R}^n, \tilde{g}(t))$ via $\pi$, which is then a
$\Gamma$-invariant ancient $\kappa$-solution on a smooth manifold. Here the $\kappa$-noncollapsing property of $(\mathbb{R}^n, \tilde{g}(t))$ follows from the $\kappa$-noncollapsing assumption for $(\mathcal{O}, g(t))$ and the fact that $\pi|_{B(\tilde{p},t,r)}:B(\tilde{p},t,r)\rightarrow  B(p,t,r)$ is surjective for any $p\in \mathcal{O}$,  $\tilde{p}\in \pi^{-1}(p)$, any $t$ (for which $g(t)$ is defined), and any $r>0$ (see for example the proof of
Proposition 3.7 in \cite{CTZ}). We fix  $t_0=0$. For any small $\eta>0$, let $M_\eta$ be the set of all points in $\mathbb{R}^n$ which are not the centers of  $\eta$-necks in $(\mathbb{R}^n, \tilde{g}(0))$. Note that for  $\eta$ sufficiently small, $M_\eta \neq \emptyset$ since it contains a soul of $(\mathbb{R}^n, \tilde{g}(0))$; cf. for example p. 969 in \cite{Hu15} (or by using an extension of Lemma A.10 in \cite{MT}).
 From  Step 1 in the proof of  Theorem 6.18 in \cite{B19} we  know that $M_\eta$ is compact.
Let $y$ be any point in $\partial M_\eta$. From Step 2 in the proof of  Theorem 6.18 in \cite{B19} we have
\begin{equation*}
R(x,t)\leq R(y,0)\omega(R(y,0)d_0(x,y)^2)
\end{equation*}
for all $x\in \mathbb{R}^n$ and $t\leq 0$, where $\omega:[0, \infty)\rightarrow [0,\infty)$ is a function independent of $\kappa$. (For the proof, instead of using Theorem 6.12 in \cite{B19} (compare Remark 45.4 in \cite{KL08}) one can also use an extension of Corollary 9.62 in \cite{MT} in our situation, which follows from the fact that the asymptotic volume ratio of any ancient $\kappa$-solution (see Definition 6.10 in \cite{B19}) is also zero in the higher dimensional case. This fact can be shown by using Theorem 6.1 and Proposition 6.5 in \cite{B19}.) Note that $y$ is the center of a  $2\eta$-neck in $(\mathbb{R}^n, \tilde{g}(0))$, and by \cite{H97} this neck (away from the two ends) is  foliated by almost horizontal constant mean curvature spheres (of dimension $n-1$) in a unique way. Let $\tilde{\Sigma}_y$ be the leaf of this foliation which passes through $y$.  By Theorem G1.1 in \cite{H97} (whose proof applies to any dimension $n\geq 4$) (of course, we can also use the solution of
Schoenflies conjecture in dimensions $n\geq 5$ here), $\tilde{\Sigma}_y$ bounds a compact domain $\tilde{\Omega}_y$ (in
$\mathbb{R}^n$) diffeomorphic to a standard  disk
$D^n$. By Step 3 in the proof of Theorem 6.18 in \cite{B19} $R(y,0) \hspace*{0.4mm} \text{diam}_{\tilde{g}(0)}(\tilde{\Omega}_y)^2\leq C$, where $C=C(n,\theta,\eta)$ is independent of $\kappa$.

As in the proof of Proposition 3.4 in \cite{CZ}, there is a point $y\in \partial M_\eta$ such that $M_\eta \subset \tilde{\Omega}_y$. Fix such a point $y$.
  We rescale the solution $\tilde{g}(t)$ so that
$R(y, 0)=1$ after  rescaling; the rescaled solution will still be denoted by $\tilde{g}(t)$. Note that $\mathbb{R}^n$ ($n\geq 5$) has only one end. For $\eta >0$ sufficiently small (depending on $\varepsilon$), we
can use Lemma \ref{gluing2} and Proposition \ref{gluing3} (see also Lemma C2.1, Corollary C2.3 and Theorem C2.4 in \cite{H97}) repeatedly to get a Hamilton's canonical parametrization $\tilde{\phi}:
\mathbb{S}^{n-1} \times (a, +\infty) \rightarrow \widetilde{\mathcal{T}} \subset \mathbb{R}^n $ which gives an $\varepsilon$-tube  $\widetilde{\mathcal{T}}$  in $(\mathbb{R}^n, \tilde{g}(0))$ containing $\mathbb{R}^n \setminus \text{Int} \hspace*{0.5mm} \tilde{\Omega}_y$ such that $y$ is a center of an $\varepsilon$-neck (with Hamilton's canonical parametrization) contained in $\widetilde{\mathcal{T}}$.   Recall that by  Proposition 4.3 in \cite{DL} (see also p. 967-968 in \cite{Hu15}), there exists a point in $\mathbb{R}^n$ which is a soul of $(\mathbb{R}^n, \tilde{g}(0))$ and which is fixed by $\Gamma$. We denote this point by $O$.
As in p. 969 of \cite{Hu15} one can show  that $O$ is not in the nine-tenths of any $2\eta$-neck in $(\mathbb{R}^n, \tilde{g}(0))$ if $\eta$ is sufficiently small. In particular, $O \in \text{Int} \hspace*{0.5mm} M_\eta \subset \text{Int} \hspace*{0.5mm} \tilde{\Omega}_y$. We will fix an $\eta$ satisfying the desired properties above. We take the Hamilton's canonical parametrization $\tilde{\phi}$  to be maximal  in the sense that $\widetilde{\mathcal{T}}$ cannot be properly included  in the image of any  Hamilton's canonical parametrization which gives an $\varepsilon$-tube in $(\mathbb{R}^n, \tilde{g}(0))$.

 As in the proof of Lemma \ref{gluing2}, the
  parametrization $\tilde{\phi}$ above can
be pushed down  to get a Hamilton's canonical parametrization $\phi: \mathbb{S}^{n-1}/
\Gamma \times (a, +\infty) \rightarrow \mathcal{T} \subset \mathcal{O}$  which gives an $\varepsilon$-tube in $(\mathcal{O}, g(0))$.  (Compare the proof of Theorem 3.9 in \cite{CTZ}; here our proof is somewhat different.)

 For any point $x \in \mathcal{O}\setminus \pi (\text{Int} \hspace*{0.5mm} \tilde{\Omega}_y)$,  we can use the $\phi$ above to give an $\varepsilon$-neck  centered at $x$.

For any $z\in \mathbb{R}^n \setminus \text{Int} \hspace*{0.5mm} \tilde{\Omega}_y$  we denote the almost horizontal constant mean curvature
$(n-1)$-sphere passing through $z$ by $\tilde{\Sigma}_z$, which bounds a compact domain $\tilde{\Omega}_z$ (in
$\mathbb{R}^n$) diffeomorphic to a closed, standard $n$-disk
$D^n$. Let $\tilde{W}$ be the $\varepsilon$-neck with Hamilton's canonical parametrization which is contained in $\mathbb{R}^n \setminus \tilde{\Omega}_y$ such that $\tilde{\Sigma}_y$ is one connected  component of the boundary of the closure of $\tilde{W}$. Let $z_0$ be any point in the connected  component of the boundary of the closure of $\tilde{W}$ which is different from $\tilde{\Sigma}_y$.
 From above we know that $\tilde{\Omega}_{z_0}$ is $\Gamma$-invariant, and $\text{Int} \hspace*{0.5mm} \tilde{\Omega}_{z_0} // \Gamma$ contains the $\varepsilon$-neck $W:=\pi(\tilde{W})$ near its end with $\partial \overline{W} \cap \text{Int} \hspace*{0.5mm} \tilde{\Omega}_{z_0}// \Gamma$
 ($=\pi (\tilde{\Sigma}_y)$) a central cross-section of an $\varepsilon$-neck. Now we can use the heat equation method  (cf. Appendixes  F and G in \cite{H97} and the proof of
Theorem 3.9 in \cite{CTZ}) to smooth the Busemann function  based at  the $\Gamma$-invariant soul $O$ a bit to a $\Gamma$-invariant smooth function $\tilde{b}$, and fuse $\tilde{b}$ with the $\Gamma$-invariant height function on the $\varepsilon$-neck centered at $z_0$ and parametrized by certain restriction of the map $\tilde{\phi}$ (or this height function multiplied by a constant and/or shifted by a constant if necessary) using a suitable cut-off function,  and  an equivariant Morse theoretical  argument using (roughly speaking) the gradient flow of the resulting function to show that $\tilde{\Omega}_{z_0}$ is
$\Gamma$-equivariantly diffeomorphic to $D^n$ (if $\varepsilon$ is small enough, which we may assume); compare the proof of Theorem G1.1 in \cite{H97}, cf. also the proof of Theorem 3.9 in \cite{CTZ}, Proposition 3.6 in \cite{Hu15} and Proposition \ref{gluing1} in the appendix.
(By the way, we would like to point out that the reaction term on the RHS of the evolution equation for the Hessian of $u$ (where $u$ is a solution to the heat equation coupled with the Ricci flow) on p. 56 of \cite{CTZ} satisfies the so called null-eigenvector assumption (Definition 10.6 in \cite{C+}) under the PIC2 condition.) So $\text{Int} \hspace*{0.5mm} \tilde{\Omega}_{z_0} // \Gamma$ is diffeomorphic to $\mathbb{R}^n // \Gamma$, and  is an $\varepsilon$-cap; in the case that $\Gamma$ is trivial,  that  $\text{Int} \hspace*{0.5mm} \tilde{\Omega}_{z_0}$  has the additional property in the definition of $\varepsilon$-cap is guaranteed by Proposition \ref{prop 3.6}.

From above, the Harnack inequality
(Corollary 6.3 in \cite{B19}), and the fact that $\pi|_{B(\tilde{p},0,r)}:B(\tilde{p},0,r)\rightarrow  B(p,0,r)$ is surjective for any $p\in \mathcal{O}$,  $\tilde{p}\in \pi^{-1}(p)$ and any $r>0$, we see that $C'^{-1}R(y,0)\leq R(\bar{z},0)\leq C'R(y,0)$ and $\text{diam}_{g(0)}(\tilde{\Omega}_{z_0} // \Gamma) \leq C'R(\bar{z},0)^{-\frac{1}{2}}$  for any $\bar{z}\in \tilde{\Omega}_{z_0} // \Gamma$, where  $C'=C'(n,\theta,\varepsilon)$ is a positive constant independent of $\kappa$. In particular, $R(\bar{z},0)$ and $R(\bar{z}',0)$ are comparable if both $\bar{z}$ and $\bar{z}'$ are points in $\tilde{\Omega}_{z_0} // \Gamma$.  We also observe that
if $\tilde{z}$ is a center of  an $\varepsilon$-neck (for the metric $\tilde{g}(0)$) contained in the end of $\tilde{\Omega}_{z_0}$,
\begin{equation*}
\text{vol}_{\tilde{g}(0)}(\tilde{\Omega}_{z_0}) > \text{vol}_{\tilde{g}(0)}(B_{\tilde{g}(0)}(\tilde{z}, R(\tilde{z},0)^{-1/2})) \geq \mu(n) R(\tilde{z},0)^{-n/2},
\end{equation*}
where $\mu(n)$ is a positive constant depending on $n$.  If $x_0\in \pi(\text{Int} \hspace*{0.5mm} (\tilde{\Omega}_y))$, $\text{Int} \hspace*{0.5mm} \tilde{\Omega}_{z_0}// \Gamma$ is an $\varepsilon$-cap centered at $x_0$.
Now the existence of the desired constants $C_1$ and $C_2$ is clear.
 \hfill{$\Box$}

 \vspace*{0.2cm}

 The following proposition is a higher dimensional analogue  of  Proposition 3.7 in \cite{CTZ}, and extends Corollary 6.20 in \cite{B19}; for the 3-dimensional case see  Section 1.5 in \cite{P2}.

\begin{prop} \label{prop 3.8} \ \  There exists a   positive function $\omega: [0,\infty)\rightarrow (0,\infty)$ and a positive constant $C_3$ both depending only on $n$ and $\theta$ with the following property. Assume that $(\mathcal{O},g(t)_{t\in (-\infty, 0]})$  is an orbifold ancient $\kappa$-solution
 of dimension $n\geq 5$ (with at most isolated singularities) which satisfies $Rm-\theta R \hspace*{1mm} \text{id} \owedge \text{id} \in C_{\text{PIC}}$.  Then

 (i) for any $x,y\in \mathcal{O}$ and any  $t\in (-\infty,0]$ we have
 \begin{equation*}
  R(x,t)\leq R(y,t)\omega(R(y,t)d_t(x,y)^2);
\end{equation*}

(ii) for any $x\in \mathcal{O}$ and any  $t\in (-\infty,0]$ we have
\begin{equation*}
|\nabla R|(x,t)< C_3 R(x,t)^{\frac{3}{2}} \hspace*{8mm} \mbox{and} \hspace*{8mm}
|\frac{\partial R}{\partial t}|(x,t)< C_3 R(x,t)^2.
\end{equation*}

 \end{prop}

\noindent {\bf Proof}. \ \ For (i), we use  Proposition 3.4 in \cite{Hu15}, an orbifold version of the sphere theorem in \cite{BS} (see the proof of Proposition \ref{prop 3.5} above) and Theorems 6.13 (in the proof of this theorem, instead of Theorem 6.12 in \cite{B19} one can also use an extension of Corollary 9.62 in \cite{MT} in our situation) and 6.19 in \cite{B19}, and argue as in the proof of Proposition 3.7 in \cite{CTZ}.  The conclusion in (ii) follows from (i) and Shi's interior derivative estimates. \hfill{$\Box$

\vspace*{0.2cm}

In summary, we have the following description of the canonical neighborhood property of orbifold ancient $\kappa$-solutions.

\begin{prop}\label{prop 3.9} (cf. \cite{P2}, Corollaries 6.20 and 6.22 in \cite{B19}, and Theorem 3.10 in \cite{CTZ}) \ \
Given   a small positive constant $\varepsilon$  and a constant $\theta >0$, there exist positive constants $C_1=C_1(n, \theta,\varepsilon)$, $C_2=C_2(n,\theta,\varepsilon)$ and $C_3=C_3(n,\theta)$ with the following property:
Suppose $(\mathcal{O},g(t))$ is an orbifold ancient $\kappa$-solution of dimension $n\geq 5$ (with at most isolated singularities) which satisfies
$Rm-\theta R \hspace*{1mm} \text{id} \owedge \text{id} \in C_{\text{PIC}}$. Then either $(\mathcal{O},g(t))$ is  compact and strictly PIC2 for any $t$ (hence diffeomorphic to a spherical orbifold), or for each space-time point $(x_0,t_0)$, there
is  an open subset $U$ of $\mathcal{O}$ with
$\overline{B(x_0,t_0,C_1^{-1}R(x_0,t_0)^{-\frac{1}{2}})} \subset U \subset B(x_0,t_0, C_1R(x_0,t_0)^{-\frac{1}{2}})$ and with the values of the scalar curvature in $U$ at time $t_0$ lying (strictly) between $C_2^{-1}R(x_0,t_0)$ and $C_2R(x_0,t_0)$, which falls into
one of the following two categories:

(a) $U$ is a strong  $\varepsilon$-neck diffeomorphic to $\mathbb{S}^{n-1}/\Gamma \times (-1,1)$  centered at $(x_0,t_0)$, or

(b) $U$ is an $\varepsilon$-cap centered at $(x_0,t_0)$ with any cross-section of the $\varepsilon$-neck contained in it diffeomorphic to $\mathbb{S}^{n-1}/\Gamma$

\noindent for some finite subgroup $\Gamma$ of $O(n)$ acting freely on $\mathbb{S}^{n-1}$, and
\begin{equation*}
vol_{g(t_0)}(U) > (C_2|\Gamma|)^{-1} R(x_0,t_0)^{-n/2};
\end{equation*}
moreover, the scalar curvature in $U$ at time $t_0$ satisfies
the derivative estimates
\begin{equation*}
|\nabla R|< C_3 R^{\frac{3}{2}} \hspace*{8mm} \mbox{and} \hspace*{8mm}
|\frac{\partial R}{\partial t}|< C_3 R^2.
\end{equation*}

\end{prop}

\noindent {\bf Proof}. \ \ Let $(\mathcal{O}, g(t))$ be as in the assumption. If $\mathcal{O}$ is compact,  by the $\kappa$-noncollapsing assumption on all scales it cannot be isometric to a metric quotient of the round cylinder. Then by an orbifold version of Corollary 6.7 in \cite{B19}  in this case $(\mathcal{O},g(t))$ is strictly PIC2 for any $t$, so $\mathcal{O}$  is diffeomorphic to a spherical orbifold.  If $\mathcal{O}$ is noncompact, the result except for the  derivative estimates follows from Propositions \ref{prop 3.5} and \ref{prop 3.7} and Lemma \ref{lem 2.5}, and the derivative estimates
follows from Proposition \ref{prop 3.8} (ii).
 \hfill{$\Box$

\vspace*{0.2cm}

Note that the constants $C_1$, $C_2$ and $C_3$ in Proposition \ref{prop 3.9} are independent of $\kappa$.

\vspace*{0.2cm}

It turns out  that orbifold standard solutions (see the appendix) have canonical neighborhood  property similar to that of orbifold ancient  $\kappa$-solutions, except that an $\varepsilon$-neck in an orbifold standard solution may not be strong.

\begin{prop}\label{prop 3.10} (cf. \cite{P2}, Corollary A.2 in \cite{CZ}, and Corollary 9.3 in \cite{B19}) \ \
Given   a small positive constant $\varepsilon$,  there exist positive constants $C_1'=C_1'(n,\varepsilon)$, $C_2'=C_2'(n,\varepsilon)$ and $C_3'=C_3'(n)$ with the following property: For each space-time point $(x_0,t_0)$ on an orbifold standard solution
$(\mathbb{R}^n//\Gamma,\hat{g}_\Gamma(t))$ (with at most an isolated singularity), there
is  an open subset $U$ of $\mathbb{R}^n//\Gamma$ with
$\overline{B(x_0,t_0,C_1'^{-1}R(x_0,t_0)^{-\frac{1}{2}})} \subset U \subset B(x_0,t_0, C_1'R(x_0,t_0)^{-\frac{1}{2}})$ and with the values of the scalar curvature in $U$ at time $t_0$ lying (strictly) between $C_2'^{-1}R(x_0,t_0)$ and $C_2'R(x_0,t_0)$, which falls into
one of the following two categories:

(a) $U$ is an   $\varepsilon$-neck diffeomorphic to $\mathbb{S}^{n-1}/\Gamma \times (-1,1)$ centered at $(x_0,t_0)$,  and on $U \times [t_0-\min \{R(x_0,t_0)^{-1},t_0\},t_0]$ the solution $\hat{g}_\Gamma(t)$ is, after scaling with $R(x_0,t_0)$ and shifting  $t_0$ to zero, $\varepsilon$-close to the corresponding subset of the evolving round cylinder $\mathbb{S}^{n-1}/\Gamma \times \mathbb{R}$ over the time interval $[-\min \{t_0R(x_0,t_0), 1\}, 0]$ with scalar curvature 1 at the time zero, and $U$ is disjoint from the surgery cap (in particular, $x_0 \notin B(p_\Gamma, 0, \varepsilon^{-1})$, where, as in the appendix, $p_\Gamma$ denotes the tip of the orbifold standard solution) when $t_0\leq R(x_0,t_0)^{-1}$, or

(b) $U$ is an $\varepsilon$-cap diffeomorphic to $\mathbb{R}^n//\Gamma$ centered at $(x_0,t_0)$,

\noindent and $vol_{\hat{g}_\Gamma(t_0)}(U) > (C_2'|\Gamma|)^{-1} R(x_0,t_0)^{-n/2}$; moreover, the scalar curvature in $U$ at time $t_0$ satisfies
the derivative estimates
\begin{equation*}
|\nabla R|< C_3' R^{\frac{3}{2}} \hspace*{8mm} \mbox{and} \hspace*{8mm}
|\frac{\partial R}{\partial t}|< C_3' R^2.
\end{equation*}

\end{prop}

\noindent {\bf Proof}. \ \   The conclusion in the case that $\Gamma$ is  trivial follows from Corollary A.2 in \cite{CZ}, Corollary 9.3 in \cite{B19}, and Proposition \ref{prop 3.6}.  Now we consider the case  that $\Gamma$ is not trivial. By Corollary A.2 in \cite{CZ} and Corollary 9.3 in \cite{B19}, for any $\eta >0$  there exist positive constants $\tilde{C}_1=\tilde{C}_1(n,\eta)$, $\tilde{C}_2=\tilde{C}_2(n,\eta)$ and $C_3'=C_3'(n)$ with the following property: For each space-time point $(\tilde{x}_0,t_0)$ on the standard solution
$(\mathbb{R}^n,\hat{g}(t))$, there
is  an open subset $\tilde{U}$ of $\mathbb{R}^n$ with
$\overline{B(\tilde{x}_0,t_0,\tilde{C}_1^{-1}R(\tilde{x}_0,t_0)^{-\frac{1}{2}})} \subset \tilde{U} \subset B(\tilde{x}_0,t_0, \tilde{C}_1R(\tilde{x}_0,t_0)^{-\frac{1}{2}})$ and with the values of the scalar curvature in $\tilde{U}$ at time $t_0$ lying between $\tilde{C}_2^{-1}R(\tilde{x}_0,t_0)$ and $\tilde{C}_2R(\tilde{x}_0,t_0)$, which falls into
one of the following two categories:

($\tilde{a}$) $\tilde{U}$ is an   $\eta$-neck  centered at $(\tilde{x}_0,t_0)$,  and on $\tilde{U} \times [t_0-\min \{R(\tilde{x}_0,t_0)^{-1},t_0\},t_0]$ the solution $\hat{g}(t)$ is, after scaling with $R(\tilde{x}_0,t_0)$ and shifting  $t_0$ to zero, $\eta$-close to the corresponding subset of the evolving round cylinder $\mathbb{S}^{n-1} \times \mathbb{R}$ over the time interval $[-\min \{t_0R(\tilde{x}_0,t_0), 1\}, 0]$ with scalar curvature 1 at the time zero, and $\widetilde{U}$ is disjoint from the surgery cap when $t_0\leq R(\tilde{x}_0,t_0)^{-1}$, or

($\tilde{b}$) $\tilde{U}$ is an $\eta$-cap centered at $(\tilde{x}_0,t_0)$,

\noindent and $\text{vol}_{\hat{g}(t_0)}(\tilde{U}) > \tilde{C}_2^{-1} R(\tilde{x}_0,t_0)^{-n/2}$; moreover, the scalar curvature in $\tilde{U}$ at time $t_0$ satisfies
the derivative estimates
\begin{equation*}
|\nabla R|< C_3' R^{\frac{3}{2}} \hspace*{8mm} \mbox{and} \hspace*{8mm}
|\frac{\partial R}{\partial t}|< C_3' R^2.
\end{equation*}

Now by above and Lemma C2.1, Corollary C2.3 and Theorem C2.4 in \cite{H97} (see also the appendix) one can get a Hamilton's canonical parametrization
\begin{equation*}
\tilde{\phi}:\mathbb{S}^{n-1} \times (a,+\infty) \rightarrow  (\mathbb{R}^n \setminus \overline{B(p_0,t_0, C(n,\eta)R(p_0,t_0)^{-\frac{1}{2}})}, \hat{g}(t_0)),
 \end{equation*}
 where $p_0$ denotes the tip of the manifold standard solution as in the appendix, and $C(n,\eta)$ is a large number. As $(\mathbb{R}^n, \hat{g}(t_0))$ is $O(n)$-invariant, $\tilde{\phi}$ is  $\Gamma$-equivariant, where the $\Gamma$-action on the second factor of $\mathbb{S}^{n-1} \times (a,+\infty)$ is trivial.  If $\eta$ is sufficiently small (depending on $\varepsilon$ and $n$), the quotient of $\tilde{\phi}$ by $\Gamma$, denoted by
\begin{equation*}
\phi: \mathbb{S}^{n-1}/\Gamma \times (a,+\infty) \rightarrow  (\mathbb{R}^n//\Gamma \setminus \overline{B(p_\Gamma,t_0, C(n,\eta)R(p_0,t_0)^{-\frac{1}{2}})}, \hat{g}_\Gamma(t_0)),
\end{equation*}
 will give a $\beta \varepsilon$-tube, where $\beta$ is from Lemma \ref{lem 2.5}. (Compare  Corollary C2.3 in \cite{H97} and the proof of Lemma \ref{gluing2}.) Moreover, if   $(x_0,t_0)$ is a   center of a   $\beta \varepsilon$-neck   contained in this $\beta \varepsilon$-tube with $t_0 > \frac{n}{n+2}$, which implies that $t_0R(x_0,t_0)\geq t_0\frac{1}{1-\frac{2}{n}t_0} >1$, by using Lemma \ref{lem 2.5} (choosing $b=0$ there) we see that $(x_0,t_0)$ will fall into the category (a);  if the constant $C(n,\eta)$ above is sufficiently large (depending on $\varepsilon$) and $(x_0,t_0)$ is a   center of a   $\beta  \varepsilon$-neck   contained in this $\beta  \varepsilon$-tube with $t_0 \leq \frac{n}{n+2}$, by using the fact that the orbifold standard solution $(\mathbb{R}^n//\Gamma,\hat{g}_\Gamma(\cdot))$  restricted to the time interval $[0, \frac{n}{n+2}]$ is asymptotic to the evolving round quotient cylinder $\mathbb{S}^{n-1}/\Gamma \times \mathbb{R}$ (with scalar curvature 1 at time 0) restricted to the same time interval  at infinity, which can be shown, for example,  by working $O(n)$-equivariantly with a manifold standard solution and using the arguments in the last paragraph on p. 262 of \cite{CZ} and in Section 3.2 of \cite{BHZ}, we see that  $(x_0,t_0)$ will still fall into the category (a).

 Note that for some constant $C'(n, \varepsilon)> C(n,\eta)$,  $\partial \overline{B(p_\Gamma,t_0, C'(n,\varepsilon)R(p_0,t_0)^{-\frac{1}{2}})}$, which has constant mean curvature,  is the central cross section of an $\varepsilon$-neck with Hamilton's canonical parametrization  $\phi_\Gamma$.   Furthermore, since the metric ball closure
$\overline{B(p_0,t_0, C'(n,\varepsilon)R(p_0,t_0)^{-\frac{1}{2}})}$ is $O(n)$-invariant, there exists a $\Gamma$-equivariant diffeomorphism between it and a standard disk $D^n$. So  $\overline{B(p_\Gamma,t_0, C'(n,\varepsilon)R(p_0,t_0)^{-\frac{1}{2}})}$ is diffeomorphic to $ D^n//\Gamma$.   Now the desired result follows.    Compare the proof of Proposition \ref{prop 3.7}.
\hfill{$\Box$

\vspace*{0.2cm}

Note that the constants $C_1'$, $C_2'$ and $C_3'$ in Proposition \ref{prop 3.10} are independent of $\Gamma$.

\section{Existence of $(r, \delta)$-surgical solutions}

The following lemma should be known (for  relevant results in the manifold case see for example Proposition 7.7 in \cite{B}, Theorem 2.2 in Chapter 8 of \cite{Hi}, Proposition 7.6.6 in  \cite{Mu} and p. 76 of \cite{CTZ}).
\begin{lem} \label{lem 4.1}
Let $P$ be a compact  orbifold  with one of its boundary components $\partial_1 P$ a smooth manifold, and $\psi: P\supset \partial_1 P \rightarrow    \partial_1 P \times \{0\} \subset \partial_1 P \times [0,1]$ be a diffeomorphism.  Then $P\cup_\psi \partial_1 P \times [0,1]$ is diffeomorphic to $P$.
\end{lem}

\noindent {\bf Proof}.\ \ Let $\varphi: \partial_1 P \times [0,1)\rightarrow P$ be a collar neighborhood of $\partial_1 P$ in $P$, that is, $\varphi$ is an embedding with $\varphi (x,0)=x$ for any $x\in \partial_1 P$.  Let $A=\varphi (\partial_1 P \times [0, \frac{1}{2}))$. Then there is a diffeomorphism $f_1: P\setminus A \rightarrow P$ with $f_1(\varphi(x, \frac{1}{2}))=x$ for any $x\in \partial_1 P$.  For convenience we also view $\psi$ as a self-diffeomorphism of $\partial_1 P $  via the obvious identification of  $\partial_1 P $ with  $\partial_1 P \times \{0\}$  given by $\partial_1 P  \ni x \mapsto (x,0)\in \partial_1 P \times \{0\}$.  Let $f_2: \bar{A}=\varphi (\partial_1 P \times [0, \frac{1}{2}])\rightarrow \partial_1 P \times [0,1]$ be given by $f_2(\varphi (x,t))=(\psi(x), 1-2t)$ for $x \in \partial_1 P$ and $t\in [0,\frac{1}{2}]$.  Note that  $f_2$ is also a diffeomorphism. Let $f_3: P\rightarrow  P\cup_\psi \partial_1 P \times [0,1]$ be given by

\begin{equation*}
f_3(y)=\left\{
\begin{aligned}
f_1(y), \hspace*{8mm}   & y  \in P\setminus A,  \\
f_2(y),  \hspace*{8mm}  & y  \in \bar{A}. \\
\end{aligned}
\right.
\end{equation*}
Then $f_3$ is a homeomorphism.  Using the ambient tubular neighborhood theorem (actually here it suffices to apply the ambient tubular neighborhood theorem to the manifold $P\cup_\psi \partial_1 P \times [0,1]\setminus \{\text{singularities}\}$) and arguing as in the proof of Theorem 1.9 in Chapter 8 of \cite{Hi} and Proposition 7.6.4 in  \cite{Mu}, we see that $f_3$
can be modified to a diffeomorphism   $f_4: P\rightarrow  P\cup_\psi \partial_1 P \times [0,1]$.
\hfill{$\Box$}

 \vspace*{0.2cm}

Now we try to identify the topology of a compact orbifold which is covered by $\varepsilon$-caps and/or $\varepsilon$-necks.

\begin{prop} \label{prop 4.3}\ \  Let $\varepsilon_0=\min \{\varepsilon_3, \varepsilon_4\}$, where the constant $\varepsilon_3$ is as in Proposition \ref{gluing3},  and the constant $\varepsilon_4$ is as in Proposition \ref{prop 3.5}.  Fix $0 < \varepsilon <  \tilde{\varepsilon}_1(\varepsilon_0)$, where $\tilde{\varepsilon}_1(\cdot)$ is as in Lemma \ref{gluing2}.
Let $(\mathcal{O},g)$ be a closed, connected $n$-orbifold with at most isolated singularities. If each point of $\mathcal{O}$ is a center of an
$\varepsilon$-neck or an $\varepsilon$-cap, then $\mathcal{O}$ is diffeomorphic  to a spherical orbifold or a weak connected sum of at most two spherical orbifolds with at most isolated singularities.  In particular, if each point of $\mathcal{O}$ is a center of an
$\varepsilon$-neck, then $\mathcal{O}$ is diffeomorphic to  a  quotient manifold of $\mathbb{S}^{n-1} \times \mathbb{R}$ by standard isometries.
\end{prop}

\noindent {\bf Proof}.  Compare the proof of Proposition 8.1 in \cite{B19}. Let $\mathcal{O}$ satisfy the assumption of our proposition.
 If $\mathcal{O}$ contains (at least) an  $\varepsilon$-cap of type $\mathbb{R}^n$  but no $\varepsilon$-caps of the other types, by using  the definition of $\varepsilon$-cap of type $\mathbb{R}^n$ (and our choice of $\varepsilon$),  Hamilton's canonical parametrization of an $\varepsilon$-tube (cf. Section C.2 in \cite{H97}; see also Lemma \ref{gluing2} and Proposition \ref{gluing3} in the appendix), and Theorem 1.5 in Chapter 8 of \cite{Hi} (a corollary of the isotopy extension theorem), we see that it is diffeomorphic to $D^n\cup_h D^n$, where $h$ is an isometry of $\mathbb{S}^{n-1}$, so it is diffeomorphic to $\mathbb{S}^n$. (To be more precise, here we also need to use, for example, a connectedness argument. Let $W$ be a topological cap which is  the union of an $\varepsilon$-cap and an $\varepsilon_0$-tube (with Hamilton's canonical parametrization). Then a cross section $S$ of the $\varepsilon_0$-tube is separating in $W$, and it divides $W$ into two (connected) parts, both with boundary $S$. Let $X$ be an  embedded compact, connected suborbifold (of codimension 0) in $W$ with boundary  $S$. Then $X\setminus S$ is still path-connected. It follows that $X$ is contained in one of the two parts of $W$ divided by $S$. Since $X$ is compact (and closed), and open relative to the part it lies in, it must occupy a whole part (of course it is the compact part).) If $\mathcal{O}$ contains an $\varepsilon$-cap of type $\mathbb{R}^n$ and an  $\varepsilon$-cap diffeomorphic to $\mathbb{R}P^n\setminus \bar{B}$,  by using  the definition of these two $\varepsilon$-caps A cross section of Hamilton's canonical parametrization should separate   the union of an $\varepsilon$-cap and an $\varepsilon_0$-tube with Hamilton's canonical parametrization. If $\mathcal{O}$ contains an $\varepsilon$-cap of type $\mathbb{R}^n$ and an  $\varepsilon$-cap diffeomorphic to $\mathbb{R}P^n\setminus \bar{B}$,  by using  the definition of these two $\varepsilon$-caps, Hamilton's canonical parametrization of an $\varepsilon$-tube, and Theorem 1.5 in Chapter 8 of \cite{Hi},  we see that it is diffeomorphic to $(\mathbb{R}P^n\setminus B)\cup_h D^n$, where $h: \partial(\mathbb{R}P^n\setminus B) \rightarrow \partial D^n $ is a homothety, so it is diffeomorphic to
 $\mathbb{R}P^n$. Similarly, if $\mathcal{O}$ contains a smooth $\varepsilon$-cap of type $\mathbb{R}^n$ and an orbifold  $\varepsilon$-cap of type II, it  is diffeomorphic to   $\mathbb{S}^n// (x,\pm x')$.

  If $\mathcal{O}$  contains (at least) an $\varepsilon$-cap diffeomorphic to $\mathbb{R}P^n\setminus \bar{B}$ but no caps of the other types, similarly  we see that it is diffeomorphic to $(\mathbb{R}P^n\setminus B)\cup_h  (\mathbb{R}P^n\setminus B)$, where $h$ is an isometry of $ \partial(\mathbb{R}P^n\setminus B)$, hence it is diffeomorphic to  $\mathbb{R}P^n \sharp \mathbb{R}P^n$; in this case  it can be viewed as  an orbifiber bundle over $\mathcal{I}$ with generic fiber $\mathbb{S}^{n-1}$ (and two exceptional fibers both  diffeomorphic to $\mathbb{R}P^{n-1}$), and has universal covering  diffeomorphic to $\mathbb{S}^{n-1}\times \mathbb{R}$, where  $\mathcal{I}$ is the one dimensional closed orbifold with two singular points both with  local group $\mathbb{Z}_2$, and  $|\mathcal{I}|$ is a closed interval.

  Similarly, if $\mathcal{O}$  contains an $\varepsilon$-cap diffeomorphic to $\mathbb{R}P^n\setminus \bar{B}$
  and an orbifold $\varepsilon$-cap of type II,  it  is diffeomorphic to   $\mathbb{R}P^n  \sharp \mathbb{S}^n// (x,\pm x')$; by the way,  in this case it can be viewed as  another orbifiber bundle over $\mathcal{I}$ with generic fiber $\mathbb{S}^{n-1}$ (and two exceptional fibers, one is diffeomorphic to $\mathbb{R}P^{n-1}$, the other is diffeomorphic to a $(n-1)$-dimensional spherical orbifold with two singularities).

  If $\mathcal{O}$  contains (at least) an orbifold  $\varepsilon$-cap  of type II  but no $\varepsilon$-caps  of the other types, similarly it is diffeomorphic to $(\mathbb{S}^n// (x,\pm x') \setminus B) \cup_h (\mathbb{S}^n// (x,\pm x') \setminus B)$, where $h$ is an isometry of $\partial (\mathbb{S}^n// (x,\pm x') \setminus B)$, hence it is diffeomorphic to a connected sum  $\mathbb{S}^n// (x,\pm x')\sharp \mathbb{S}^n// (x,\pm x') $ (where the connected sum occurs at two regular points); in this case  it can also be viewed as an orbifiber bundle over $\mathcal{I}$ with generic fiber $\mathbb{S}^{n-1}$.

   If $\mathcal{O}$ contains (at least) an  $\varepsilon$-cap of type  $C_\Gamma$ ($|\Gamma|\geq 2$) but no $\varepsilon$-caps of the other types, by using Proposition  \ref{gluing1} and Lemma \ref{lem 4.1}  we see that it  is diffeomorphic to  $D^n// \Gamma \cup_h D^n// \Gamma$, where $h$ is a diffeomorphism of $\mathbb{S}^{n-1}/\Gamma$, and   can be viewed as  a weak connected sum of the form $\mathbb{S}^n//\Gamma \tilde{\sharp} \mathbb{S}^n//\Gamma $ (where the weak connected sum occurs at two singular points).

Similarly, if $\mathcal{O}$ contains  an  $\varepsilon$-cap of type  $C_\Gamma$ ($|\Gamma|\geq 2$) and an $\varepsilon$-cap of  type $C_\Gamma^\sigma$,  we see that it  is diffeomorphic to $(\mathbb{S}^n// \langle\Gamma, \hat{\sigma} \rangle \setminus B)\cup_h D^n//\Gamma$ (the notation is as in Section 2), where $h$ is a diffeomorphism between the boundaries.   It can be viewed as  a weak connected sum of the form  $\mathbb{S}^n// \langle\Gamma, \hat{\sigma}\rangle \tilde{\sharp} \mathbb{S}^n//\Gamma$ (where the weak connected sum occurs at two singular points).

If $\mathcal{O}$ contains  (at least)  an  $\varepsilon$-cap of type   $C_\Gamma^\sigma$ ($|\Gamma|\geq 2$) but no caps of the other types,  by using Proposition  \ref{gluing1} and Lemma \ref{lem 4.1}  we see that it  is a smooth manifold diffeomorphic to $(\mathbb{S}^n// \langle\Gamma, \hat{\sigma} \rangle \setminus B)\cup_h (\mathbb{S}^n// \langle\Gamma, \hat{\sigma'} \rangle \setminus B')$ (the notation is as in Section 2), where $h$ is a diffeomorphism between the boundaries, and has the structure of an orbifiber bundle over  $\mathcal{I}$ with generic fiber diffeomorphic to $\mathbb{S}^{n-1}/\Gamma$.   Clearly it can be viewed as a weak connected sum of the form  $\mathbb{S}^n// \langle\Gamma, \hat{\sigma}\rangle \tilde{\sharp} \mathbb{S}^n// \langle\Gamma, \hat{\sigma'}\rangle$ (where the weak connected sum occurs at two singular points).

  By using  Lemma \ref{gluing2}, Proposition \ref{gluing3} in the appendix and Theorem C2.5 in \cite{H97} we see that  if each point of $\mathcal{O}$ is a center of an $\varepsilon$-neck, then $\mathcal{O}$ is diffeomorphic  to  a  quotient manifold  of $\mathbb{S}^{n-1} \times \mathbb{R}$ by standard isometries, and has the structure of   a $\mathbb{S}^{n-1}/\Gamma$-bundle over $\mathbb{S}^1$ with structure group $\text{Isom}(\mathbb{S}^{n-1}/\Gamma)$  (by the way, when the fiber is diffeomorphic to $\mathbb{S}^{n-1}$ or $\mathbb{R}P^{n-1}$, one can also consult  \cite{FGKO} or \cite{CGGK});   in this case it can  be viewed as a  connected sum on a single orbifold $\mathbb{S}^n//\Gamma$, that is, a  connected sum of the form $\mathbb{S}^n// \Gamma \sharp_f$.   \hfill{$\Box$}

\vspace *{0.2cm}

 Fix a compact Riemannian orbifold  $(\mathcal{O},g_0)$  of dimension $n\geq 12$ with positive isotropic curvature and with at most isolated singularities. Let $(\mathcal{O},g(t))$, $t\in [0,T_{\max})$, be the maximal solution to the (smooth)  Ricci flow starting with $(\mathcal{O},g_0)$.
Let $\hat{T}=\frac{n}{2\inf_{x\in \mathcal{O}}R(x,0)}$. Then $T_{\max}\leq \hat{T}$. By  Theorem 1.2 in \cite{B19} there is a continuous family of closed, convex, $O(n)$-invariant sets $\{\mathcal{F}_t ~|~ t\in [0,\hat{T}]\}$ in the vector space $\mathcal{C}_B(\mathbb{R}^n)$ of algebraic curvature tensors in dimension $n$ such that the curvature tensor  of $(\mathcal{O},g_0)$ lies in $\mathcal{F}_0$, the family is invariant under the Hamilton ODE $\frac{d}{dt}\text{Rm}=Q(\text{Rm})$, and
\begin{equation}
 \mathcal{F}_t \subset \{\text{Rm} \hspace*{1mm}|\hspace*{1mm} \text{Rm}-\theta R \hspace*{1mm} \text{id}\owedge \text{id} \in C_{\text{PIC}}\} \cap  \{\text{Rm} \hspace*{1mm}|\hspace*{1mm} \text{Rm}+f(R)\hspace*{1mm} \text{id}\owedge \text{id} \in C_{\text{PIC}2}\}
\end{equation}
 for any  $t \in [0,\hat{T}]$, where  $f:[0,\infty)\rightarrow [0,\infty)$ is an increasing concave  function satisfying $\text{lim}_{s\rightarrow \infty}\frac{f(s)}{s}=0$, and $\theta$ is a positive number. By (a version of) Hamilton's maximum principle the curvature tensor of $(\mathcal{O},g(t))$ lies in $\mathcal{F}_t$ for any  $t\in [0,T_{\max})$.

 This motivates the following definition.

\vspace *{0.2cm}

\noindent {\bf Pinching assumption} (see Definition 7.1 in \cite{B19}):  Let $f:[0,\infty)\rightarrow [0,\infty)$ be an increasing concave  function satisfying $\text{lim}_{s\rightarrow \infty}\frac{f(s)}{s}=0$, and $\theta$ be a positive number. A Riemannian orbifold satisfies the $(f,\theta)$-pinching assumption  (or has $(f,\theta)$-pinched curvature) if $\text{Rm}+f(R)\hspace*{1mm} \text{id}\owedge \text{id} \in C_{\text{PIC}2}$ and  $\text{Rm}-\theta R \hspace*{1mm} \text{id}\owedge \text{id} \in C_{\text{PIC}}$.  An evolving Riemannian orbifold $(\mathcal{O}(t), g(t))$, $t\in I$, has  $(f,\theta)$-pinched curvature if it has $(f,\theta)$-pinched curvature at each time $t\in I$.

\vspace *{0.2cm}

Now we give the definition of $(\varepsilon,C)$-canonical neighborhood.

\vspace *{0.2cm}

\noindent {\bf Definition} (cf. \cite{P2} and p. 841 in \cite{CH}).  Let $0< \varepsilon < \frac{1}{4}\varepsilon_0$ be  as in the definition of $\varepsilon$-caps in Section 2, and $C$ be a positive constant.  A space-time point point $(x,t)$ in a surgical solution $g(\cdot)$ to the $n$-dimensional Ricci flow is said to
have an $(\varepsilon, C)$-canonical neighborhood if either  the point $x$ is contained in a  compact component which is strictly PIC2  at time $t$, or there is   an open
neighborhood $U$ of $x$ satisfying $\overline{B(x,t,C^{-1}R(x,t)^{-\frac{1}{2}})} \subset U\subset
B(x,t,CR(x,t)^{-\frac{1}{2}})$ and with  the values of the scalar curvature in $U$ at time $t$  lying (strictly) between $C^{-1}R(x,t)$ and $CR(x,t)$, which falls into one of the following two
types:

(a) $U$ is a strong  $\varepsilon$-neck diffeomorphic to $\mathbb{S}^{n-1}/\Gamma \times (-1,1)$  centered at $(x,t)$, or

(b) $U$ is an $\varepsilon$-cap centered at $(x,t)$ with any cross-section of the $\varepsilon$-neck contained in it diffeomorphic to $\mathbb{S}^{n-1}/\Gamma$

\noindent for some finite subgroup $\Gamma$ of $O(n)$ acting freely on $\mathbb{S}^{n-1}$,
 and
 \begin{equation*}
 \text{vol}_{g(t)}(U) > (C|\Gamma|)^{-1} R(x,t)^{-n/2};
 \end{equation*}
 moreover, the scalar curvature in $U$ at time $t$  satisfies
the derivative estimates
\begin{equation*}
|\nabla R|< C R^{\frac{3}{2}} \hspace*{8mm} \mbox{and} \hspace*{8mm}
|\frac{\partial R}{\partial t}|< C R^2.
\end{equation*}

 From now on we fix  $\varepsilon_0=\min \{\varepsilon_3, \varepsilon_4\}$, where the constant $\varepsilon_3$ is as in Proposition \ref{gluing3},  and the constant $\varepsilon_4$ is as in Proposition \ref{prop 3.5}.  Choose $0 < \varepsilon <  \frac{1}{4} \min \{\frac{1}{10^4n}, \tilde{\varepsilon}_1(\varepsilon_0)\}$, where $\tilde{\varepsilon}_1(\cdot)$ is as in Lemma \ref{gluing2},  and $\theta>0$. Let $\beta=\beta(\varepsilon)$ be the constant given in Lemma \ref{lem 2.5}, and choose
 \begin{equation*}
 C = \max \{100 \varepsilon^{-1}, C_1(n,\theta,\varepsilon), C_2(n,\theta,\varepsilon), C_3(n,\theta), C_1'(n,\beta\varepsilon), C_2'(n,\beta\varepsilon), C_3'(n)\},
 \end{equation*}
 where  the constants on the RHS are from  Propositions \ref{prop 3.9}  and \ref{prop 3.10}.   We have the following definition.

\vspace *{0.2cm}

\noindent {\bf Canonical neighborhood assumption}:  Fix $\varepsilon$ and $C$ as above. Let $r>0$. An evolving Riemannian $n$-orbifold $\{(\mathcal{O}(t), g(t))\}_{t \in I}$
 satisfies the canonical neighborhood assumption  $(CN)_r$ with $(4\varepsilon, 4C)$-control if  any  space-time point $(x,t)$ with  $R(x,t)\geq
r^{-2}$ has a  $(4\varepsilon, 4C)$-canonical neighborhood.

\vspace *{0.2cm}

\noindent {\bf Remark}. The pinching assumption and the canonical neighborhood assumption are the so called a priori assumptions for the surgical solutions to the Ricci flow, which we need to justify.

\vspace *{0.2cm}

The so called bounded curvature at bounded distance property is one of the key ingredients in
Perelman's work \cite{P1}, \cite{P2}; compare Theorem 10.2 in \cite{MT},  Theorem
6.1.1 in \cite{BBB+} and Theorem 6.4 in \cite{BBM}. 4-dimensional versions have appeared
in \cite{CTZ}, \cite{CZ} (though not  displayed as a separate proposition there), \cite{Hu13} and \cite{Hu15}. For the higher dimensional case  see \cite{B19} (though not  displayed as a separate proposition there).

\begin{prop} \label{prop 4.4} (cf. \cite{BBB+}) \ \  Let $n\geq 5$, $0<\varepsilon' <  \varepsilon_1$, where $\varepsilon_1$ is from Proposition \ref{gluing1}, $f$ and $\theta$ in the pinching assumption be fixed,  and $A$ and $ C'$ be two positive numbers.  Then there  exist numbers $Q>0$ and $\Lambda>0$ with the following
property. Assume that  $(\mathcal{O}(t), g(t))$, $t\in [a,b]$, is a complete surgical solution (with at most isolated orbifold singularities) to the $n$-dimensional Ricci flow with
 positive isotropic curvature  and  $(f,\theta)$-pinched, bounded curvature.
Let $(x_0, t_0)$ be a space-time point (in $\mathcal{O}(t_0)$) which is not in a compact component with strictly PIC2, such that:

1. $R(x_0, t_0)\geq Q$;

2. For each point $y\in B(x_0, t_0, AR(x_0, t_0)^{-1/2})$, if $R(y,
t_0)\geq 4R(x_0, t_0)$, then $(y, t_0)$ has an $(\varepsilon',
C')$-canonical neighborhood.

\noindent Then for any $y\in B(x_0, t_0, AR(x_0, t_0)^{-1/2})$, we
have
\begin{equation*}
\frac{R(y, t_0)}{R(x_0, t_0)}\leq \Lambda.
\end{equation*}
\end{prop}

\noindent {\bf  Proof} \ \ We will adapt the proof of Theorem 6.1.1 in \cite{BBB+}  to our situation, incorporating  arguments  in the proof of Propositions 4.2 and 4.4 in
\cite{CTZ}  and Step 3 in the proof of Theorem 10.10 in  \cite{B19}. (Compare the proof of Proposition 4.1 in \cite{Hu15}.) We argue by contradiction. Suppose the result is not true.
Then there exist    sequences $Q_k \rightarrow +\infty$ and $\Lambda_k \rightarrow +\infty$,  a sequence of complete surgical solutions (with at most isolated orbifold singularities) $(\mathcal{O}_k(t), g_k(t))$ to the $n$-dimensional Ricci flow
with  positive isotropic curvature  and  $(f,\theta)$-pinched, bounded curvature, and a sequence of space-time points $(x_k, t_k)$  in  $\mathcal{O}_k(t_k)$  which are not in  components with strictly PIC2,  such that:

1. $R(x_k, t_k)\geq Q_k$;

2. for each point $y\in B(x_k, t_k, AR(x_k, t_k)^{-1/2})$, if $R(y,
t_k)\geq 4R(x_k, t_k)$, then $(y, t_k)$ has an $(\varepsilon',
C')$-canonical neighborhood;

3. for each $k$, there exists $z_k\in B(x_k, t_k, AR(x_k, t_k)^{-1/2})$ with

\begin{equation*}
\frac{R(z_k, t_k)}{R(x_k, t_k)}> \Lambda_k.
\end{equation*}

\noindent For each $k$, consider the parabolic rescaling
\begin{equation*}
\bar{g}_k(\cdot):=R(x_k,t_k)g_k(t_k+\frac{\cdot}{R(x_k, t_k)}).
\end{equation*}
We will adopt the convention in \cite{BBB+} and
 \cite{BBM}  to put a bar on the points when the relevant geometric quantities are computed w.r.t. the metric $\bar{g}_k$.

Define
\begin{equation*}
\rho:=\sup \{s>0 ~|~ \exists C(s)>0, \forall k\in \mathbb{N}, \forall\bar{y}\in B(\bar{x}_k,0,s), R(\bar{y},0)\leq C(s)\}.
\end{equation*}
Then $0<\rho \leq A$. It is easy to see that there exists, up to an extraction, a sequence of points $\bar{y}_k \in B(\bar{x}_k,0,\rho)$ such that
\begin{equation*}
R(\bar{y}_k,0)\rightarrow +\infty \hspace*{8mm}  \mbox{and} \hspace*{8mm}   d_0(\bar{x}_k,\bar{y}_k)\rightarrow \rho \hspace*{8mm}  \mbox{as} \hspace*{8mm} k\rightarrow \infty.
\end{equation*}

As in \cite{BBB+}, we choose points $\bar{x}_k'$ and $\bar{y}_k'$ in  the (minimizing) geodesic segment $[\bar{x}_k\bar{y}_k]$ (w.r.t. $\bar{g}_k(0)$) for large $k$ such
that $R(\bar{x}_k',0)=2C'$, $R(\bar{y}_k',0)=R(\bar{y}_k,0)/(2C')$,
and $[\bar{x}_k'\bar{y}_k'] \subset [\bar{x}_k\bar{y}_k]$ is a
maximal subsegment on which
\begin{equation*}
2C'\leq R(\cdot,0)\leq \frac{R(\bar{y}_k,0)}{2C'},
\end{equation*}
with $\bar{x}_k'$ closest to  $\bar{x}_k$. Note that $d_0(\bar{x}_k, \bar{y}_k')\rightarrow \rho$ as $k\rightarrow \infty$.
Each point $\bar{z} \in
[\bar{x}_k'\bar{y}_k']$ has an $(\varepsilon', C')$-canonical
neighborhood, denoted by $U(\bar{z})$, which contains neither $\bar{x}_k$ nor $\bar{y}_k$.
Using the minimizing property of the geodesic segment $[\bar{x}_k\bar{y}_k]$  we see  that $U(\bar{z})$ can not be an $\varepsilon'$-cap, so it must be a  strong $\varepsilon'$-neck. For any $k$ we have $d_0(\bar{x}_k, \bar{x}_k')\geq (\varepsilon' \sqrt{2C'})^{-1}$. Moreover,  for $k$ sufficiently large we have $d_0(\bar{x}_k', \bar{y}_k')\geq (\varepsilon' \sqrt{2C'})^{-1}$, and $d_0(\bar{x}_k, \bar{x}_k')\leq \rho-(\varepsilon' \sqrt{2C'})^{-1}$.  Up to an extraction we may assume that $\lim_{k\rightarrow \infty}d_0(\bar{x}_k, \bar{x}_k')=a$, where $0<a<\rho$. Let $\bar{z}_k$ be the point on the geodesic segment $[\bar{x}_k\bar{y}_k]$  with $d_0(\bar{x}_k, \bar{z}_k)=\frac{a+\rho}{2}$. Let $U_k$ be the union of
 $U(\bar{z})$'s for all $\bar{z} \in
[\bar{x}_k'\bar{y}_k']$, and $(\widetilde{U}_k, \pi_k)$ be the universal cover of $U_k$.
 Let $\bar{\gamma}_k=[\bar{z}_k\bar{y}_k] \cap U_k$.

   We pull back
the rescaled  solutions $(U_k, \bar{g}_k(\cdot))$ to $\widetilde{U}_k$ via the covering maps $\pi_k: \widetilde{U}_k \rightarrow U_k$.  Choose a point $\tilde{z}_k \in \pi_k^{-1}(\bar{z}_k)$, and let $\tilde{\gamma}_k$ be the component of $\pi_k^{-1}(\bar{\gamma}_k)$ which contains the point $\tilde{z}_k $. By using a local compactness theorem for the Ricci flow (see Theorem C.3.3 in \cite{BBB+}, cf. also Theorem 4.1.5 in \cite{CaZ}, Theorem 16.1 in \cite{H95b}, and Appendix E in \cite{KL08})  we see that up to an extraction, the sequence of pointed
 solutions  $(B(\tilde{z}_k,0,\frac{\rho-a}{2}), \pi_k^*\bar{g}_k(\cdot), (\tilde{z}_k,0))$ (where $B(\tilde{z}_k,0,\frac{\rho-a}{2})=B_{\pi_k^*\bar{g}_k(0)}( \tilde{z}_k, \frac{\rho-a}{2}))$ converges smoothly to a partial Ricci flow  $(B(z_\infty,0,\frac{\rho-a}{2}),  g_\infty(\cdot), (z_\infty,0))$, and the sequence of geodesics $\tilde{\gamma}_k$ converges to a (minimal) geodesic $\gamma_\infty$
 (w.r.t. $g_\infty(0)$) starting from $z_\infty$. The geodesic $\gamma_\infty$ is covered by strong $2\varepsilon'$-necks, and the scalar curvature goes to infinity as a point approaches the open end of the image of $\gamma_\infty$. By Lemma \ref{gluing1.5},  $ B(z_\infty,0,\frac{\rho-a}{2})$ contains a $2\varepsilon'$-horn, say $H$, diffeomorphic to $\mathbb{S}^{n-1}\times (0,1)$, which contains the open end of the image of $\gamma_\infty$. For any $z\in B(z_\infty,0,\frac{\rho-a}{2})$, the evolving metric $g_\infty(z,\cdot)$ exists at least on the time interval $[-\frac{R(z,0)^{-1}}{2},0]$. In addition, by using the pinching assumption we see that $(B(z_\infty,0,\frac{\rho-a}{2}), g_\infty(\cdot))$ is weakly PIC2 and strictly PIC.

 Let $(\overline{B}_\infty, \hat{d})$ be the completion of the metric space $(B(z_\infty,0,\frac{\rho-a}{2}),  d_\infty)$, where the metric  $d_\infty$ is induced by the Riemannian metric $g_\infty(0)$. Let $y_\infty$ be the limit point of
 $\gamma_\infty$ which lies in $\overline{B}_\infty\setminus B(z_\infty,0,\frac{\rho-a}{2})$, and $\widehat{U}_\infty=H \cup \{y_\infty\}$. The pointed metric space $(\widehat{U}_\infty, \hat{d}, y_\infty)$ is a locally complete Alexandrov space with nonnegative curvature, and there exists an $n$-dimensional tangent cone of  $(\widehat{U}_\infty, \hat{d})$ at $y_\infty$ with small aperture, denoted by $C_\infty$, which is a nonnegatively curved metric cone. Denote the vertex of the cone $C_\infty$ by 0. Then a sequence of  suitable rescalings of certain pieces of $(H,  g_\infty(\cdot))$ converges smoothly to a Ricci flow of weakly PIC2 and strict PIC whose final time slice lies in  $C_\infty\setminus \{0\}$. (For more details, see  Step 3 in the proof of  Theorem 6.1.1 in \cite{BBB+} and Step 2 in the proof of Theorem 4.1 in \cite{CZ}; compare also \cite{B19}, \cite{KL} and \cite{P1}.) This contradicts   Proposition 6.5 in \cite{B19}.   \hfill{$\Box$}

\vspace *{0.2cm}

The following proposition is similar to Theorem 6.2.1 in \cite{BBB+}, Theorem 6.5 in \cite{BBM},  Proposition 2.3 in \cite{Hu13},  and
Proposition 4.2 in \cite{Hu15}; compare Proposition 11.1 in \cite{B19},  Propositions 4.2 and 4.4 in \cite{CTZ},  and Lemma 4.3 in \cite{P2}.

\begin{prop} \label{prop 4.5}(cf. \cite{BBM}) \ \  Let $n\geq 5$,  $f$ and $\theta$ in the pinching assumption be fixed.  Let  $\varepsilon>0$  and $C>0$ be chosen as  above. Given $r, \delta>0$,
there exist  numbers  $h \in (0, \delta r)$ and $D> 10$, such that if
$(\mathcal{O}(\cdot),g(\cdot))$ is a complete surgical solution to the Ricci flow in dimension $n$ with
 positive isotropic curvature and with bounded curvature, defined on a time interval
$[0, b]$  and satisfying the $(f,\theta)$-pinching
assumption  and the canonical neighborhood assumption $(CN)_r$ with  $(4\varepsilon, 4C)$-control, then
the following holds:

 \noindent Let $t \in (0,b]$ and  $x, y, z \in \mathcal{O}(t)$ such that $R(x,t) \leq
2/r^2$, $R(y,t)=h^{-2}$  and $R(z,t)\geq D/h^2$. Suppose that there is a  curve $\gamma$  in $\mathcal{O}(t)$  connecting $x$ to $z$ and containing $y$, such that each point in $\gamma$
of scalar curvature in $[8Cr^{-2}, (4C)^{-1}Dh^{-2}]$ is a
center of a $4\varepsilon$-neck. Then
$(y,t)$ is a center of a strong $\delta$-neck.
\end{prop}

\noindent {\bf Proof}.\ \ We essentially follow the proof of Theorem 6.2.1 in \cite{BBB+}, Theorem 6.5 in \cite{BBM},  and Proposition 11.1 in \cite{B19} with some necessary modifications. (Compare the proof of
Propositions 4.2 and 4.4 in \cite{CTZ}, and Proposition 4.2 in \cite{Hu15}.) We argue by contradiction. Otherwise, there exist numbers
$r, \delta>0$,  sequences $h_k \rightarrow 0$, $D_k\rightarrow
+\infty$, a sequence of complete surgical solutions
$(\mathcal{O}_k(\cdot),g_k(\cdot))$ in dimension $n$ with
positive isotropic curvature and bounded curvature and
satisfying the  $(f,\theta)$-pinching assumption and $(CN)_r$ with  $(4\varepsilon, 4C)$-control, and sequences $ t_k>0$, $x_k, y_k, z_k \in \mathcal{O}_k(t_k)$ with $R(x_k,t_k)\leq 2r^{-2}$, $R(y_k,t_k)=h_k^{-2}$
and $R(z_k,t_k)\geq D_kh_k^{-2}$, and finally a sequence of curves
$\gamma_k$ in $\mathcal{O}_k(t_k)$ connecting $x_k$ to $z_k$ and containing $y_k$, whose points
of scalar curvature in $[8Cr^{-2}, (4C)^{-1}D_kh_k^{-2}]$ are
centers of $4\varepsilon$-necks, but  $(y_k,t_k)$  is not  a center of a strong
$\delta$-neck.

Let $\bar{g}_k(t)=h_k^{-2}g_k(t_k+h_k^2t)$ for each $k$. For any $\rho>0$,  as in the proof of Theorem 6.5 in \cite{BBM}, using Proposition \ref{prop 4.4} we see that when $k$ is sufficiently large, $x_k \notin B_{\bar{g}_k(0)}(y_k,\rho)$, $z_k \notin B_{\bar{g}_k(0)}(y_k,\rho)$, and   $B_{\bar{g}_k(0)}(y_k,\rho)$ is contained in the union of some 4$\varepsilon$-necks, which must be strong 4$\varepsilon$-necks by our canonical neighborhood assumption.  Now let $U_k \subset \mathcal{O}_k(t_k)$ be the maximal $4\varepsilon$-tube containing  $y_k$ (cf. Proposition \ref{gluing1} and Lemma \ref{lem 4.1}), in the sense that $U_k $ cannot be properly included  in  any  $4\varepsilon$-tube in $\mathcal{O}_k(t_k)$ containing $y_k$.  Let $(\widetilde{U}_k, \pi_k)$ be the universal cover of $U_k$.
  Choose a point $\tilde{y}_k\in \pi_k^{-1}(\tilde{y}_k)$.
  Then we pull back the parabolically rescaled solutions  $(U_k, \bar{g}_k(t))$ to  $\widetilde{U}_k$
 via the covering maps $\pi_k$. As in the proof of Theorem 6.2.1 in \cite{BBB+}  and Proposition 11.1 in \cite{B19}, using Proposition \ref{prop 4.4} above and Theorem 6.4 in \cite{B19} we can show that the sequence of pointed pulled-back solutions $(\widetilde{U}_k, \pi_k^*\bar{g}_k(t),(\tilde{y}_k,0))$ subconverges smoothly to a pointed standard cylindrical flow on $\mathbb{S}^{n-1}\times \mathbb{R}$ in the time interval $(-\infty,0]$. So given any $\tilde{\delta}>0$, when $k$ is sufficiently large, $(\tilde{y}_k,0)$  is a center of a  $\tilde{\delta}$-neck with Hamilton's canonical parametrization, which can be pushed down to give a  $\tilde{\delta}$-neck centered at $(y_k,t_k)$. (Compare  Corollary C2.3 in \cite{H97} and the proof of Lemma \ref{gluing2}.) By Lemma \ref{lem 2.5}, $(y_k,t_k)$  is a center of a strong $\delta$-neck  for $k$  large enough. A contradiction.
\hfill{$\Box$}

\vspace *{0.2cm}

Given the surgery parameters $r,
\delta >0$, let $h:=h(r,\delta), D:=D(r,\delta)$, and $ \Theta:=2Dh^{-2}$  be the associated cutoff parameters (following the convention in \cite{BBM}) as
determined in Proposition \ref{prop 4.5}, which also depend on $\varepsilon$, $n$, $f$,  and $\theta$. The parameter $ \Theta$ will serve as  the curvature threshold for the surgery process (as in \cite{BBM}, \cite{Hu13} and \cite{Hu15}), that is, we'll do
surgery on $\delta$-necks when the  supremum $R_{\text{max}}(t)$ of the scalar curvature at time $t$ reaches $\Theta$.

\vspace *{0.2cm}

Now we describe Hamilton's metric surgery on a $\delta$-neck $N$ (with $\delta$ small) in a Riemannian orbifold $(\mathcal{O},g)$ with a diffeomorphism $\psi:\mathbb{S}^{n-1}/\Gamma \times (-\frac{1}{\delta}, \frac{1}{\delta}) \rightarrow  N$ which is a $\delta$-homothety with scaling factor $Q$. Roughly speaking, the surgery is removing the open subset $\psi(\mathbb{S}^{n-1}/\Gamma \times (-(10\delta)^{-1}, (10\delta)^{-1}))$ from $\mathcal{O}$ and gluing in a Riemannian orbifold  (called a surgery cap) diffeomorphic to $D^n//\Gamma$ (which has an isolated singular point when $\Gamma$ is nontrivial) along each boundary component thus produced. More precisely,  the gluing at the boundary component $\psi(\mathbb{S}^{n-1}/\Gamma \times \{-(10\delta)^{-1}\})$  of  $\mathcal{O} \setminus \psi(\mathbb{S}^{n-1}/\Gamma \times (-(10\delta)^{-1}, (10\delta)^{-1}))$ is as follows; the gluing at $\psi(\mathbb{S}^{n-1}/\Gamma \times \{(10\delta)^{-1}\})$ is similar. As in the appendix and the proof of Proposition \ref{prop 3.5}, we choose a smooth function $\varphi: (-(10\delta)^{-1}, b_0-(10\delta)^{-1}]\rightarrow \mathbb{R}$, where $b_0 > \frac{1}{10}$, (in particular,
  $\varphi(z)= ce^{-1/(z+(10\delta)^{-1})}$ when $z \in (-(10\delta)^{-1}, \frac{1}{10}-(10\delta)^{-1}]$, where  $0< c < \frac{1}{\sqrt{(n-1)(n-2)}}$ is a small positive constant,) and a smooth cutoff function $\chi: (\frac{1}{20}-(10\delta)^{-1}, b_0-(10\delta)^{-1}]\rightarrow [0,1]$   with $\chi(z)=1$ for $\frac{1}{20}-(10\delta)^{-1} < z \leq \frac{1}{18}-(10\delta)^{-1}$ and $\chi(z)=0$ for $\frac{1}{12}-(10\delta)^{-1} \leq z \leq b_0-(10\delta)^{-1}$, and  deform the metric $g$ on   $\psi(\mathbb{S}^{n-1}/\Gamma \times (-(10\delta)^{-1},  b_0-(10\delta)^{-1}])$ to a new metric $g'$, with
\begin{equation*}
 g'(x)= \begin{cases}
    g(x), & {z\in (-\frac{1}{\delta}, -\frac{1}{10\delta}];}  \\
    e^{-2\varphi(z)}g(x), & {z\in (-\frac{1}{10\delta},\frac{1}{20}-\frac{1}{10\delta}];} \\
    e^{-2\varphi(z)}(\chi(z)g(x)+(1-\chi(z))Q^{-1}(\psi_*\bar{g}_{\text{cyl}})(x)), & {z\in (\frac{1}{20}-\frac{1}{10\delta},  b_0-\frac{1}{10\delta}],}
     \end{cases}
\end{equation*}
 where  for $x\in \psi(\mathbb{S}^{n-1}/\Gamma \times (-\frac{1}{\delta}, b_0-\frac{1}{10\delta}])$ we denote  $z(x)= p_2\circ \psi^{-1}(x)$ by $z$,  here $p_2: \mathbb{S}^{n-1}/\Gamma \times (-\delta^{-1},  \delta^{-1}) \rightarrow (-\delta^{-1},  \delta^{-1})$ is the projection onto the second factor, and $\bar{g}_{\text{cyl}}$ is the standard metric on the (quotient) cylinder $\mathbb{S}^{n-1}/\Gamma \times (-\delta^{-1},  \delta^{-1})$ with scalar curvature 1. (We require that $\frac{1}{\delta} \gg b_0$.) In particular, the function $\varphi$ is chosen such that the metric $g'$ restricted to the portion $\psi(\mathbb{S}^{n-1}/\Gamma \times [b'-(10\delta)^{-1}, b_0-(10\delta)^{-1}])$ has constant positive sectional curvature, here $\frac{1}{10} < b' < b_0$  as in the appendix. Then we smoothly attach a  quotient of a closed metric $n$-ball of constant positive sectional curvature  by $\Gamma$ to $\mathcal{O} \setminus \psi(\mathbb{S}^{n-1}/\Gamma \times (b_0-(10\delta)^{-1}, (10\delta)^{-1}-b_0))$ (with the new metric) along the boundary component  $\psi(\mathbb{S}^{n-1}/\Gamma \times  \{b_0-(10\delta)^{-1}\})$, as we do in the construction of the orbifold standard solution in the appendix.  Note that there is a positive number $\delta'=\delta'(\delta)$  with $\lim_{\delta \rightarrow 0} \delta'=0$ such that the union of $\psi(\mathbb{S}^{n-1}/\Gamma \times (-\delta^{-1}, -(10\delta)^{-1}])$ and the surgery cap with the metric just described  is $\delta'$-homothetic to the union of $\mathbb{S}^{n-1}/\Gamma \times (-\frac{9}{10\delta}, 0]$ and the surgery cap with the metric $\hat{g}_\Gamma(0)$  defined in the appendix. (Here, $\delta'$ also depends on the choice of the functions $\varphi$ and $\chi$ above.)

\vspace *{0.2cm}

We adapt two more definitions from \cite{BBM}.

\vspace *{0.2cm}

\noindent {\bf Definition} (compare \cite{BBM}, \cite{Hu15}). \ \ Given an interval
$I\subset [0,+\infty)$, fix surgery parameters $r$, $\delta>0$
and let
$h$, $D$, $\Theta=2Dh^{-2}$ be the associated cutoff parameters. Let
$(\mathcal{O}(t),g(t))$ ($t \in I$) be an evolving complete Riemannian orbifold with  at most isolated singularities. Let $t_0 \in I$ and $(\mathcal{O}_+,g_+)$ be a
(possibly empty) Riemmanian $n$-orbifold. We say that
$(\mathcal{O}_+,g_+)$ is obtained from
$(\mathcal{O}(\cdot),g(\cdot))$ by $(r,\delta)$-surgery at time
$t_0$ if

i. $R_{\text{max}}(g(t_0))=\Theta$, and there is a
collection    of pairwise disjoint strong $\delta$-necks in $\mathcal{O}(t_0)$, centered at some points with scalar curvature equal to $h^{-2}$ at time $t_0$, such that   $\mathcal{O}_+$ is obtained from $\mathcal{O}(t_0)$ by doing
  Hamilton's surgery along these necks as described above, and removing each of the following components:

  (a)  a component which is strictly PIC2,

  (b) a component  where each point  is a center of a
$4\varepsilon$-neck or a $4\varepsilon$-cap;

ii. $R_{\text{max}}(g_+)\leq
\Theta/2$ when $\mathcal{O}_+\neq \emptyset$.

\vspace *{0.2cm}

\noindent {\bf Definition} (cf. \cite{BBM} and \cite{Hu15}). \ \  A surgical solution
$(\mathcal{O}(\cdot),g(\cdot))$ to the Ricci flow defined on some time interval
$I\subset [0,+\infty)$  is an $(r,\delta)$-surgical solution  if it has the  following
properties:

i.  It satisfies the $(f,\theta)$-pinching assumption  at any time $t\in I$, and $R(x,t) \leq \Theta$ for all $(x,t)$;

ii. At each singular time $t_0\in I$,
$(\mathcal{O}_+(t_0),g_+(t_0))$ is obtained from
$(\mathcal{O}(\cdot),g(\cdot))$ by $(r,\delta)$-surgery at time
$t_0$;

iii. The canonical neighborhood assumption  $(CN)_r$ with $(4\varepsilon, 4C)$-control holds.

\vspace *{0.2cm}

The following proposition  is
analogous to Proposition A in \cite{BBB+}, Proposition A in \cite{BBM},  Proposition 2.7 in \cite{Hu13}  and Proposition 4.4 in \cite{Hu15}.

\begin{prop} \label{prop 4.6}\ \  Let $n\geq 12$, the family $\mathcal{F}_t$, $t\in [0,\hat{T}]$, be as above, and $f$ and $\theta$ be as in (3.1). Let $\varepsilon$  and $C$ be chosen as before. Let $\delta$ be sufficiently small such that $\delta \ll \varepsilon$ and (an orbifold version of) Proposition 8.2 in \cite{B19} holds.
Let $(\mathcal{O}(t), g(t))$, $t\in [0,b]$,
be a compact $(r, \delta)$-surgical solution to the Ricci flow of dimension $n$ (with $(f,\theta)$-pinched curvature), where $b<\hat{T}$.  Suppose
that the curvature tensor  at time $b$ lies in  $\mathcal{F}_b$, and $R_{max}(b)=\Theta$, which is defined above.
 Then there exists a Riemannian orbifold
$(\mathcal{O}_+,g_+)$ which is obtained from $(\mathcal{O}(\cdot),g(\cdot))$ by
$(r,\delta)$-surgery at time $b$, such that the curvature tensor of  $g_+$ still lies in the set $\mathcal{F}_b$,
and $R_{min}(g_+)\geq R_{min}(g(b))$.
\end{prop}

\noindent {\bf Proof}.\ \  Using Proposition \ref{prop 4.5} and an orbifold version of Proposition 8.2 in \cite{B19},  the proof is similar to that of Proposition A in \cite{BBB+}, Proposition A in \cite{BBM} and Proposition 2.7 in \cite{Hu13}. First we remove from $\mathcal{O}(b)$ all components  which are strictly PIC2 at time $b$, and  get a Riemannian orbifold, denoted by $\mathcal{O}'$.  We divide  $\mathcal{O}'$ into three parts $X$, $Y$ and $Z$, such that the scalar curvature at each point in $X$ is less than $2r^{-2}$, the scalar curvature at each point in $Y$ lies  in the interval $[2r^{-2}, \frac{\Theta}{2})$, and the scalar curvature at each point in $Z$ is not less than $\frac{\Theta}{2}$.  Note that there exists a maximal collection $\{N_i\}$ of pairwise disjoint cutoff necks in $\mathcal{O}'$ (a cutoff neck in $\mathcal{O}'$ is  a strong $\delta$-neck in $\mathcal{O}'$ centered at some point   of scalar
curvature $h^{-2}$).

\vspace *{0.2cm}

\noindent {\bf Claim}. Any connected component of $\mathcal{O}'\setminus \cup_iN_i$ is contained in either  $X\cup Y$ or  $Y\cup Z$.

\vspace *{0.2cm}

\noindent {\bf Proof of Claim}.  Otherwise there is some component $W$ of $\mathcal{O}'\setminus \cup_iN_i$ containing  a point $x\in X$ and a point $z\in Z$.
We can choose  a shortest path $\gamma$ in $W$ connecting $x$ and $z$ w.r.t. the distance induced by the metric $g(b)$; cf. for example Theorem 2.5.23 in \cite{BBI}.  It is not hard to see that the intersection of $\gamma$ with $\partial W$ is empty. For example, one can argue as the following: Suppose the intersection of $\gamma$ with $\partial W$ is not empty. Note that any point $y\in \gamma \cap \partial W$ has a canonical neighborhood $U$ in $\mathcal{O}'$. If $U$ is a $4\varepsilon$-cap, by the  curvature assumption on $x$  and $z$, the fact that $\frac{1}{2}h^{-2} \leq R(y,b) \leq \frac{3}{2}h^{-2}$,  and the scalar curvature bounds in the definition of $(4\varepsilon,4C)$-canonical neighborhood of the cap type, neither $x$ nor $z$ is contained in $U$. Write $U=N\cup C$, where $N$
is a $4\varepsilon$-neck in the end of the cap $U$, $N\cap C=\emptyset$, $\overline{N}\cap
C=\partial C$ and $y\in \text{Int} \hspace*{0.5mm} C$. Let $N_{i_0}$ be the $\delta$-neck whose closure contains the point $y$. We may assume that the point $y$ is near the negative end of the neck $N_{i_0}$  (w.r.t. the height function). First observe that $C \cap N_{i_0} \neq \emptyset$. Otherwise any neighborhood of $y$  (in $\mathcal{O}'$) must contain a point in $N_{i_0}$, which is not in $C$. This is impossible, as $\text{Int} \hspace*{0.5mm} C$ is an open subset of $\mathcal{O}'$.  Then it is easy to see that the supremum of the height function on the neck $N_{i_0}$ restricted to $C \cap N_{i_0}$ is attainable, and is attained at a point in $\partial C$ instead of $\text{Int} \hspace*{0.5mm} C$ by an argument as before using again the fact that $\text{Int} \hspace*{0.5mm} C$ is an open subset of $\mathcal{O}'$. In particular, $\partial C \cap N_{i_0} \neq \emptyset$. Note that $\partial C$ is a central cross-section of a $4\varepsilon$-neck by definition of $4\varepsilon$-cap, and $\delta \ll \varepsilon$. From the proof of Proposition \ref{gluing1} we see that the central cross-section $S$ of the neck $N$ is contained in $N_{i_0}$. But this is impossible, since on the one hand, as $S$ separates $\mathcal{O}'$ (i.e. $\mathcal{O}'\setminus S$ is not connected),  the path $\gamma$ must intersect $S$; on the other hand,  $\gamma$ can not intersect $N_{i_0}$, in particular, $\gamma$ can not intersect $S$. (Alternatively, one can also argue as follows. Observe that for any $i$, the middle nine-tenths of the neck $N$  can not intersect $N_i$, for otherwise suppose that the middle nine-tenths of the neck $N$ intersects some $N_i$, then by
the proof of Proposition \ref{gluing1} there would exist a cross-section of $N$, say $S'$, which is contained in $N_i$. But this is impossible as seen above. So the middle nine-tenths of the neck $N$ is contained in $W$. But this is also impossible, as the path $\gamma$ must intersect the central cross-section of the neck $N$ at least twice, $\gamma$ can not be minimizing in the situation that the middle nine-tenths of the neck $N$ is contained in $W$,  which is a contradiction.)  If $U$ is a $4\varepsilon$-neck, denote by $A$ the boundary component of $W$ which contains $y$, then $\text{diam}~A \ll \text{diam}~ U$. Note that $\gamma$ would enter and exit $U \cap W$  from the same end of $U$.  Again this contradicts the minimizing property of $\gamma$ in $W$.  Compare the argument for the first  claim in the proof of
Lemma 7.7 in \cite{BBM}.

By the definition of $(r,\delta)$-surgical solution, at time $b$ any point $y$ in $\gamma$ of scalar curvature in $[8Cr^{-2}, (4C)^{-1}Dh^{-2}]$ has a $(4\varepsilon, 4C)$-canonical neighborhood $U$ (in $\mathcal{O}'$). $U$ cannot be a component of strictly PIC2, because we have thrown
away all such components in $\mathcal{O}(b)$.
 We will show that $U$ cannot be a
$4\varepsilon$-cap either.    Otherwise similarly as above write $U=N\cup C$, where $N$
is a $4\varepsilon$-neck, $N\cap C=\emptyset$, $\overline{N}\cap
C=\partial C$ and $y\in \text{Int} \hspace*{0.5mm} C$. Note again that neither $x$ nor $z$ is contained in $U$. Let $\psi: \mathbb{S}^{n-1}/\Gamma \times
(-(4\varepsilon)^{-1}, (4\varepsilon)^{-1}) \rightarrow N$ be the
diffeomorphism associated to  the neck $N$, and
$S=\psi( \mathbb{S}^{n-1}/\Gamma \times \{0\})$. Assume that as the ``height" $s\rightarrow (4\varepsilon)^{-1}$  the point in $N$ approaches the end of $U$.  Clearly $\gamma$ is
not minimizing in $U$, since if $x'$ (resp. $z'$) is an intersection
of $\gamma$ with $S$ between $x$ and $y$ (resp. $y$ and $z$), then
$d(x',z')\ll d(x',y)+d(y,z')$. The geodesic segment (in $U$)
$[x'z']$ is not contained in $W$ by the minimality of $\gamma$ in
$W$. So $ [x'z']\cap
\partial W \neq \emptyset$. By definition of $W$, the corresponding component of
$\partial W$ is a boundary component, denoted by $S_i^+$, of the closure of some
cutoff neck $N_i$. Then $d(S_i^+,S)< \text{diam}(S)$ since $ [x'z']\cap
S_i^+ \neq \emptyset$. It follows that the ratio of the scalar curvature on the neck $N_i$ near $S_i^+$ to that on the neck $N$ near $S$ is very close to 1 (cf. p. 32 in \cite{BBB+}  and Proposition A.11 in \cite{MT}), and $\text{diam}(S_i^+) \ll \text{diam}(N)$.  Then  we see that $\psi(\mathbb{S}^{n-1}/\Gamma \times
(-(4\varepsilon)^{-1}, -(8\varepsilon)^{-1}))\subset W$.
 This contradicts the minimality of $\gamma$ in $W$.
 Compare the argument for the second claim in the proof of
Lemma 7.7 in \cite{BBM} and  Claim 2 in the proof of Proposition 2.7 in \cite{Hu13}.  (By the way, the inequality in lines 7-8 on p. 1225 in \cite{Hu13} should read
\begin{equation*}
\begin{split}
&d(p', S) \leq  d(p',\partial_+V')+ d(\partial_+V',S_i^+)+d(S_i^+,S)+\text{diam}(\partial_+V')+ \text{diam}(S_i^+)  \\
< & 0.2\varepsilon_1^{-1}+ 0.03\varepsilon_1^{-1}+\text{diam}(S) +\text{diam}(\partial_+V')+ \text{diam}(S_i^+) <0.3\varepsilon_1^{-1}.
\end{split}
\end{equation*}
Note that $\text{diam}(S)$, $\text{diam}(\partial_+V')$, and $\text{diam}(S_i^+)$ are much smaller than $\varepsilon_1^{-1}$ there.)

So any point  in $\gamma$ of scalar curvature in $[8Cr^{-2}, (4C)^{-1}Dh^{-2}]$ is the center of  a $4\varepsilon$-neck.  Then we can use Proposition \ref{prop 4.5} to find a strong $\delta$-neck $N$ centered at some point in $\gamma$ of scalar curvature $h^{-2}$. Note that $N$ is disjoint from any $N_i$, because otherwise $\gamma$ would intersect such $N_i$ with $N_i\cap N\neq \emptyset$ similarly as argued above. This contradicts the maximality of $\{N_i\}$.  \hfill{$\Box$}

\vspace *{0.2cm}

Now we do Hamilton's surgery on $\mathcal{O}'$ along these strong $\delta$-necks $N_i$, and after surgery we throw away all the components  which are covered by $(4\varepsilon, 4C)$-canonical neighborhoods. Then we get a  Riemannian orbifold, denoted by $(\mathcal{O}_+,g_+)$, which has all the desired properties by the Claim above and an orbifold version of Proposition 8.2 in \cite{B19}.   \hfill{$\Box$}

\vspace *{0.2cm}

For the notation in the following lemma  see  Section 2.

\begin{lem} \label{lem 4.7}(cf. Theorem 8.1.3 in \cite{BBB+}, and Proposition 5.2 in \cite{Hu15}) \ \  For any $n\geq 3$, $A>0$, $\Lambda>0$  and $0<T_0<\frac{n-1}{2}$, there exists
$\rho=\rho(n, A,\Lambda, T_0)>A$ with the following
property. Let $(\mathbb{R}^n//\Gamma, \hat{g}_\Gamma(\cdot))$ be an orbifold standard solution,  $U$ be an open subset of $\mathbb{R}^n//\Gamma$  such that the ball $B(p_\Gamma,0,\rho):=B_{\hat{g}_{\Gamma}(0)}(p_\Gamma,\rho)\subset U$ is
relatively compact, and $g(\cdot)$ be a Ricci flow defined on
$U\times [0,T]$   ($T\in (0,T_0]$). Assume that

(i) $||Rm(g(\cdot))||_{0,U\times [0,T],g(\cdot)}\leq \Lambda$, and

(ii) $g(0)$ is $\rho^{-1}$-close to $\hat{g}_\Gamma(0)$ on $B(p_\Gamma, 0, \rho)$.

\noindent Then $g(\cdot)$ is $A^{-1}$-close to $\hat{g}_\Gamma(\cdot)$ on
$B(p_\Gamma,0, A)\times [0,T]$.

\vspace *{0.2cm}

\noindent Here, $||Rm(g(\cdot))||_{0,U\times [0,T],g(\cdot)}:=
\sup_{U \times [0,T]} \{|Rm_{g(t)}(x)|_{g(t)}\}$.
\end{lem}

\noindent {\bf Proof}.\ \ See the proof of Theorem 8.1.3 in \cite{BBB+} and Proposition 5.2 in \cite{Hu15}. We argue by contradiction. Otherwise there exist numbers $n\geq 3$,
$A>0$, $\Lambda >0$, $0<T_0<\frac{n-1}{2}$, a sequence
$\rho_k\rightarrow +\infty$ as $k\rightarrow \infty$, a sequence of open subsets $U_k\subset \mathbb{R}^n//\Gamma_k$, where
$\Gamma_k$ are finite subgroups of $O(n)$ acting freely on $\mathbb{S}^{n-1}$, such that
$B(p_k,0,\rho_k):=B_{\hat{g}_k(0)}(p_k,\rho_k) \subset U_k$ is relatively compact,  where $p_k:=p_{\Gamma_k}$, and  $\hat{g}_k:=\hat{g}_{\Gamma_k}$, and a sequence of Ricci flows $g_k(\cdot)$ defined
on $U_k \times [0,T_k]$ with $0<T_k\leq T_0$,  such that

(i) $|\text{Rm}_{g_k(t)}|_{g_k(t)}\leq \Lambda$ on $U_k \times [0,T_k]$,

(ii) $g_k(0)$ is $\rho_k^{-1}$-close to $\hat{g}_k(0)$ on  $B(p_k,0,\rho_k)$,

\noindent but for some $t_k\in (0,T_k]$, $g_k(t_k)$ is not $A^{-1}$-close to  $\hat{g}_k(t_k)$  on
$B(p_k,0,A)$. We assume that $t_k$ is minimal for this property.

Now we pull back the solutions  $g_k(\cdot)$  to $\pi_k^{-1}(U_k) \subset \mathbb{R}^n$ via the natural projections $\pi_k:=\pi_{\Gamma_k}: \mathbb{R}^n \rightarrow
\mathbb{R}^n//\Gamma_k$.
 We may assume that the sequence $t_k$ subconverges to a number $t_\infty>0$. Then we can argue as in the proof of Theorem 8.1.3 in \cite{BBB+}, using a  version of Shi's derivative
 estimates  (see Theorem 3.29 in \cite{MT}),
  a local version of Hamilton's compactness theorem for the Ricci flow \cite{H95}  (cf. Theorem 4.1.5 in \cite{CaZ}, Theorem 16.1 in \cite{H95b}, and Appendix E in \cite{KL08}) and Chen-Zhu's
  uniqueness theorem  for complete Ricci flows \cite{CZ06}, to get that the sequence $(B(p_0,0,\rho_k), \pi_k^*g_k(t), p_0)$, $t\in [0,t_\infty)$, subconverges to the smooth standard solution $(\mathbb{R}^n, \hat{g}(t), p_0)$ on the time interval $[0,t_\infty)$. As observed  in p. 93 in \cite{BBB+}, since we have assumption (ii), by the argument in Section 2 of \cite{H95},  the diffeomorphisms involved in the subconvergence here can be chosen to be the identity maps.  Using Shi's derivative estimates and Lemma 8.2.1 in \cite{BBB+} we see that for $k$ sufficiently large $\pi_k^*g_k(\cdot)$  is  $A^{-1}$-close to $\hat{g}(\cdot)$ on  $B(p_0,0,A)\times [0,t_k]$.
 Since $\pi_k^*g_k(\cdot)$ is $\Gamma_k$-invariant and the smooth standard solution is $O(n)$-invariant, it follows that for $k$ sufficiently large $g_k(t_k)$ is  $A^{-1}$-close to $\hat{g}_k(t_k)$ on   $B(p_k,0,A)$. This is a contradiction.  \hfill{$\Box$}

\vspace *{0.2cm}

The following result is important in constructing $(r,\delta)$-surgical solutions to the Ricci flow.

\begin{prop} \label{prop 4.8}(cf. Lemma 4.5 in \cite{P2}, Theorem 8.1.2 in \cite{BBB+}, Theorem 8.1 in \cite{BBM}, Lemma 7.3.6 in \cite{CaZ},  Proposition 10.5 in \cite{B19}, Lemma 9.1.1 in \cite{Z11},  and Proposition 5.1 in \cite{Hu15})\ \  Fix  $n\geq 5$, and let $\varepsilon$  and $C$ be chosen as before. For any $A>0$ and $\alpha \in
(0,\frac{n-1}{2})$, there exists
$\bar{\delta}=\bar{\delta}(A,\alpha)>0$ with the following
property. Let  $r$ and $\delta$ be two positive numbers  with  $\delta \leq
\bar{\delta}$, and
$(\mathcal{O}(\cdot),g(\cdot))$ be a complete surgical solution of dimension $n$ with positive isotropic curvature defined on some time interval $[0,b]$, which is an $(r,\delta)$-surgical  solution (and in particular, satisfies $(CN)_r$ with $(4\varepsilon, 4C)$-control) on $[0,b)$. Let $t_0\in (0,b)$ be a singular time. Consider the
restriction of $(\mathcal{O}(\cdot),g(\cdot))$ to $[t_0,b]$. Let
$p\in \mathcal{O}_+(t_0)$ be the tip of some surgery cap of scale
$h=h(r,\delta)$, and let $t_1=\min \{b,t_0+\alpha h^2\}$. Then either

(i) the parabolic neighborhood $P(p,t_0,Ah,t_1-t_0)$ is
unscathed, and is, after scaling with the factor $h^{-2}$ and
shifting time $t_0$ to zero, $A^{-1}$-close to   the corresponding
subset of an orbifold standard solution,
 or

(ii) the assertion (i) holds with $t_1$ replaced by some $t^+ \in
[t_0, t_1)$; moreover  $B(p,t_0,Ah)$ is removed by a surgery
at time $t^+$.
\end{prop}

\noindent {\bf Proof}.\ \ Using Lemma \ref{lem 4.7} and the fact that  $|\text{Rm}|\leq C(n)R$ for any Riemannian orbifold of dimension $n\geq 5$ with positive isotropic curvature, where $C(n)$ is a positive constant depending only on $n$, the proof is almost the same as that  of Theorem 8.1.2 in \cite{BBB+} and Theorem 8.1 in \cite{BBM}. Note that here we do not need to use the pinching assumption in the definition of $(r,\delta)$-surgical  solution.
\hfill{$\Box$}

\vspace *{0.2cm}

Note that the constant $\bar{\delta}$ in Proposition \ref{prop 4.8} is independent of  the diffeomorphism
type of the orbifold standard solution appearing in (i) above.

\vspace *{0.2cm}

Now we can establish the $\kappa$-noncollapsing property under the a priori assumptions.

\begin{prop} \label{prop 4.9}(cf. Proposition 10.9 in \cite{B19}, Lemma 4.5 in \cite{CTZ},  Lemma 5.4 in \cite{Hu15}, Proposition 16.1 in \cite{MT}, and  Lemma 5.2 in \cite{P2}) Let $(\mathcal{O},g_0)$ be a compact Riemannian orbifold of dimension $n\geq 5$ with at most isolated singularities and with positive isotropic curvature.  Let $\varepsilon$  and $C$ be chosen as before.
Then there exists a positive constant $\kappa$ and a positive function $\tilde{\delta}(\cdot)$ with the following property: If we have an $(r,\delta)$-surgical solution $(\mathcal{O}(t), g(t))$, $t\in [0,T]$, to the Ricci flow with $\delta\leq \tilde{\delta}(r)$ starting with $(\mathcal{O},g_0)$, then for any space-time point $(x_0,t_0)$, either the flow is $\kappa$-noncollapsed at $(x_0,t_0)$ on all scales less than $\varepsilon$, or $x_0$ is contained in a component of strictly PIC2 at time $t_0$.
\end{prop}

\noindent {\bf Proof}.\ \ The proof is an  adaption  of that of Proposition 10.9 in \cite{B19};  the necessary modifications are as in the proof of  in  Lemma 4.5 in \cite{CTZ} and Lemma 5.4 in \cite{Hu15},  with the help of Propositions \ref{prop 4.8} and  \ref{gluing1}.

 Fix any $0< r_0 \leq \varepsilon$. Choose
$t_0\in (0,T]$ and $x_0\in \mathcal{O}(t_0)$ such that $x_0$ is not contained in a component of strictly PIC2 at time $t_0$. Assume
$|\text{Rm}(\cdot,\cdot)| \leq r_0^{-2}$ on $P(x_0, t_0,r_0, -r_0^2)$, we
want to bound $\text{vol}_{t_0}(B(x_0,t_0,r_0))/r_0^n$ from below. We
will consider two  cases, where Case 2 is divided further into two subcases.

\vspace *{0.2cm}

Case 1. $r_0\geq \frac{r}{C(n, \theta, \varepsilon)}$, where $C(n, \theta, \varepsilon)$ is a constant depending only on $n, \theta$ and $\varepsilon$, and is to be determined in Case 2. Using Perelman's reduced volume \cite{P1}, \cite{P2}, and Proposition \ref{prop 4.8}, which is utilized to find the desired function $\tilde{\delta}(\cdot)$,
 the argument of Case 3  in the proof of Proposition 10.9  in \cite{B19} can be adapted to our case
without essential changes; cf. also  Step 1 in the proof of Lemma 4.5 in \cite{CTZ}.

\vspace *{0.2cm}

Case 2.  $r_0 < \frac{r}{C(n, \theta, \varepsilon)}$.  By an argument using orbifold Bishop-Gromov theorem (see \cite{B93} and \cite{L}) we may assume w.l.o.g. that there is some point $(x', t') \in \overline{P(x_0, t_0,r_0, -r_0^2)}$ such that $|\text{Rm}(x', t')|=r_0^{-2}$. (Compare Lemma 10.1.2 in  \cite{BBB+}, and p.232-233 in \cite{CZ}.)  As before, note that $|\text{Rm}|\leq C(n)R$ for any Riemannian orbifold of dimension $n\geq 5$ with positive isotropic curvature, where $C(n)$ is a positive constant depending only on $n$. So if $ C(n, \theta, \varepsilon)$ is large compared to $C(n)$,  $(x', t')$ admits a canonical neighborhood which is not a component of strictly PIC2 (by using our assumption on $x_0$). Then using  the derivative estimates in the definition of the canonical neighborhood we can get $R(x_0, t_0) > r^{-2}$ by choosing $C(n, \theta, \varepsilon)$ sufficiently large. In particular, $(x_0,t_0)$ is a center of a $4\varepsilon$-neck or $4\varepsilon$-cap.  We consider two subcases.

\vspace *{0.2cm}

Subcase 2.1:  there is a point $x$ in
the connected component containing $x_0$ such that $R(x,t_0)\leq
r^{-2}$.  In this subcase one can use Proposition \ref{gluing1} and Lemma \ref{lem 4.1} to find a tube  which is diffeomorphic to $\mathbb{S}^{n-1}/\Gamma \times \mathbb{R}$,  where $\Gamma$ is a  finite subgroup of $O(n)$ acting freely on $\mathbb{S}^{n-1}$, such that one end of the tube is adjacent to the $4\varepsilon$-neck or cap centered at $(x_0, t_0)$, and $|\text{Rm}|\leq C(n, \theta, \varepsilon)^2r^{-2}$ on the (strong) $4\varepsilon$-neck near the other end. By  definition of the canonical neighborhood, to  bound $\text{vol}_{t_0}(B(x_0,t_0,r_0))/r_0^n$ from below we only need   to
 control $|\Gamma|$. Then we can reduce this subcase to Case 1; see Step 2 in the proof of Lemma 4.5 in \cite{CTZ}.

\vspace *{0.2cm}

Subcase 2.2:   every point $x$ in the
connected component containing $x_0$ has $R(x,t_0)> r^{-2}$.
Then by our assumption this component, denoted by $Y$, is covered by $4\varepsilon$-caps and/or $4\varepsilon$-necks. In particular $Y$
contains a tube which is diffeomorphic to $\mathbb{S}^{n-1}/\Gamma \times \mathbb{R}$, where $\Gamma$ is a  finite subgroup of $O(n)$ acting freely on $\mathbb{S}^{n-1}$.  Again we only need to control $|\Gamma|$.
Let $t_1=\inf \{t \hspace*{1mm} | \hspace*{1mm}  0\leq t < t_0,  Y\times (t, t_0) \hspace*{1mm} \text{is} \hspace*{1mm} \text{unscathed}\}$.
 Note that  $\inf_{x \in Y} R_{g_+(t_1)}(x)< r^{-2}$.  Let $t_0'= \sup \{t \hspace*{1mm}  | \hspace*{1mm}  t_1 < t < t_0, \hspace*{1mm}  \inf_{x \in Y} R(x,t)=r^{-2}\}$. By the  definition of the canonical neighborhood and the inequality  $|\text{Rm}|\leq C(n)R$  (as noted above),  we  can find $x_0'\in Y$  such that $(x_0',t_0')$ is a center of a strong $4\varepsilon$-neck (diffeomorphic to $\mathbb{S}^{n-1}/\Gamma \times \mathbb{R}$) on which  $|\text{Rm}| \leq
C(n, \theta, \varepsilon)^2r^{-2}$. Then we are in a similar situation
as in Case 1, and we are done; compare Subcase 2.1 and the argument for Case 2 in Section 10.4 of \cite{BBM}.
\hfill{$\Box$}

\vspace *{0.2cm}

Note that the constant $\kappa$ in Proposition \ref{prop 4.9} is independent of the surgery parameters $r,\delta$ and the associated cutoff parameters $h, D, \Theta$.

\vspace *{0.2cm}

Now we can show,  roughly  speaking, a self-improving property of the canonical neighborhood assumption with certain parameters, which is crucial in the justification of the canonical neighborhood assumption. The manifold case is explicitly stated in Theorem 10.10 in \cite{B19}, which has the origin in Perelman's work \cite{P2}.

\begin{prop} \label{prop 4.10} (cf. Theorem 10.10 in \cite{B19})\ \ Let $(\mathcal{O},g_0)$ be a compact Riemannian orbifold of dimension $n$ ($\geq 5$) with at most isolated singularities and with positive isotropic curvature.  Let   $\varepsilon$  and $C$ be chosen as before, and $(f,\theta)$ be given. Then there exist positive numbers $\hat{r}$ and $\hat{\delta}$ with the following property: If we have an $(\hat{r},\hat{\delta})$-surgical solution starting with $(\mathcal{O},g_0)$ which is defined on some time interval $[0,T]$ (in particular, the solution satisfies the $(f,\theta)$-pinching assumption  and the canonical neighborhood assumption
$(CN)_{\hat{r}}$ with $(4\varepsilon,4C)$-control by definition), then this solution satisfies the canonical neighborhood assumption $(CN)_{2\hat{r}}$ with $(2\varepsilon,2C)$-control.
\end{prop}
\noindent {\bf Proof}.\ \  See the proof of Theorem 10.10 in \cite{B19}. Suppose the result is not true. Then there exists a sequence of $(\hat{r}_j, \hat{\delta}_j)$-surgical solutions to the Ricci flow, $(\mathcal{O}^{(j)}(t), g^{(j)}(t))$, $t\in [0,T_j]$,  with $\hat{r}_j\leq \frac{1}{j}$ and $\hat{\delta}_j\leq \text{min}\{\tilde{\delta}(\hat{r}_j), \frac{1}{j}\}$, where $\tilde{\delta}(\cdot)$ is the function in Proposition \ref{prop 4.9}, (in particular, the solutions satisfy the $(f,\theta)$-pinching assumption  and the canonical neighborhood assumptions $(CN)_{\hat{r}_j}$ with $(4\varepsilon,4C)$-control,) starting with $(\mathcal{O}, g_0)$, and a sequence of space-time points $(x_j, t_j)$ ($x_j\in \mathcal{O}^{(j)}(t_j)$) with $Q_j:=R(x_j,t_j)\geq (2\hat{r}_j)^{-2}$ but having no $(2\varepsilon, 2C)$-canonical neighborhoods (in particular, $x_j$ is not contained in a component which is strictly PIC2 at time $t_j$). By Proposition \ref{prop 4.9},  the solution $(\mathcal{O}^{(j)}(t), g^{(j)}(t))$ is $\kappa$-noncollapsed on scales less than $\varepsilon$ for some constant $\kappa>0$ independent of $j$.
We rescale the Riemannian orbifolds $(\mathcal{O}^{(j)}(t_j), g^{(j)}(t_j))$  by the factors $Q_j$. By Proposition \ref{prop 4.4}  the pointed  Riemannian orbifolds $(\mathcal{O}^{(j)}(t_j), Q_jg^{(j)}(t_j), x_j)$  satisfy the so called bounded curvature at bounded distance property.  We also have curvature derivative estimates by (using the argument in the proof of) Proposition 10.4 in \cite{B19}. By Proposition 4.1 in \cite{KL} (see also Theorem 3.5 in \cite{L}) we see that the pointed Riemannian orbifolds $(\mathcal{O}^{(j)}(t_j), Q_jg^{(j)}(t_j), x_j)$  subconverge  in the pointed Cheeger-Gromov sense to a pointed complete Riemannian orbifold $(\mathcal{O}^\infty, g^\infty, x_\infty)$.  By Proposition 6.1 in \cite{CTZ} (see also Proposition 3.18 in \cite{KL}) $(\mathcal{O}^\infty, g^\infty)$ has bounded curvature.

Using Propositions \ref{prop 3.6}, \ref{prop 3.10}, \ref{prop 4.8}, \ref{gluing1} and Lemmas \ref{lem 2.6}, \ref{lem 2.7}, \ref{lem 2.5} (consult the proof of Lemma 9.2.1 in \cite{BBB+} for the use of Lemma \ref{lem 2.5}), we can extend $(\mathcal{O}^\infty, g^\infty)$ backwards in time to an orbifold Ricci flow solution $(\mathcal{O}^\infty, g^\infty(t))$ defined on some time interval $(\tau^*,0]$  (here $\tau^*< 0$), which has bounded curvature for each $t\in (\tau^*,0]$. As in Step 6 in the proof of Theorem 7.2 in \cite{B19}, by using the Harnack inequality in \cite{B09} (which extends the work \cite{H93}), Lemma 8.3(b) in \cite{P1}, the $(f,\theta)$-pinching assumption, the canonical neighborhood assumptions
$(CN)_{\hat{r}_j}$ with $(4\varepsilon,4C)$-control, (the proof of) Proposition \ref{prop 4.4} above, Proposition 4.1 in \cite{KL} and Proposition 6.1 in \cite{CTZ} (see also Proposition 3.18 in \cite{KL}) we can show that $\tau^*=-\infty$. So the parabolical rescalings of the solutions $(\mathcal{O}^{(j)}(t), g^{(j)}(t))$  around $(x_j, t_j)$ with factors $Q_j$ subconverge to an orbifold ancient  $\kappa$-solution satisfying $\text{Rm}-\theta R \hspace*{1mm} \text{id} \owedge \text{id} \in \text{PIC}$.  By using  Propositions \ref{prop 3.5}, \ref{prop 3.6}, \ref{prop 3.9}, \ref{gluing1} and Lemmas \ref{lem 2.6}, \ref{lem 2.7}, we see that $(x_j,t_j)$ has a $(2\varepsilon, 2C)$-canonical neighborhood for $j$ sufficiently large.  A contradiction.
 \hfill{$\Box$}

\vspace *{0.2cm}

Finally we can finish the construction of an $(r,\delta)$-surgical solution starting with a  compact Riemannian orbifold of dimension $n\geq 12$ with at most isolated singularities and with positive isotropic curvature.

\begin{thm} \label{thm 4.11}(cf. Theorem 11.2 in \cite{B19}  and Theorem 4.6 in \cite{CTZ}) \ \  Let $(\mathcal{O},g_0)$ be a compact Riemannian orbifold of dimension $n\geq 12$ with at most isolated singularities and with positive isotropic curvature.  Let  $\varepsilon$  and $C$ be chosen as in the paragraph before the definition of canonical neighborhood assumption. Then  we can find  positive numbers $\hat{r}$ and $\hat{\delta}$   such that there exists an $(\hat{r}, \hat{\delta})$-surgical solution (with $(4\varepsilon, 4C)$-control) to the Ricci flow starting with $(\mathcal{O},g_0)$, which  becomes extinct in finite time.
\end{thm}

\noindent {\bf Proof}.\ \ See the proof of Theorem 11.2 in \cite{B19}.  Let $(\mathcal{O},g_0)$ be as in the assumption. Suppose $|\text{Rm}(g_0)|\leq K$. We run the Ricci flow starting with $(\mathcal{O},g_0)$. Let $\hat{T}=\frac{n}{2\inf_{x\in \mathcal{O}}R(x,0)}$. By  Theorem  1.2  in \cite{B19} and (a version) of Hamilton's maximum principle, there is a continuous family of closed, convex, $O(n)$-invariant sets $\{\mathcal{F}_t ~|~ t\in [0,\hat{T}]\}$ in the vector space $\mathcal{C}_B(\mathbb{R}^n)$ of algebraic curvature tensors in dimension $n$ such that the curvature tensor  of $(\mathcal{O},g_0)$ lies in $\mathcal{F}_0$, the family is invariant under the Hamilton ODE $\frac{d}{dt}\text{Rm}=Q(\text{Rm})$ and satisfies the inclusion relation (3.1) for some suitable $(f,\theta)$ (as mentioned before), and the curvature tensor of the  solution  (to this smooth Ricci flow) at time
$t$  ($< \hat{T}$)  lies in the set $\mathcal{F}_t$, hence also satisfies the pinching assumption.  Choose $\hat{r}$ and $\hat{\delta}$ as in Proposition \ref{prop 4.10}. We may assume that  $R_{g_0}(x)< (2\hat{r})^{-2}$  for any $x\in \mathcal{O}$ and that for $\delta=\hat{\delta}$ the conclusion of (an orbifold version of) Proposition 8.2 in \cite{B19} holds. We may also assume that $\hat{\delta}\leq \tilde{\delta}(\hat{r})$, where $\tilde{\delta}$ is the function in Proposition \ref{prop 4.9}.  We can choose $\alpha_1>0$ such that $(2\hat{r})^{-2}< \max_{(x,t)\in \mathcal{O}\times [0,\alpha_1]}R(x,t) < \hat{r}^{-2}$. Then $(\mathcal{O},g(t))$, $t\in [0, \alpha_1]$, satisfies (vacuously) the canonical neighborhood assumption $(CN)_{\hat{r}}$ with $(4\varepsilon, 4C)$-control.   By Proposition \ref{prop 4.10}, the flow satisfies $(CN)_{2\hat{r}}$ with $(2\varepsilon, 2C)$-control on the time interval $[0, \alpha_1]$.
By continuity, there exists a positive number $\alpha_2$ such that the solution satisfies the canonical neighborhood assumption $(CN)_{\hat{r}}$ with $(4\varepsilon, 4C)$-control for  $t\in [0, \alpha_1+\alpha_2]$; cf. for example, Lemmas 5.3.2 and 5.3.3 in \cite{BBB+}, Proposition 9.79 in \cite{MT}, Claim 1 on p. 1229 of \cite{Hu13}, and Lemma \ref{lem 2.7}.  By Proposition \ref{prop 4.10}, in fact the solution satisfies the canonical neighborhood assumption $(CN)_{2\hat{r}}$ with $(2\varepsilon, 2C)$-control  for  $t\in [0, \alpha_1+\alpha_2]$.

Arguing in this way and using Lemma 9.2.3 in \cite{BBB+} and Claim 17.6 in \cite{MT} we  get a flow which satisfies $(CN)_{2\hat{r}}$ with $(2\varepsilon, 2C)$-control, up until the time $t_1:=\inf \{t~|~ R_{\text{max}}(t)=\Theta\}$, where  $\Theta:=2D(\hat{r},\hat{\delta})h(\hat{r},\hat{\delta})^{-2}$, and $D(\hat{r},\hat{\delta})$ and $h(\hat{r},\hat{\delta})$ are the associated parameters  as
determined in Proposition \ref{prop 4.5}.
  Using Propositions \ref{prop 4.5} we can find strong $\hat{\delta}$-necks centered at points of scalar curvature $h(\hat{r},\hat{\delta})^{-2}$ at this time. Using Proposition \ref{prop 4.6} we  do $(\hat{r},\hat{\delta})$-surgeries on these strong $\hat{\delta}$-necks, throw away all components covered by $(4\varepsilon, 4C)$-canonical neighborhoods, and  get an orbifold $(\mathcal{O}_+,g_+)$  which has curvature tensor lying in the set $\mathcal{F}_{t_1}$  and has $R_{\text{min}}$  non-decreasing after the surgery and $R_{\text{max}}(g_+) \leq \frac{\Theta}{2}$  (if $\mathcal{O}_+$ is not empty). Moreover the post-surgery solution still satisfies the canonical neighborhood assumption $(CN)_{2\hat{r}}$ with $(2\varepsilon, 2C)$-control.
  Then we run the Ricci flow starting with $(\mathcal{O}_+,g_+)$, and iterate the process as before again and again.

  By  Proposition \ref{prop 4.3} the solution that we constructed is a surgical solution, and  furthermore,  an $(\hat{r}, \hat{\delta})$-surgical solution by the previous discussion. Since $\hat{\delta} \leq \tilde{\delta}(\hat{r})$, by Proposition \ref{prop 4.9} there exists $\kappa>0$ such that the $(\hat{r},\hat{\delta})$-surgical solution is $\kappa$-noncollapsed on all scales less than $\varepsilon$ at all points not contained in a component of strictly PIC2.  By the $\kappa$-noncollapsing property each surgery reduces the volume by at least a fixed amount, so the total number of surgeries is bounded above. (By the way, there is another way to see that the total number of surgeries is bounded above; cf. for example
  Lemma 5.9 in \cite{BBM}.) By Proposition 10.2 in \cite{B19} (which is also true in our orbifold situation) the flow with surgery becomes extinct in finite time  ($\leq \hat{T}$).
 \hfill{$\Box$}

\section {Proof of Theorems 1.1 and 1.3}

The following result describes the topological effect  of a Ricci flow surgery.

\begin{lem} \label{lem 5.1}\ \
 At any singular time, the pre-surgery orbifold is diffeomorphic to a  connected sum among the  post-surgery orbifold and the removed components.
\end{lem}
\noindent {\bf Proof}.  Let $t_0$ be a singular time, and suppose that we do a Hamilton's surgery as described before on a $\delta$-neck $N$ with diffeomorphism $\psi:\mathbb{S}^{n-1}/\Gamma \times (-\frac{1}{\delta}, \frac{1}{\delta}) \rightarrow  N$ in a  component, say $\mathcal{O}_1$,  of $\mathcal{O}(t_0)$. After the surgery on $N$ we get a Riemannian orbifold, say $(\mathcal{O}_1',g')$, which may not be  connected.  Let $\Sigma_1= \psi(\mathbb{S}^{n-1}/\Gamma \times \{-(10\delta)^{-1}\})$, and  $\Sigma_2= \psi(\mathbb{S}^{n-1}/\Gamma \times \{(10\delta)^{-1}\})$.  $\Sigma_1$ (resp. $\Sigma_2$) bounds a closed domain, say $D_1$ (resp. $D_2)$, in $\mathcal{O}_1'$ which is diffeomorphic to $D^n//\Gamma$  and is not contained in $\mathcal{O}_1$.   We view  $\psi|_{\mathbb{S}^{n-1}/\Gamma \times \{-(10\delta)^{-1}\}}: \mathbb{S}^{n-1}/\Gamma \times \{-(10\delta)^{-1}\} \rightarrow \Sigma_1 $ as a map $f_1: \mathbb{S}^{n-1}/\Gamma \rightarrow \Sigma_1 $, and  $\psi|_{\mathbb{S}^{n-1}/\Gamma \times \{(10\delta)^{-1}\}}: \mathbb{S}^{n-1}/\Gamma \times \{(10\delta)^{-1}\} \rightarrow \Sigma_2 $ as a map $f_2: \mathbb{S}^{n-1}/\Gamma \rightarrow \Sigma_2 $.   By  inspecting the surgery procedure we see that  $f_i$  extends to a diffeomorphism $\tilde{f}_i: D^n//\Gamma \rightarrow D_i$, $i=1,2$.

Now we  remove from $\mathcal{O}_1'$  the interior of $D_i$, $i=1,2$, and identify the two boundary components thus produced (i.e., $\Sigma_1$ and $\Sigma_2$) using the diffeomorphism $f_2\circ f_1^{-1}$.  The orbifold thus obtained, denoted by $\tilde{\mathcal{O}_1}$, is  the result of a  connected sum operation on $\mathcal{O}_1'$.   Using (the proof of) Theorem 1.9 in Chapter 8 of \cite{Hi}  we see that there is a diffeomorphism
from $\psi(\mathbb{S}^{n-1}/\Gamma \times [-(5\delta)^{-1}, (5\delta)^{-1}])$ to $\psi(\mathbb{S}^{n-1}/\Gamma \times [-(5\delta)^{-1},-(10\delta)^{-1}] \cup \mathbb{S}^{n-1}/\Gamma \times [(10\delta)^{-1}, (5\delta)^{-1}])/\sim$ which is the identity map near the boundary, where  the equivalence relation ``$\sim$'' is
generated by $\psi(p,-(10\delta)^{-1}) \sim \psi(p,(10\delta)^{-1})$ for any $p\in \mathbb{S}^{n-1}/\Gamma$.  It follows that $\mathcal{O}_1$ is diffeomorphic to $\tilde{\mathcal{O}_1}$.
\hfill{$\Box$}

\vspace *{0.2cm}

\noindent {\bf Proof of Theorem 1.3}.  Let $(\mathcal{O},g_0)$ be a compact Riemannian orbifold of dimension $n\geq 12$ with at most isolated singularities and with positive isotropic curvature. By Theorem \ref{thm 4.11} we can construct
an $(r,\delta)$-surgical solution to the Ricci flow starting with $(\mathcal{O}, g_0)$,  which becomes extinct in finite time. Recall that each point in any component that is removed in the process of surgery is contained in a canonical neighborhood, so  by the surgery procedure, definition of canonical neighborhood, and Proposition \ref{prop 4.3}, each such  component  is diffeomorphic to  a spherical orbifold or a  weak connected sum of at most two spherical orbifolds with at most isolated singularities  which admits a metric with positive isotropic curvature.
 Now Theorem 1.3 follows from Lemma \ref{lem 5.1}.  \hfill{$\Box$}

\vspace *{0.2cm}

\noindent {\bf Proof of Theorem 1.1}. Compare the proof of the Main Theorem in \cite{CTZ}. Let a Riemannian manifold $(M, g_0)$ satisfy the assumption of
Theorem 1.1. From  Theorem 1.3 we see that  $M$ is diffeomorphic to  a connected sum of a finite number of orbifolds, called components, where  each component admits a metric with positive isotropic curvature, and is diffeomorphic to either a spherical orbifold with at most isolated singularities or a weak connected sum of  at most two spherical orbifolds with at most isolated singularities. We denote these components by  $\mathcal{O}_1, \mathcal{O}_2, \cdot\cdot\cdot, \mathcal{O}_k$. Note that if a component is a smooth manifold, from the proof of Proposition \ref{prop 4.3} we know that it is diffeomorphic to a spherical manifold or  the total space of an orbifiber bundle over $\mathbb{S}^1$ or $\mathcal{I}$ with generic fiber diffeomorphic to $\mathbb{S}^{n-1}/\Gamma$ (here $\Gamma$  may be trivial).

 First we do  connected sum operations to  resolve the orbifold singularities of the components, which are created by the Ricci flow surgeries. (If a component is a smooth manifold, we don't need to worry about it at this step.) If a component $\mathcal{O}_i$ is diffeomorphic to a spherical orbifold with at least a singular point, by using Lemma \ref{lem 3.2} above, Lemma 5.2 in \cite{CTZ} and Corollary 2.4 in Chapter VI of \cite{B72} we see that it has at most two singularities; we   remove a suitable open neighborhood of each
 singularity and (by the ambient isotopy uniqueness of closed tubular
neighborhoods of compact suborbifolds, which follows from the Remark on p. 312 in Bredon \cite{B72}, cf. p.443 of \cite{BSi}; in fact, here we only need to use the ambient isotopy uniqueness of regular
neighborhoods of compact suborbifolds, cf. p. 32 of \cite{BSi85}) get a manifold  diffeomorphic to some $S^{n-1}/ \Gamma_i \times [-1, 1]$ or  some  $\mathbb{S}^n// \langle\Gamma_i, \hat{\sigma}_i \rangle \setminus B$, which has boundary diffeomorphic to $\mathbb{S}^{n-1}/\Gamma_i$  and can be viewed as an orbifiber bundle over a one-dimensional orbifold with one singular point whose underlying space is a closed interval.

 If a  component  is of the form
 $\mathbb{S}^n// (x,\pm x')  \hspace*{1mm} \sharp \hspace*{1mm} \mathbb{S}^n// (x,\pm x')$  (where the connected sum occurs  at two regular points), we  undo this connected sum,  and get two orbifolds both diffeomorphic to
  $\mathbb{S}^n// (x,\pm x') $. Then  we remove a suitable open neighborhood for each of the four singularities and get  two manifolds both diffeomorphic to  $\mathbb{R}P^{n-1} \times [-1,1]$. If it is of the form  $\mathbb{R}P^n  \sharp \mathbb{S}^n// (x,\pm x')$, we deal with it in a similar way.

 For a component  of the form $\mathbb{S}^n//\Gamma \tilde{\sharp} \mathbb{S}^n//\Gamma $ with $|\Gamma|\geq 2$ (where the weak connected sum occurs at two singular points), it has two singularities; we remove  a suitable open neighborhood for each of the two  singularities and (by for example, ambient isotopy uniqueness of closed tubular
neighborhoods of compact suborbifolds and Proposition 7.7 in \cite{B}) get a manifold diffeomorphic to $\mathbb{S}^{n-1}/\Gamma \times [-1,1]$.

 For a component of the form  $\mathbb{S}^n// \langle\Gamma, \hat{\sigma}\rangle \tilde{\sharp} \mathbb{S}^n//\Gamma$ with $|\Gamma|\geq 2$ (where the weak connected sum occurs at two singular points), it has a unique singularity; we remove a suitable open neighborhood of the unique singularity  and  get a manifold diffeomorphic to $\mathbb{S}^n// \langle\Gamma, \hat{\sigma} \rangle \setminus B$.

 Now we glue all these  manifolds with boundaries   along their boundaries pairwise to finish  the connected sum operations which resolve the orbifold singularities created by the Ricci flow surgeries.  The closed manifolds we thus get, denoted by $Y_j, 1\leq j \leq m$,  are  diffeomorphic to  (the total spaces of) orbifiber bundles over $\mathbb{S}^1$ or $\mathcal{I}$ with generic fibers diffeomorphic to $\mathbb{S}^{n-1}/\Gamma_j$ such that the total spaces  admit  metrics with positive isotropic curvature by an obvious extension of Theorem 1.1 in \cite{MW} to the case of orbifold connected sum (cf. also \cite{CH}, in particular, see Remark 2.2, the paragraph immediately after Proposition 2.3, and Remark 5.6 in \cite{CH}); but here I don't know whether the structure groups can always be reduced to the isometry groups Isom($\mathbb{S}^{n-1}/\Gamma_j$).

 To understand the topology of the total space $Y$ of an orbifiber bundle over $\mathbb{S}^1$ or $\mathcal{I}$ with generic fiber diffeomorphic to $\mathbb{S}^{n-1}/\Gamma$, we construct its universal covering as follows: First note that if $Y$ is diffeomorphic to (the total space of) an orbifiber bundle over  $\mathcal{I}$ with generic fiber diffeomorphic to $\mathbb{S}^{n-1}/\Gamma$, we can take a 2-fold covering to get a $\mathbb{S}^{n-1}/\Gamma$-bundle over $\mathbb{S}^1$; cf. for example  Section 2.7 of \cite{CHK}. Actually this 2-fold covering can be constructed by pulling back the orbifiber bundle over  $\mathcal{I}$ via the 2-fold orbifold covering $\mathbb{S}^1 \rightarrow \mathcal{I}$. (By the way, although in a local model (uniformizing system) of an orbifold considered in \cite{CR} any element of the group  is required to either act trivially or have fixed-point set of codimension at least two (see p. 65 of \cite{CR}), one can extend their notions of orbifolds and good maps slightly so that the 2-fold orbifold covering $\mathbb{S}^1 \rightarrow \mathcal{I}$ is a good map in the extended sense.  Moreover, the pull-back bundle can also be viewed as an orbifold fiber product;  cf. for example p. 83 in \cite{CR} and Chapter 13 in \cite{T}.)  For a $\mathbb{S}^{n-1}/\Gamma$-bundle over $\mathbb{S}^1$,   we can first construct a covering by pulling back this bundle using the universal covering $\mathbb{R}\rightarrow \mathbb{S}^1$, thus we get    a $\mathbb{S}^{n-1}/\Gamma$-bundle over $\mathbb{R}$, which must be  trivial.  Then we see that the universal covering of  $Y$ is diffeomorphic to $\mathbb{S}^{n-1}\times \mathbb{R}$.   So both the manifold $(\mathbb{S}^n// \langle\Gamma, \hat{\sigma} \rangle \setminus B_1)\cup_h (\mathbb{S}^n// \langle\Gamma, \hat{\sigma'} \rangle \setminus B_2)$  appearing in the proof of Proposition \ref{prop 4.3}  and $Y_j$ above are diffeomorphic to  quotients of $\mathbb{S}^{n-1}\times \mathbb{R}$ by the corresponding deck transformations, which are  certain diffeomorphisms.

Then we perform  manifold connected sums to undo the Ricci flow surgeries which do not introduce orbifold singularities. So we do connected sum among $Y_j$, and a finite number of spherical manifolds and a finite number of manifolds diffeomorphic to  the total space of an orbifiber bundle over $\mathbb{S}^1$ or $\mathcal{I}$ with generic fiber diffeomorphic to $\mathbb{S}^{n-1}/\Gamma$  which  are untouched in the previous step. We also  do manifold connected sums which undo the decomposition   of the orbifolds diffeomorphic to $\mathbb{S}^n// (x,\pm x')  \hspace*{1mm} \sharp \hspace*{1mm} \mathbb{S}^n// (x,\pm x')$ or $\mathbb{S}^n// (x,\pm x')  \hspace*{1mm} \sharp \hspace*{1mm}  \mathbb{R}P^n $ performed in the previous step. Thus we recover the original manifold $M$ as a desired manifold connected sum.  \hfill{$\Box$}

\vspace *{0.2cm}

\noindent {\bf Remark}. The proof of Theorem 1.1 implies the following result, which is a very slight extension of Theorem 1.4 in \cite{B19}, and is a higher-dimensional analogue of the main theorem in \cite{H97} (see also \cite{CZ}):
 Let $(M,g_0)$ be a compact
manifold of dimension $n\geq 12$ with  positive isotropic curvature  and containing no  incompressible $(n-1)$-dimensional space forms  $\mathbb{S}^{n-1}/\Gamma$ with $|\Gamma|>2$ or two-sided incompressible $\mathbb{R}P^{n-1}$. Then $M$ is diffeomorphic  to a spherical space form, or $\mathbb{S}^{n-1} \times \mathbb{S}^1$, or $\mathbb{S}^{n-1} \tilde{\times} \mathbb{S}^1$, or a connected sum of a finite number of such manifolds.  Note that the corresponding assumption in Theorem 1.4 in \cite{B19}, ``$M$ does not contain any nontrivial  incompressible $(n-1)$-dimensional space forms'', which means that $M$ contains no   incompressible $(n-1)$-dimensional space forms  $\mathbb{S}^{n-1}/\Gamma$ with $\Gamma \neq \{e\}$, is slightly more restrictive than ours in this remark.  Under our assumption in this remark  the manifold may contain one-sided incompressible $\mathbb{R}P^{n-1}$, so we may have  $\varepsilon$-cap diffeomorphic to  $\mathbb{R}P^n \setminus \bar{B}$,  in addition to $\varepsilon$-cap diffeomorphic to  $\mathbb{R}^n $.  But with our assumption in this remark,  we still get a manifold instead of an orbifold after each surgery.

\vspace *{0.2cm}

\noindent {\bf Proof of Corollary  1.2}.  Compare the proof of Corollary 1 in \cite{CTZ}.  Note that by using Theorem 1.1 and if necessary adding $\mathbb{S}^{n-1}\times \mathbb{S}^1$'s or $\mathbb{S}^{n-1} \tilde{\times} \mathbb{S}^1$'s  to remove the self-connected sums, we can assume that $M$ is diffeomorphic to $\mathbb{S}^n/G_1 \sharp \cdot\cdot\cdot \sharp \mathbb{S}^n/G_k \sharp
\mathbb{S}^{n-1}\times \mathbb{R}/H_1 \sharp \cdot\cdot\cdot \sharp \mathbb{S}^{n-1}\times \mathbb{R}/H_l$, where $G_i$'s are finite subgroups of $O(n+1)$ acting freely on $\mathbb{S}^n$, and $H_j$'s are cocompact discrete subgroups of $\text{Diff}(\mathbb{S}^{n-1}\times \mathbb{R})$ acting freely  and properly  on $\mathbb{S}^{n-1}\times \mathbb{R}$  (for definition of a free and proper action, see \cite{Le}).  (In the case $k=0$ or $l=0$ we understand that there are no $\mathbb{S}^n/G_i$'s or $\mathbb{S}^{n-1}\times \mathbb{R}/H_j$'s.) So
\begin{equation*}\pi_1(M) \cong G_1\ast \cdot\cdot\cdot \ast G_k \ast H_1 \ast \cdot\cdot\cdot \ast H_l.
\end{equation*}
Recall from the proof of Theorem 1.1 that here the manifold $\mathbb{S}^{n-1}\times \mathbb{R}/H_j$  or its certain double cover is diffeomorphic to a $\mathbb{S}^{n-1}/\Gamma_j$-bundle over the circle  (for some finite subgroup $\Gamma_j< O(n)$ acting freely on $\mathbb{S}^{n-1}$), whose fundamental group is an extension of $\Gamma_j$  by $\mathbb{Z}$  (of course this extension must be split). It follows that $H_j$ is virtually infinite cyclic. It is well-known that a subgroup $B$ of a group $A$ of finite index contains a normal subgroup of $A$ of finite index.  So each $H_j$ has a normal subgroup, denoted by $N_j$ ($N_j \lhd H_j)$, which is isomorphic to $\mathbb{Z}$ and is of finite index ($[H_j:N_j]< \infty$).

Let
\begin{equation*}\varphi:G_1\ast \cdot\cdot\cdot \ast G_k \ast H_1 \ast \cdot\cdot\cdot \ast H_l \rightarrow G_1 \times \cdot\cdot\cdot \times G_k\times H_1/N_1 \times \cdot\cdot\cdot \times H_l/N_l
\end{equation*}
 be the surjective homomorphism induced from the identity maps $\text{id}_{G_i}:G_i\rightarrow G_i$ and the natural projections $p_j:H_j\rightarrow H_j/N_j$.
Note that $[G_1\ast \cdot\cdot\cdot \ast G_k \ast H_1 \ast \cdot\cdot\cdot \ast H_l :\text{ker} (\varphi)]=|G_1|\cdot\cdot\cdot|G_k||H_1/N_1|\cdot\cdot\cdot|H_l/N_l|< \infty$. Moreover, by using the Kurosh subgroup theorem we see that  $\text{ker} (\varphi)$ is torsion free. (This type of argument is standard.) By using Theorem 1.4 in \cite{B19} (see also the Remark above) we see that if we take the finite covering space $\tilde{M}$ of $M$ corresponding to the normal subgroup of $\pi_1(M)$ which is isomorphic to $\text{ker} (\varphi)$,  and if necessary we take a further  orientation double cover $\tilde{\tilde{M}}$ of $\tilde{M}$, we get a manifold which is diffeomorphic to $\mathbb{S}^n$, or $\mathbb{S}^{n-1}\times \mathbb{S}^1$, or a connected sum of finite copies of $\mathbb{S}^{n-1}\times \mathbb{S}^1$. (In the case that $\tilde{M}$ is nonorientable, by using Theorem 1.4 in \cite{B19}  we see that $\tilde{M}$ is diffeomorphic to $\mathbb{S}^{n-1} \tilde{\times} \mathbb{S}^1$ or a  connected sum of finite copies of $\mathbb{S}^{n-1} \tilde{\times} \mathbb{S}^1$  and $\mathbb{S}^{n-1} \times \mathbb{S}^1$.) Note that $\pi_1(\tilde{M})$ (and $\pi_1(\tilde{\tilde{M}})$) is a free group  and $[\pi_1(M):\pi_1(\tilde{M})]< \infty$ (and $[\pi_1(M):\pi_1(\tilde{\tilde{M}})]< \infty$).
\hfill{$\Box$}

\vspace *{0.2cm}

\noindent {\bf Remark}.  From  the proof of Corollary 1.2  we see that the fundamental group of a compact manifold of dimension $n\geq 12$ with positive isotropic curvature has a free  subgroup which is normal and of finite index.

\appendix

\section {Appendix: gluing and strengthening of  $\varepsilon$-necks}

The main purpose of  this appendix is to fix some notation and conventions, and collect some technical results on the composition of $\varepsilon$-isometries and the gluing and strengthening of $\varepsilon$-necks.

We follow the conventions and notation in \cite{BBB+} and \cite{BBM} in general. We refer the readers to \cite{BBB+} and \cite{BBM} for some notions used but not defined  here. We also refer the readers to  \cite{BMP}, \cite{CHK}, \cite{KL}, \cite{R}, \cite{Sc}, \cite{T} etc for various notions and properties related to (effective)  orbifolds. We only consider effective orbifolds with $C^\infty$ differential structure.

 To describe the structure of orbifold ancient $\kappa$-solutions and orbifold standard solutions we need  notions on necks and caps.

   Let $\Gamma$ be a  finite subgroup of $O(n)$ acting freely on $\mathbb{S}^{n-1}$.  (Note that if $n$ is odd, then it is well-known that $\Gamma$ must be trivial or $\mathbb{Z}_2$ acting antipodally, and in this case $\mathbb{S}^{n-1}/\Gamma \cong \mathbb{S}^{n-1}$ or $\mathbb{R}P^{n-1}$.)
   Suppose   $\sigma$ is  an isometric involution of the spherical manifold $\mathbb{S}^{n-1}/\Gamma$ with at most isolated fixed points,  let $\hat{\sigma}$ be the
involution on the manifold $\mathbb{S}^{n-1}/\Gamma \times \mathbb{R}$
defined by $\hat{\sigma}(x,s)=(\sigma(x),-s)$ for $x\in
\mathbb{S}^{n-1}/\Gamma $ and $s\in \mathbb{R}$, consider the quotient orbifold $(\mathbb{S}^{n-1}/\Gamma \times
\mathbb{R})//\langle \hat{\sigma}\rangle$, which has at most  isolated singularities. We also denote this
orbifold by $\mathbb{S}^{n-1}/\Gamma \times_{\mathbb{Z}_2}  \mathbb{R}$. By the
way, note that we can consider $\Gamma$ and $\hat{\sigma}$ as isometries of $\mathbb{S}^n$ in a natural way, by lifting $\sigma$ to an isometry of  $\mathbb{S}^{n-1}$ (which is always possible) and viewing $\mathbb{S}^n$ as a suspension of
$\mathbb{S}^{n-1}$. We'll use the same notation for these isometries of $\mathbb{S}^n$. (So $\gamma([(x,s)]):=[(\gamma(x),s)]$ and $\hat{\sigma}([(x,s)]):=[(\tilde{\sigma}(x),-s)]$ for $\gamma\in \Gamma$, $x\in \mathbb{S}^{n-1}$ and $s\in [-1,1]$, where $\tilde{\sigma}\in O(n)$ is a lift of $\sigma$.  Note that the spherical orbifold $\mathbb{S}^n// \langle\Gamma, \hat{\sigma} \rangle$ is independent of the choice of  lift of $\sigma$.)

The $(\mathbb{S}^{n-1}/\Gamma \times
\mathbb{R})//\langle \hat{\sigma}\rangle$  above is a smooth manifold if and only if $\sigma$ has no any fixed points in $\mathbb{S}^{n-1}/\Gamma$; if  this is the case,  we denote this smooth manifold  by $C_\Gamma^\sigma$ as in \cite{CTZ},
 which is diffeomorphic to  $\mathbb{S}^n// \langle\Gamma, \hat{\sigma} \rangle \setminus \bar{B}$,  where $B$ is a small, open metric ball around the unique singular point (in the case that $\Gamma$ is not trivial) of the spherical orbifold $\mathbb{S}^n// \langle\Gamma, \hat{\sigma}\rangle$, and  $\bar{B}$ is the closure of $B$, which is diffeomorphic to $D^n//\Gamma$; of course, when $\Gamma$ is trivial  and $\sigma$ has no any fixed points in $\mathbb{S}^{n-1}$, $\mathbb{S}^n// \langle\Gamma, \hat{\sigma}\rangle$ is just $\mathbb{R}P^n$, which has no singularities, and $B$ can be chosen to be a small, open metric ball around any point  in $\mathbb{R}P^n$   such that the closure  $\bar{B}$  of $B$  is  diffeomorphic to $D^n$.  (The claim that  $C_\Gamma^\sigma$ is diffeomorphic to  $\mathbb{S}^n// \langle\Gamma, \hat{\sigma} \rangle \setminus \bar{B}$ follows, for example, from a fact about $\mathbb{S}^n// \langle\Gamma, \hat{\sigma} \rangle $ shown in the proof of Proposition 3.1 in \cite{Hu23}.) Note that when $n$ is odd, $(\mathbb{S}^{n-1}/\Gamma \times
\mathbb{R})//\langle \hat{\sigma}\rangle$  is a smooth manifold if and only if $\Gamma$ is  trivial   and  $\langle\sigma \rangle$   acts on $\mathbb{S}^{n-1}$  antipodally; see the proof of Proposition \ref{prop 3.5}.

If $\mathbb{S}^{n-1}/\Gamma \times_{\mathbb{Z}_2}  \mathbb{R}$ has nonempty isolated singularities,  by Proposition \ref{prop 3.5}, it must be diffeomorphic to $\mathbb{S}^n// (x,\pm x') \setminus \bar{B}$, where, as in \cite{CTZ},   $\mathbb{S}^n// (x,\pm x')$ denotes the quotient orbifold   $\mathbb{S}^n// \langle\iota \rangle$ with $\iota$  the isometric involution on  $\mathbb{S}^n$ given by
 $(x_1,x_2,\cdot\cdot\cdot,x_{n+1}) \mapsto
(x_1,-x_2,\cdot\cdot\cdot,-x_{n+1})$, which  has two  singular points,  and $B$ is a  small, open  metric ball around  a regular point in $\mathbb{S}^n// (x,\pm x')$ such that the closure  $\bar{B}$  of $B$  is disjoint from the two singular points and  diffeomorphic to $D^n$.  We may choose the regular point to be the image of the north pole under the natural projection $\mathbb{S}^n \rightarrow  \mathbb{S}^n// (x,\pm x')$.

  Let $\mathcal{O}$ be an $n$-dimensional orbifold with at most isolated singularities,  $U$  an open subset, and $\Gamma$ as above.   If $U$ is  diffeomorphic to $\mathbb{S}^{n-1}/\Gamma \times \mathbb{R}$, we call it  a topological neck. As in \cite{CTZ}, we divide topological caps into several types.  If $U$ is diffeomorphic to  $\mathbb{R}^n$  or $C_\Gamma^\sigma$, we call it a smooth cap; if $U$ is diffeomorphic to $\mathbb{R}^n//\Gamma$  with $|\Gamma|\geq 2$, we call it an orbifold cap of type I, and denote it by  $C_\Gamma$; if $U$ is diffeomorphic to $\mathbb{S}^n// (x,\pm x') \setminus \bar{B}$, we call it an orbifold cap of type II.

Now we adapt some definitions from \cite{BBB+} and \cite{BBM}.

\vspace*{0.2cm}

\noindent {\bf Definition} (cf. \cite{BBM}). \ \  Given an interval $I\subset
\mathbb{R}$, an evolving Riemannian orbifold is a family of pairs
$\{(\mathcal{O}(t),g(t))\}_{t \in I}$, where for any $t\in I$,  $\mathcal{O}(t)$ is a
(possibly disconnected or empty) orbifold  with at most isolated singularities, and $g(t)$ is a
Riemannian metric on $\mathcal{O}(t)$.  We say that it is
piecewise $C^1$-smooth if there exists a discrete subset $J$ of $ I$,
 such that the
following conditions are satisfied:

i. On each connected component of $I\setminus J$, $t \mapsto
\mathcal{O}(t)$ is constant, and $t \mapsto g(t)$ is $C^1$-smooth;

ii. For each $t_0\in J$, $\mathcal{O}(t_0)=\mathcal{O}(t)$ for any
$t< t_0$ sufficiently close to $t_0$, and $t\mapsto g(t)$ is left
continuous at $t_0$;

iii. For each $t_0 \in J\setminus \{\sup \hspace*{0.4mm} I\}$, $t\mapsto
(\mathcal{O}(t),g(t))$ has a right limit at $t_0$, denoted by
$(\mathcal{O}_+(t_0),g_+(t_0))$.

\vspace *{0.2cm}

As in \cite{BBM}, a time $t\in I$ is regular if $t$ has a neighborhood in
$I$ where $\mathcal{O}(\cdot)$ is constant and $g(\cdot)$ is
$C^1$-smooth. Otherwise it is singular.

\vspace *{0.2cm}

  \noindent {\bf Definition  }(compare \cite{BBM}, \cite{Hu15}). \ \  A piecewise $C^1$-smooth
  evolving complete Riemannian $n$-orbifold $\{(\mathcal{O}(t), g(t))\}_{t \in I }$ with at most isolated singularities is a
  surgical solution to the Ricci flow if it has the following
  properties.

  i. The equation $\frac{\partial g}{\partial t}=-2 \hspace*{0.4mm} \text{Ric}$ is satisfied
  at all regular times;

  ii. For each singular time $t_0$ there is a   collection
  $\mathcal{S}$ of disjoint embedded $\mathbb{S}^{n-1}/\Gamma$'s in $\mathcal{O}(t_0)$
  (where $\Gamma$'s are finite subgroups of $O(n)$ acting freely on $\mathbb{S}^{n-1})$, and an Riemannian orbifold $\mathcal{O}'$ such that

  (a) $\mathcal{O}'$ is obtained from  $\mathcal{O}(t_0)$ by removing a suitable  open tubular neighborhood of each  element of $\mathcal{S}$ and
  gluing in a Riemannian orbifold  diffeomorphic to $D^n//\Gamma$ along each boundary component thus produced which is diffeomorphic to $\mathbb{S}^{n-1}/\Gamma$;

 (b) $\mathcal{O}_+(t_0)$ is a union of some connected components of $\mathcal{O}'$ and
 $g_+(t_0)=g(t_0)$ on $\mathcal{O}_+(t_0)\cap \mathcal{O}(t_0)$;

(c) each component of $\mathcal{O}'\setminus \mathcal{O}_+(t_0)$ is
diffeomorphic to a spherical orbifold with at most isolated singularities, or a neck, or a cap, or a weak connected sum of at most two spherical orbifolds with at most isolated singularities.

\vspace *{0.2cm}

 Let $\{(\mathcal{O}(t), g(t))\}_{t\in I}$ be a surgical solution to the Ricci flow
and $t_0\in I$. As in \cite{BBM}, if $t_0$ is singular, we set
$\mathcal{O}_{\text{reg}}(t_0):= \mathcal{O}(t_0)\cap \mathcal{O}_+(t_0)$, and  $\mathcal{O}_{\text{sing}}(t_0):=
\mathcal{O}(t_0)\setminus \mathcal{O}_{\text{reg}}(t_0)$. If $t_0$ is regular,
$\mathcal{O}_{\text{reg}}(t_0)=\mathcal{O}(t_0)$ and $\mathcal{O}_{\text{sing}}(t_0)=\emptyset$. Let $t_0\in
[a,b]\subset I$ be a time, and $Y$ be a subset of $\mathcal{O}(t_0)$. If  $Y\subset \mathcal{O}_{\text{reg}}(t)$  for each $t \in [a,b)$, then as in
\cite{BBM}, we say that the set $Y\times [a,b]$ is unscathed.

    Let  $\varepsilon$ be a small positive constant. We adopt the notions of $\varepsilon$-closeness of Riemannian metrics and evolving Riemannian metrics, and $\varepsilon$-isometry and $\varepsilon$-homothety in \cite{BBB+}, except that we replace the $C^{[\varepsilon^{-1}]+1}$-diffeomorphism $\psi$ on p.26 and p.28 in \cite{BBB+} by diffeomorphism $\psi$,  and extend them to the orbifold case.
Following \cite{BBB+}, \cite{BBM}, \cite{CTZ}, \cite{KL08}, \cite{MT} and \cite{P2}, we define $\varepsilon$-neck,  $\varepsilon$-tube,  $\varepsilon$-horn, and strong
$\varepsilon$-neck.   Let   $(\mathcal{O},g)$  be  a Riemannian $n$-orbifold  with at most isolated singularities.
 Given a
point $x_0\in \mathcal{O}$  and $Q>0$, an open subset  $U$  containing $x_0$ is an $\varepsilon$-neck  centered at $x_0$ with scaling factor $Q$ if
there is a diffeomorphism $\psi: \mathbb{S}^{n-1}/\Gamma
  \times (-\varepsilon^{-1},\varepsilon^{-1})  \rightarrow U$   such that the pulled back metric
  $\psi^*g$, rescaled with the factor $Q$, is $\varepsilon$-close (in
  $C^{[\varepsilon^{-1}]}$ topology) to the standard metric on $ \mathbb{S}^{n-1}/\Gamma
  \times (-\varepsilon^{-1},\varepsilon^{-1}) $ with scalar curvature 1  and length $2\varepsilon^{-1}$  for the interval   $(-\varepsilon^{-1},\varepsilon^{-1}) $, and  $x_0 \in \psi(\mathbb{S}^{n-1}/\Gamma \times \{0\})$.   An open subset $U\subset \mathcal{O}$ is called an $\varepsilon$-tube if it  is a union of some $\varepsilon$-necks  and is diffeomorphic to $\mathbb{S}^{n-1}/\Gamma \times (0,1)$. Beware that our convention on $\varepsilon$-tube is slightly different from that (Definition 3.1.3) in \cite{BBB+}.
 If an $\varepsilon$-tube has scalar curvature bounded on one end and tending to $+\infty$ on the other end, we call it an $\varepsilon$-horn.

\vspace *{0.2cm}

\noindent {\bf Definition}. (cf. Definition 4.4 in \cite{BBM} and Definition 3.6 in \cite{CTZ})
Given  an evolving $n$-dimensional Riemannian orbifold   $\{(\mathcal{O}(t), g(t))\}_{t\in I}$ and a space-time point $(x_0,t_0)$, an open subset
$U \subset \mathcal{O}(t_0)$ is a strong
  $\varepsilon$-neck  centered at $(x_0,t_0)$ if  the parabolic region $\{(x,t)| x \in U, t\in [t_0-Q^{-1},t_0]\}$ (for some $Q>0$) is unscathed and  there
  is a diffeomorphim $\psi: \mathbb{S}^{n-1}/\Gamma
  \times (-\varepsilon^{-1},\varepsilon^{-1}) \rightarrow U$ such that, the evolving metric  $g(\cdot)$ on the parabolic region $\{(x,t)| x \in U, t\in [t_0-Q^{-1},t_0]\}$, after pulling back via $\psi$ and
  parabolically rescaling with factor $Q$ at time $t_0$, that is, the evolving metric $Q\psi^*g(t_0+\frac{t}{Q})$ on the region $\{(x,t)| x \in U, t\in [-1,0]\}$, is $\varepsilon$-close (in $C^{[\varepsilon^{-1}]}$
  topology) to the round cylinder solution to the Ricci flow on the space-time region  $\mathbb{S}^{n-1}/\Gamma
  \times (-\varepsilon^{-1},\varepsilon^{-1}) \times [-1,0]$, with scalar curvature one and length $2\varepsilon^{-1}$  for the interval   $(-\varepsilon^{-1},\varepsilon^{-1}) $ at time zero, and
   $x_0 \in \psi(\mathbb{S}^{n-1}/\Gamma \times \{0\})$.

\vspace *{0.2cm}

The following elementary result is essentially (3.5) in \cite{B07}, except that we point out that the constant $E$ there does not depend on the dimension $n$.
\begin{lem} \label{lem ba} (cf. \cite{B07}) Let $g_1$ and $g_2$ be two Riemannian metrics on an orbifold $\mathcal{O}$, and $T$ be a $(p,q)$-tensor field on $\mathcal{O}$ with $p+q > 0$. Then for any $k\geq 1$ we have
\begin{equation*}
\begin{split}
 & |\nabla_{g_2}^kT-\nabla_{g_1}^kT|_{g_1} \\
\leq &  C(k,p,q)\sum_{\substack{i_0+\cdot\cdot\cdot+i_m=k\\[3pt]  \hspace*{0.4mm} i_0<k,  \hspace*{0.4mm} i_1, \cdot\cdot\cdot, i_m\geq 1,   \hspace*{0.4mm} m\geq 1}}|\nabla_{g_2}^{i_0}T|_{g_1}|\nabla_{g_2}^{i_1}g_1|_{g_1}\cdot\cdot\cdot |\nabla_{g_2}^{i_m}g_1|_{g_1},
  \end{split}
\end{equation*}
where $C(k,p,q)$ is a constant depending only on $k$, $p$, and $q$.
\end{lem}
\noindent {\bf Proof}. We follow the proof of (3.5) in \cite{B07}; compare the proof of Lemma 4.5 in \cite{CC}. We induct on $k$.  For $k=1$ we have
\begin{equation*}
|\nabla_{g_2}T-\nabla_{g_1}T|_{g_1}\leq(p+q)|\Gamma_{g_2}-\Gamma_{g_1}|_{g_1}|T|_{g_1}\leq \frac{3}{2}(p+q)|T|_{g_1} |\nabla_{g_2}g_1|_{g_1},
\end{equation*}
where $\Gamma_{g_i}$ are the Christoffel symbols of the Levi-Civita connections $\nabla_{g_i}$ w.r.t. the metrics $g_i$, $i=1, 2$.

Assume the inequality holds for $k$. Then
\begin{equation*}
\begin{split}
 & |\nabla_{g_2}^{k+1}T-\nabla_{g_1}^k\nabla_{g_2}T|_{g_1} \\
\leq &  C(k, p,q+1)\sum_{\substack{i_0+\cdot\cdot\cdot+i_m=k\\[3pt]  \hspace*{0.4mm} i_0<k,  \hspace*{0.4mm} i_1, \cdot\cdot\cdot, i_m\geq 1,   \hspace*{0.4mm} m\geq 1}}|\nabla_{g_2}^{i_0+1}T|_{g_1}|\nabla_{g_2}^{i_1}g_1|_{g_1}\cdot\cdot\cdot |\nabla_{g_2}^{i_m}g_1|_{g_1},
  \end{split}
\end{equation*}
and for $T'=\nabla_{g_2}T-\nabla_{g_1}T$,
\begin{equation*}
\begin{split}
 & |\nabla_{g_2}^kT'-\nabla_{g_1}^{k}T'|_{g_1} \\
\leq &  C(k,p,q+1)\sum_{\substack{i_0+\cdot\cdot\cdot+i_m=k\\[3pt]  \hspace*{0.4mm} i_0<k,  \hspace*{0.4mm} i_1, \cdot\cdot\cdot, i_m\geq 1,   \hspace*{0.4mm} m\geq 1}}|\nabla_{g_2}^{i_0}T'|_{g_1}|\nabla_{g_2}^{i_1}g_1|_{g_1}\cdot\cdot\cdot |\nabla_{g_2}^{i_m}g_1|_{g_1}.
  \end{split}
\end{equation*}
Note that
\begin{equation*}
\begin{split}
 & |\nabla_{g_2}^j(\Gamma_{g_2}-\Gamma_{g_1})|_{g_1} \\
\leq & C(j)\sum_{\substack{i_1+\cdot\cdot\cdot+i_m=j+1\\[3pt]  \hspace*{0.4mm}  \hspace*{0.4mm} i_1, \cdot\cdot\cdot, i_m\geq 1,   \hspace*{0.4mm} m\geq 1}}|\nabla_{g_2}^{i_1}g_1|_{g_1}\cdot\cdot\cdot |\nabla_{g_2}^{i_m}g_1|_{g_1},
  \end{split}
\end{equation*}
and
\begin{equation*}
\begin{split}
 & |\nabla_{g_2}^iT'|_{g_1} \\
\leq & (p+q)\sum_{j=0}^i \binom{i}{j}|\nabla_{g_2}^j(\Gamma_{g_2}-\Gamma_{g_1})|_{g_1}|\nabla_{g_2}^{i-j}T|_{g_1} \\
\leq & (p+q) C(i)\sum_{\substack{i_0+\cdot\cdot\cdot+i_m=i+1\\[3pt]  \hspace*{0.4mm} i_0\leq i,  \hspace*{0.4mm} i_1, \cdot\cdot\cdot, i_m\geq 1,   \hspace*{0.4mm} m\geq 1}}|\nabla_{g_2}^{i_0}T|_{g_1}|\nabla_{g_2}^{i_1}g_1|_{g_1}\cdot\cdot\cdot |\nabla_{g_2}^{i_m}g_1|_{g_1}.
  \end{split}
\end{equation*}
It follows that
\begin{equation*}
\begin{split}
 & |\nabla_{g_2}^{k+1}T-\nabla_{g_1}^{k+1}T|_{g_1} \\
\leq & |\nabla_{g_2}^{k+1}T-\nabla_{g_1}^k\nabla_{g_2}T|_{g_1}+ |\nabla_{g_1}^k\nabla_{g_2}T-\nabla_{g_1}^{k+1}T|_{g_1} \\
\leq & |\nabla_{g_2}^{k+1}T-\nabla_{g_1}^k\nabla_{g_2}T|_{g_1}+ |\nabla_{g_2}^kT'-\nabla_{g_1}^{k}T'|_{g_1} + |\nabla_{g_2}^kT'|_{g_1} \\
\leq &  C(k+1,p,q)\sum_{\substack{i_0+\cdot\cdot\cdot+i_m=k+1\\[3pt]  \hspace*{0.4mm} i_0\leq k,  \hspace*{0.4mm} i_1, \cdot\cdot\cdot, i_m\geq 1,   \hspace*{0.4mm} m\geq 1}}|\nabla_{g_2}^{i_0}T|_{g_1}|\nabla_{g_2}^{i_1}g_1|_{g_1}\cdot\cdot\cdot |\nabla_{g_2}^{i_m}g_1|_{g_1},
  \end{split}
\end{equation*}
and we are done.
 \hfill{$\Box$}

\vspace *{0.2cm}

The following result should be well-known.  It is used in the proof of Proposition \ref{prop 4.10}.

\begin{lem} \label{lem 2.6}  Let $n \geq 2$. For any $0<\varepsilon <0.1$, there exists $0< \tilde{\varepsilon}_2=\tilde{\varepsilon}_2(\varepsilon)\leq \varepsilon$ with the following property.  Let  $(\mathcal{O}_1, g_1) $  and   $(\mathcal{O}_2, g_2) $ be two Riemannian $n$-orbifolds, and $\psi:\mathcal{O}_1 \rightarrow \mathcal{O}_2$ be an $\eta$-isometry with $0<\eta \leq \tilde{\varepsilon}_2$.  Then $\psi^{-1}: \mathcal{O}_2 \rightarrow \mathcal{O}_1$ is an $\varepsilon$-isometry.
\end{lem}

\noindent {\bf Proof}.  Given  $0<\varepsilon <0.1$,  let $0<\eta \leq\varepsilon$, and $\psi:(\mathcal{O}_1,g_1) \rightarrow (\mathcal{O}_2,g_2)$ be an $\eta$-isometry. Then
\begin{equation*}
\sup_{x\in \mathcal{O}_1} \sum_{k=0}^{[\eta^{-1}]}|\nabla_{g_1}^k(\psi^*g_2-g_1)(x)|_{g_1}^2< \eta^2,
\end{equation*}
where $\nabla_{g_1}$ is the Levi-Civita connection of the metric $g_1$.  In particular,

\begin{equation*}
(1+2\eta)^{-1}g_1 <(1-\eta)g_1 < \psi^*g_2 < (1+\eta)g_1 <  (1+2\eta)g_1,
\end{equation*}
and by Lemma 3.13 in \cite{CC},
\begin{equation*}
(1+2\eta)^{-(p+q)/2}|T|_{g_1} \leq |T|_{\psi^*g_2} \leq  (1+2\eta)^{(p+q)/2}|T|_{g_1}
\end{equation*}
for any $(p,q)$-tensor field $T$ on $\mathcal{O}_1$.

 For $1\leq k \leq [\varepsilon^{-1}]$, by Lemma \ref{lem ba} we have
\begin{equation*}
\begin{split}
& |\nabla_{g_1}^k(\psi^*g_2)+\nabla_{\psi^*g_2}^kg_1|_{\psi^*g_2}    \\
= & |\nabla_{g_1}^k(\psi^*g_2-g_1)-\nabla_{\psi^*g_2}^k(\psi^*g_2-g_1)|_{\psi^*g_2} \\
\leq &  C(k)\sum_{\substack{i_0+\cdot\cdot\cdot+i_m=k\\[3pt]  \hspace*{0.4mm} i_0<k,  \hspace*{0.4mm} i_1, \cdot\cdot\cdot, i_m\geq 1,   \hspace*{0.4mm} m\geq 1}}|\nabla_{g_1}^{i_0}(\psi^*g_2-g_1)|_{\psi^*g_2}|\nabla_{g_1}^{i_1}(\psi^*g_2)|_{\psi^*g_2}\cdot\cdot\cdot |\nabla_{g_1}^{i_m}(\psi^*g_2)|_{\psi^*g_2},
  \end{split}
\end{equation*}
where $\nabla_{\psi^*g_2}$ is the Levi-Civita connection of the metric $\psi^*g_2$, and $C(k)$ is a constant depending only on  $k$.  Note that
\begin{equation*}
|\nabla_{g_1}^{i}(\psi^*g_2)|_{\psi^*g_2}\leq (1+2\eta)^{(2+i)/2} |\nabla_{g_1}^{i}(\psi^*g_2)|_{g_1} \leq \eta  (1+2\eta)^{(2+i)/2}
\end{equation*}
for $i>0$.
So  for $1\leq k \leq [\varepsilon^{-1}]$  we have
\begin{equation*}
\begin{split}
& |\nabla_{\psi^*g_2}^kg_1|_{\psi^*g_2}    \\
\leq & |\nabla_{g_1}^k(\psi^*g_2)|_{\psi^*g_2} + |\nabla_{g_1}^k(\psi^*g_2)+\nabla_{\psi^*g_2}^kg_1|_{\psi^*g_2} \\
 \leq & \eta  (1+2\eta)^{(2+k)/2} + C(k)\sum_{\substack{i_0+\cdot\cdot\cdot+i_m=k\\[3pt] \hspace*{0.4mm} i_0<k, \hspace*{0.4mm} i_1, \cdot\cdot\cdot, i_m\geq 1, \hspace*{0.4mm} m\geq 1}} (1+2\eta)^{(2+i_0)/2}|\nabla_{g_1}^{i_0}(\psi^*g_2-g_1)|_{g_1}\\
 & \cdot |\nabla_{g_1}^{i_1}(\psi^*g_2)|_{\psi^*g_2}\cdot\cdot\cdot |\nabla_{g_1}^{i_m}(\psi^*g_2)|_{\psi^*g_2}\\
\leq  &  C(\varepsilon) \eta.
\end{split}
\end{equation*}
 It follows that if $0< \tilde{\varepsilon}_2=\tilde{\varepsilon}_2(\varepsilon)\leq \varepsilon$ is sufficiently small and $\eta \leq \tilde{\varepsilon}_2$, there holds
\begin{equation*}
\sup_{x\in \mathcal{O}_1} \sum_{k=0}^{[\varepsilon^{-1}]}|\nabla_{\psi^*g_2}^k(\psi^*g_2-g_1)(x)|_{\psi^*g_2}^2< \varepsilon^2,
\end{equation*}
  that is,
  \begin{equation*}
\sup_{y\in \mathcal{O}_2} \sum_{k=0}^{[\varepsilon^{-1}]}|\nabla_{g_2}^k((\psi^{-1})^*g_1-g_2)(y)|_{g_2}^2< \varepsilon^2,
\end{equation*}
 and we are done.    \hfill{$\Box$}

\vspace *{0.2cm}

The following result should also be well-known; cf. for example Proposition 9.79 in \cite{MT}. It is used in the proof of Proposition \ref{prop 4.10}  and Lemma \ref{lem 2.5}.

\begin{lem} \label{lem 2.7}   Let $n \geq 2$. For any $0<\varepsilon <0.1$  and $\varepsilon'>\varepsilon$, there exists $\tilde{\varepsilon}_3=\tilde{\varepsilon}_3(\varepsilon, \varepsilon')>0$ with the following property.  Let  $(\mathcal{O}_1, g_1) $,   $(\mathcal{O}_2, g_2) $, and $(\mathcal{O}_3, g_3) $ be  Riemannian $n$-orbifolds,  $\psi_1:\mathcal{O}_1 \rightarrow \mathcal{O}_2$ be an $\varepsilon$-isometry, and $\psi_2:\mathcal{O}_2 \rightarrow \mathcal{O}_3$ be an $\eta$-isometry with $0<\eta \leq \tilde{\varepsilon}_3$.  Then $\psi_2 \circ\psi_1: \mathcal{O}_1 \rightarrow \mathcal{O}_3$ is an $\varepsilon'$-isometry.
\end{lem}
\noindent {\bf Proof}.  Compare the proof of Proposition 4.7 in \cite{CC} (here our assumption is somewhat weaker). Let $0<\varepsilon <0.1$  and $\varepsilon'>\varepsilon$, and $(\mathcal{O}_i,g_i)$, $i=1,2,3$, be Riemannian orbifolds. Let $\psi_1:\mathcal{O}_1 \rightarrow \mathcal{O}_2$ be an $\varepsilon$-isometry, and $\psi_2:\mathcal{O}_2 \rightarrow \mathcal{O}_3$ be an $\eta$-isometry with $0< \eta < \varepsilon'$.  Denote $T=\psi_2^*g_3-g_2$.  Then
\begin{equation*}
\sup_{x\in \mathcal{O}_1} \sum_{k=0}^{[\varepsilon^{-1}]}|\nabla_{g_1}^k(\psi_1^*g_2-g_1)(x)|_{g_1}^2< \varepsilon^2,
\end{equation*}
 and
\begin{equation*}
\sup_{y\in \mathcal{O}_2} \sum_{k=0}^{[\eta^{-1}]}|\nabla_{g_2}^kT(y)|_{g_2}^2< \eta^2,
\end{equation*}
where $\nabla_{g_i}$ is the Levi-Civita connection of the metric $g_i$.

As in the proof of Lemma \ref{lem 2.6}, we have
\begin{equation*}
(1+2\varepsilon)^{-1}g_1 <\psi_1^*g_2 < (1+2\varepsilon)g_1,
\end{equation*}
so
\begin{equation*}
\begin{split}
& |\psi_1^*\psi_2^*g_3-g_1|_{g_1}    \\
\leq & |\psi_1^*\psi_2^*g_3-\psi_1^*g_2|_{g_1} + |\psi_1^*g_2-g_1|_{g_1}    \\
\leq & (1+2\varepsilon)|\psi_1^*T|_{\psi_1^*g_2} + |\psi_1^*g_2-g_1|_{g_1}.    \\
  \end{split}
\end{equation*}

 For $1\leq k \leq [\varepsilon'^{-1}]$,  by Lemma \ref{lem ba}, we have
\begin{equation*}
 |\nabla_{\psi_1^*g_2}^kg_1|_{g_1}\leq C(k) \varepsilon
\end{equation*}
as in the proof of Lemma \ref{lem 2.6}, where $\nabla_{\psi_1^*g_2}$ is the Levi-Civita connection of the metric $\psi_1^*g_2$, and here and below we  use $C$ to denote various constants which may be different in different places, and
\begin{equation*}
\begin{split}
 & |\nabla_{g_1}^k(\psi_1^*T)-\nabla_{\psi_1^*g_2}^k(\psi_1^*T)|_{g_1} \\
\leq &  C(k)\sum_{\substack{i_0+\cdot\cdot\cdot+i_m=k\\[3pt]  \hspace*{0.4mm} i_0<k, \hspace*{0.4mm} i_1, \cdot\cdot\cdot, i_m \geq 1, \hspace*{0.4mm} m\geq 1}}|\nabla_{\psi_1^*g_2}^{i_0}(\psi_1^*T)|_{g_1}|\nabla_{\psi_1^*g_2}^{i_1}g_1|_{g_1}\cdot\cdot\cdot |\nabla_{\psi_1^*g_2}^{i_m}g_1|_{g_1},
  \end{split}
\end{equation*}
so
\begin{equation*}
\begin{split}
 & |\nabla_{g_1}^k(\psi_1^*T)|_{g_1} \\
\leq & |\nabla_{g_1}^k(\psi_1^*T)-\nabla_{\psi_1^*g_2}^k(\psi_1^*T)|_{g_1}+ |\nabla_{\psi_1^*g_2}^k(\psi_1^*T)|_{g_1}\\
\leq &  C(k)\sum_{\substack{i_0+\cdot\cdot\cdot+i_m=k\\[3pt] \hspace*{0.4mm} i_0<k,\hspace*{0.4mm} i_1,\cdot\cdot\cdot,i_m\geq 1, \hspace*{0.4mm} m\geq 1}} (1+2\varepsilon)^{(2+i_0)/2}|\nabla_{\psi_1^*g_2}^{i_0}(\psi_1^*T)|_{\psi_1^*g_2}|\nabla_{\psi_1^*g_2}^{i_1}g_1|_{g_1}\cdot\cdot\cdot |\nabla_{\psi_1^*g_2}^{i_m}g_1|_{g_1}\\
 & + (1+2\varepsilon)^{(2+k)/2}|\nabla_{\psi_1^*g_2}^k(\psi_1^*T)|_{\psi_1^*g_2}\\
 \leq &  C(\varepsilon, \varepsilon')\eta.
  \end{split}
\end{equation*}
Choose $\theta=\frac{(\varepsilon'/\varepsilon)^2-1}{4}$.
Then  if $ \tilde{\varepsilon}_3=\tilde{\varepsilon}_3(\varepsilon, \varepsilon')>0$ is sufficiently small and $\eta \leq \tilde{\varepsilon}_3$, we have
\begin{equation*}
\begin{split}
 &  \sup_{x\in \mathcal{O}_1} \sum_{k=0}^{[\varepsilon'^{-1}]}|\nabla_{g_1}^k(\psi_1^*\psi_2^*g_3-g_1)(x)|_{g_1}^2\\
 \leq & \sup_{x\in \mathcal{O}_1} \sum_{k=0}^{[\varepsilon'^{-1}]}  ( (1+\theta^{-1})|\nabla_{g_1}^k(\psi_1^*T)(x)|_{g_1}^2  +(1+\theta)|\nabla_{g_1}^k(\psi_1^*g_2-g_1)(x)|_{g_1}^2) \\
 \leq & (1+\theta^{-1}) C(\varepsilon, \varepsilon')\eta^2+ (1+\theta)\varepsilon^2 \\
 < &  \varepsilon'^2,
 \end{split}
\end{equation*}
  and we are done.    \hfill{$\Box$}

\vspace *{0.2cm}

\noindent {\bf Remark}.  Similarly one can show the following result: Let $n \geq 2$. For any $0<\varepsilon <0.1$  and $\varepsilon'>\varepsilon$, there exists $\tilde{\varepsilon}'_3=\tilde{\varepsilon}'_3(\varepsilon, \varepsilon')>0$ with the following property.  Let  $(\mathcal{O}_1, g_1) $,   $(\mathcal{O}_2, g_2) $, and $(\mathcal{O}_3, g_3) $ be  Riemannian $n$-orbifolds,  $\psi_1:\mathcal{O}_1 \rightarrow \mathcal{O}_2$ be an $\eta$-isometry with $0<\eta \leq \tilde{\varepsilon}'_3$, and $\psi_2:\mathcal{O}_2 \rightarrow \mathcal{O}_3$ be an $\varepsilon$-isometry.  Then $\psi_2 \circ\psi_1: \mathcal{O}_1 \rightarrow \mathcal{O}_3$ is an $\varepsilon'$-isometry.

\vspace *{0.2cm}

The following result on the gluing of $\varepsilon$-necks is a higher-dimensional analogue  of  item (5) of Proposition A.11 in \cite{MT} and Lemma 3.2.2 in \cite{BBB+}.

\begin{prop} \label{gluing1}
 Let $n\geq 4$. There exists a positive constant $\varepsilon_1=\varepsilon_1(n)$ with the following property. Let $(\mathcal{O},g)$ be a Riemannian orbifold of dimension $n$ with at most isolated singularities, $0< \varepsilon \leq 2\varepsilon_1$, and $N_i$ be an $\varepsilon$-neck centered at a point $x_i$ in $\mathcal{O}$ given by a diffeomorphism $\psi_i: \mathbb{S}^{n-1}/\Gamma_i
  \times (-\varepsilon^{-1}, \varepsilon^{-1})  \rightarrow N_i$, $i=1,2$, such that  $N_1  \cap \psi_2(\mathbb{S}^{n-1}/\Gamma_2
  \times [-0.9\varepsilon^{-1}, 0.9\varepsilon^{-1}]) \neq \emptyset$.
   Then $\Gamma_1$ and $\Gamma_2$ are conjugate in $O(n)$, and there is a cross-section of the neck $N_1$ (that is the image $\psi_1( \mathbb{S}^{n-1}/\Gamma_1
  \times \{z\})$ for some $z\in (-\varepsilon^{-1}, \varepsilon^{-1})$) and a cross-section of the neck $N_2$ which cobound a compact domain in  $N_1\cup N_2$ diffeomorphic to the product   $\mathbb{S}^{n-1}/\Gamma_1
   \times [0,1]$.
\end{prop}

\noindent {\bf Proof}.\ \    We use a Morse theoretical argument.  Note that  Lemma A.4 in \cite{MT} (there is also a statement in the proof of Lemma G1.4 in \cite{H97} which is similar to  Lemma A.4 in \cite{MT}) can be extended to the case of $n$-dimensional $\varepsilon$-necks. (For example, to extend  Lemma A.4  in \cite{MT} to our situation we can combine the argument in the proof of Lemma A.4  in \cite{MT} with techniques from the covering space theory.)     Since  $N_1  \cap \psi_2(\mathbb{S}^{n-1}/\Gamma_2
  \times [-0.9\varepsilon^{-1}, 0.9\varepsilon^{-1}]) \neq \emptyset$, if $\varepsilon_1=\varepsilon_1(n)$ is sufficiently small and $0< \varepsilon \leq 2\varepsilon_1$,  the ratio of the scales of the two $\varepsilon$-necks $N_1$ and $N_2$ is very close to 1, and there exist numbers $a,b \in (-\varepsilon^{-1},\varepsilon^{-1})$  with $a<b$ and $b-a\geq 0.099 \varepsilon^{-1}$ such that $\psi_1(\mathbb{S}^{n-1}/\Gamma_1 \times [a, b])\subset N_2$. Then one sees that there exists a cross section, say $S$, of the neck $N_2$  such that $S\subset \psi_1(\mathbb{S}^{n-1}/\Gamma_1 \times [a+0.093\varepsilon^{-1}, a+0.094\varepsilon^{-1}])$.
  We rescale the metric $g$ by a constant factor such  that  after rescaling  the scalar curvature on $\psi_1(\mathbb{S}^{n-1}/\Gamma_1 \times [a, b])$ is close to 1. Let $p_2: \mathbb{S}^{n-1}/\Gamma_1
  \times \mathbb{R}  \rightarrow \mathbb{R}$ and $p_2': \mathbb{S}^{n-1}/\Gamma_2
  \times \mathbb{R}  \rightarrow \mathbb{R}$   be the natural projections onto the second factor. The following estimates are for quantities  computed w.r.t. the rescaled metric (on $\psi_1(\mathbb{S}^{n-1}/\Gamma_1 \times [a, b])$) under the assumption that  $0< \varepsilon\leq 2\varepsilon_1(n)$ with  $\varepsilon_1(n)$ sufficiently small. We have
  $||\nabla(p_2\circ \psi_1^{-1})|-1|\ll 1$ and $||\nabla (p_2'\circ \psi_2^{-1})|-1|\ll 1$ on $\psi_1(\mathbb{S}^{n-1}/\Gamma_1 \times [a, b])$.  Using the $n$-dimensional extension of Lemma A.4 in \cite{MT} we see that after replacing  $\psi_2$ by $\psi_2$ precomposed with a   reflection (that is reversing the direction of the $\varepsilon$-neck structure of $N_2$) if necessary,  the directions of the two vector fields, that is,  the gradient vector fields $\nabla (p_2\circ \psi_1^{-1})$ and $\nabla (p_2'\circ \psi_2^{-1})$ both restricted to the region $\psi_1(\mathbb{S}^{n-1}/\Gamma_1 \times [a, b])$, are very close to each other; one can also see this by using the fact that both  $\nabla (p_2\circ \psi_1^{-1})$ and  $\nabla (p_2'\circ \psi_2^{-1})$ (restricted to the region $\psi_1(\mathbb{S}^{n-1}/\Gamma_1 \times [a, b])$),  are very close to the eigenspace of Ric (viewed as a (1,1)-tensor) corresponding to the smallest eigenvalue at each point of $\psi_1(\mathbb{S}^{n-1}/\Gamma_1 \times [a, b])$ (this fact in turn follows from the closeness of the metric $g$ (restricted to $\psi_1(\mathbb{S}^{n-1}/\Gamma_1 \times [a, b])$) to the metrics on the corresponding parts of the two quotients $\mathbb{S}^{n-1}/\Gamma_i
  \times (-\varepsilon^{-1}, \varepsilon^{-1})$ ($i=1,2$) of the round cylinder). Cf. the proof of Lemma 5.4.5 in \cite{B07}.  So $|\nabla(p_2\circ \psi_1^{-1}-p_2'\circ \psi_2^{-1})|\ll 1$ on $\psi_1(\mathbb{S}^{n-1}/\Gamma_1 \times [a, b])$.
After  translations of the parameters for the parametrizations $\psi_i$ (the shifted parametrizations will still be denoted by $\psi_i$) we may assume that $\psi_1(\mathbb{S}^{n-1}/\Gamma_1 \times [0.9\varepsilon^{-1}, 0.999\varepsilon^{-1}])\subset N_2$ and $\psi_1(\mathbb{S}^{n-1}/\Gamma_1 \times \{0.95\varepsilon^{-1}\}) \cap \psi_2(\mathbb{S}^{n-1}/\Gamma_2 \times \{0.95\varepsilon^{-1}\}) \neq \emptyset$.

   Let $W=\psi_1(\mathbb{S}^{n-1}/\Gamma_1 \times [0.9\varepsilon^{-1}, 0.999\varepsilon^{-1}]) $,   $\rho: W \rightarrow [0,1]$ be a smooth cutoff function   which  is 1 on $\psi_1( \mathbb{S}^{n-1}/\Gamma_1
  \times [0.904 \varepsilon^{-1}, 0.995\varepsilon^{-1}])$  and 0 outside  $\psi_1( \mathbb{S}^{n-1}/\Gamma_1
  \times (0.901 \varepsilon^{-1}, 0.998\varepsilon^{-1}))$,  and $\zeta: W \rightarrow [0,1]$ be a suitable smooth  function   which  is 1 on $\psi_1( \mathbb{S}^{n-1}/\Gamma_1
  \times [0.9\varepsilon^{-1}, 0.907\varepsilon^{-1}])$  and 0 on $\psi_1( \mathbb{S}^{n-1}/\Gamma_1
  \times [0.992\varepsilon^{-1}, 0.999\varepsilon^{-1}])$   with $|\nabla \zeta |< 20\varepsilon$.  We construct  a smooth   function $h=\zeta p_2\circ \psi_1^{-1}|_W+(1-\zeta)p_2'\circ \psi_2^{-1}|_W$ on $W$. Note that $|p_2\circ \psi_1^{-1}-p_2'\circ \psi_2^{-1}|\cdot |\nabla \zeta|\ll 1$ on $W $. So if  $0<\varepsilon\leq 2\varepsilon_1(n)$ with  $\varepsilon_1(n)$ sufficiently small, the gradient vector field $\nabla h\neq 0$ everywhere on $W$. Consider the vector field $X=\frac{\rho}{|\nabla h|^2}\nabla h$ on  $W$. (Cf. Milnor \cite{M69}.) Using the flow generated by the vector field $X$  we see that  the  hypersurfaces   $\psi_1(\mathbb{S}^{n-1}/\Gamma_1
  \times \{0.906\varepsilon^{-1}\})$ and $S$ (defined in the last paragraph) cobound a compact domain in $W$ which is diffeomorphic to the product   $\mathbb{S}^{n-1}/\Gamma_1
    \times [0,1]$. (Compare the proof of item (5) of Proposition A.11 in \cite{MT}, the proof of Theorem G1.1 in \cite{H97}, the last paragraph in the proof of Theorem 3.9 in \cite{CTZ},  and the second to last paragraph in the proof of Proposition 3.6 in \cite{Hu15}.) In particular, the two hypersurfaces are diffeomorphic.   By \cite{dR} and Lemma 2.5.6 in \cite{Wo}, $\Gamma_1 $ and $\Gamma_2$ are conjugate in $O(n)$.  
\hfill{$\Box$}

\vspace *{0.2cm}

\noindent {\bf Remark}.    In the special case that both $\Gamma_1$ and $\Gamma_2$ are trivial and the dimension $n\geq 5$,  we have a more topological argument to show the  conclusion of Proposition \ref{gluing1} along the lines of proof of Lemma 3.2.2 in \cite{BBB+}; see the proof of Claim 1 in the proof of Proposition \ref{prop 3.5}.

\vspace *{0.2cm}

The following result is an application of Proposition \ref{gluing1} (and the Cerf-Palais disk theorem). It is an extension of Lemma 3.2.8 in \cite{BBB+}, and is used in the proof of Proposition \ref{prop 4.4}.

\begin{lem} \label{gluing1.5} (cf.  Lemma 3.2.8 in \cite{BBB+})  Let $n\geq 4$, and $0<\varepsilon \leq 2\varepsilon_1$, where  $\varepsilon_1$ is as in Proposition \ref{gluing1}. Suppose that $M$ is a connected Riemannian manifold of dimension $n$, and $\gamma:[0,1) \rightarrow M$ is a minimal geodesic   such that  $\gamma(s)$ is a center of an $\varepsilon$-neck in $M$ diffeomorphic to $\mathbb{S}^{n-1} \times (0,1)$ for each $s\in [0,1)$,  and $R(\gamma(s))\rightarrow \infty$ as $s\rightarrow 1$.  Then there exists an $\varepsilon$-horn in $M$ diffeomorphic to $\mathbb{S}^{n-1} \times (0,1)$  which contains the open end of the image of $\gamma$.
\end{lem}

\noindent {\bf Proof}.\ \ Let  $N_0$ be the $\varepsilon$-neck in $M$ centered at $\gamma(0)$, and $S$ be the central cross-section of the neck $N_0$.

\vspace *{0.2cm}

\noindent {\bf Claim}. The hypersurface $S$ separates the manifold $M$.

\vspace *{0.2cm}

\noindent {\bf Proof of Claim}. We argue by contradiction. Suppose the claim is not true.  Choose $s_0>0$ such that $\gamma(s_0) \in N_0$, and a point $x_1 \in N_0\setminus S$ such that $x_1$ and $\gamma(s_0)$ lie in different components of $N_0\setminus  S$. Choose a path $\alpha: [0,1]\rightarrow M\setminus S$ with $\alpha(0)=\gamma(s_0)$ and $\alpha(1)=x_1$. Note that there exists $s_1 \in (0,1)$ such that for any  $s\in [s_1,1)$, the intersection of the (image of the) path $\alpha$ with the central cross section  of the $\varepsilon$-neck centered at $\gamma(s)$ is empty, since  $R(\gamma(s))\rightarrow \infty$ as $s\rightarrow 1$  but the image of $\alpha$ is compact. Denote the central cross section  of the $\varepsilon$-neck centered at $\gamma(s_1)$ by $S_1$. We may assume that $S\cap S_1=\emptyset$ by choosing $s_1$ with $1-s_1$ small. By Proposition \ref{gluing1} and Lemma \ref{lem 4.1} the hypersurface $S$ and $S_1$ cobound a compact domain in $M$, denoted by $\Omega$, which is diffeomorphic to $\mathbb{S}^{n-1} \times [0,1]$ and contains $\gamma(s_0)$ but does not contain $x_1$. As the  intersection of the (image of the) path $\alpha$ with $\partial \Omega$ is empty, and $\alpha(0)=\gamma(s_0) \in \Omega$, it follows that the image of $\alpha$ is contained in $\Omega$. But on the other hand $\alpha(1)=x_1 \notin \Omega$.  A contradiction.   \hfill{$\Box$}

\vspace *{0.2cm}

Now we remove from  $M$  the component of $M\setminus S$ which does not contain $\gamma((0, 1))$, and get a manifold, denoted by $Y$, which contains $\gamma([0, 1))$ and has $S$ as its boundary. Glue smoothly a $n$-disk $D_0$ to $Y$ along their boundaries, and denote the manifold thus obtained by $Z$. Note that $Z$ is a connected, noncompact, smooth manifold without boundary. Using Proposition \ref{gluing1} and arguing as in the proof of Lemma 3.2.8 in \cite{BBB+} we see that $Z=\cup_{i=0}^\infty D_i$, where each $D_i$ is diffeomorphic to $D^n$, and $D_i \subset \text{Int} \hspace*{0.5mm} D_{i+1}$.  Now using the Cerf-Palais disk theorem (see  for example Theorem 3.1 in Chapter 8 of \cite{Hi} or Theorem 2.1 on p. 197 of \cite{M07}) one can show that $Z$ is diffeomorphic to $\mathbb{R}^n$  (see p. 126 of \cite{M07}). By the Cerf-Palais disk theorem again $Z\setminus D_0$ is diffeomorphic to $\mathbb{S}^{n-1} \times (0,1)$,  and the desired result follows.
 \hfill{$\Box$}

\vspace *{0.2cm}

Hamilton's canonical parametrization/uniformization of an $\varepsilon$-tube (see Section C.2 in \cite{H97}) is an important tool in some situations. Below is a very slight adaption of Hamilton's definition of a normal neck  in Section C.2 in \cite{H97}. Let $g_{\text{cyl}}$ be the standard  product metric on the cylinder $\mathbb{S}^{n-1}\times (a, b)$ with scalar curvature 1  and length $b-a$ for the interval $(a,b)$.  (Beware that our normalization of the product metric on the cylinder is slightly different from that of \cite{H97}, where the radius of the factor $\mathbb{S}^{n-1}$ is 1.)
Let $(\mathcal{O}, g)$ be a Riemannian $n$-orbifold and $\Phi: \mathbb{S}^{n-1}\times
(a,b)\rightarrow U \subset \mathcal{O}$  be a smooth map whose image $U$ is a  smooth
submanifold diffeomorphic to $\mathbb{S}^{n-1}/\Gamma
\times (a,b)$, where  $\Gamma$ is a  finite subgroup of $O(n)$ acting freely on $\mathbb{S}^{n-1}$. Suppose that $\Phi=\phi\circ \pi$, where $\pi: \mathbb{S}^{n-1}\times
(a,b)\rightarrow  \mathbb{S}^{n-1}/\Gamma \times (a,b)$ is the natural projection,   and $\phi:\mathbb{S}^{n-1}/\Gamma
\times (a,b) \rightarrow U $ is a  topological neck.  For $z\in (a,b)$, we define the ratio $r(z)$ of the mean radius  of the horizontal sphere $\mathbb{S}^{n-1} \times \{z\}$ in the pull-back metric $\Phi^*g$ (restricted to the sphere) to the radius of $\mathbb{S}^{n-1} \times \{z\}$ in the  metric $g_{\text{cly}}$ (restricted to the sphere) so that the area w.r.t. $\Phi^*g$ is
\begin{equation*}
A(\mathbb{S}^{n-1}   \times \{z\}, \Phi^*g)=A(\mathbb{S}^{n-1}\times \{z\}, g_{\text{cly}}|_{\mathbb{S}^{n-1}\times \{z\}}) r(z)^{n-1}.
\end{equation*}
 We'll call $\Phi$  Hamilton's canonical uniformization  if in addition it satisfies the
conditions (a), (b), (c) and (d) listed in Section C.2  of \cite{H97}; in this case we also call $\phi$ Hamilton's canonical parametrization. For the convenience of the readers, below we will recall these conditions from \cite{H97}.

\vspace *{0.2cm}

(a) every horizontal sphere $\mathbb{S}^{n-1} \times \{z\}$ for $z\in (a,b)$ has constant mean curvature in the pull-back metric $\Phi^*g$;

(b) the identity map from every horizontal sphere in the standard  metric $g_{\text{cly}}$ on the cylinder  (restricted to the sphere) to the same sphere in the   pull-back metric $\Phi^*g$ (also restricted to the sphere) is harmonic;

(c) the volume of any subcylinder in the pull-back metric $\Phi^*g$ is given by
\begin{equation}
\text{vol}(\mathbb{S}^{n-1}   \times [w,w'], \Phi^*g)=A(\mathbb{S}^{n-1}\times \{w\}, g_{\text{cly}}|_{\mathbb{S}^{n-1}\times \{w\}}) \int_w^{w'}r(z)^ndz;
\end{equation}
 and

 (d) if $\bar{V}$ is a Killing vector field on $\mathbb{S}^{n-1} \times \{z\}$ in the metric $g_{\text{cly}}$ restricted to the sphere, and  $W$ is a unit normal vector field  to the sphere in the metric $\Phi^*g$, then
\begin{equation*}
\int_{\mathbb{S}^{n-1} \times \{z\}}g_{\text{cly}}(\bar{V},W)d\bar{a}=0,
\end{equation*}
 where $d\bar{a}$ is the surface measure on the sphere induced by $g_{\text{cly}}$.

\vspace *{0.2cm}

\noindent {\bf Remark}.    Consider the identity map $\text{id}: (\mathbb{S}^{n-1} \times  (a,b), g_{\text{cly}}) \rightarrow (\mathbb{S}^{n-1} \times  (a,b), cg_{\text{cly}})$, where $c$ is an arbitrary positive constant. If we want the map $\text{id}$ to be a Hamilton's canonical parametrization, then the exponent appearing in the integrand on the RHS of the equality (A.1) should be $n$ instead of $n-1$.   Compare Remark 3.10 (iii) in \cite{HS}.

\vspace *{0.2cm}

The following lemma is adapted from  Theorem C2.2 and Corollary C2.3 in \cite{H97}, and is a  preparation for  Proposition \ref{gluing3}.

\begin{lem} \label{gluing2} (cf. Theorem C2.2 and Corollary C2.3 in \cite{H97})
 Let $n\geq 4$. There exists a positive constant $\varepsilon_2=\varepsilon_2(n)$  such that  for any positive number $\varepsilon \leq \varepsilon_2$,  there is  a positive number $\tilde{\varepsilon}_1=\tilde{\varepsilon}_1(\varepsilon)(=\tilde{\varepsilon}_1(\varepsilon,n))\leq \varepsilon$  with the following property. Let $(\mathcal{O},g)$ be a Riemannian orbifold of dimension $n$ with at most isolated singularities, and $N$ be an $\eta$-neck  in $(\mathcal{O},g)$ given by a diffeomorphism
 $\psi: \mathbb{S}^{n-1}/\Gamma
  \times (-\eta^{-1}, \eta^{-1})  \rightarrow N$ with  $0<\eta \leq \tilde{\varepsilon}_1$. Then there exists a
  Hamilton's canonical parametrization $\phi:\mathbb{S}^{n-1}/\Gamma \times (a, b)\rightarrow \mathcal{T}\subset N$ which gives  an $\varepsilon$-tube $\mathcal{T}$ in $N$  with $\psi( \mathbb{S}^{n-1}/\Gamma
  \times (-\eta^{-1}+1, \eta^{-1}-1))\subset \mathcal{T}$, in the sense that for any $(c,d)\subset (a,b)$ with $d-c=2\varepsilon^{-1}$, $\phi|_{\mathbb{S}^{n-1}/\Gamma\times (c,d)}$ gives an $\varepsilon$-neck structure (up to a translation of the parameter) for its image.
\end{lem}

\noindent {\bf Proof}.\ \  Let $0< \varepsilon < \frac{1}{10^4n} $ be given. Suppose that $N$ is an $\eta$-neck with scaling factor $Q$ in $(\mathcal{O},g)$ given by a diffeomorphism  $\psi: \mathbb{S}^{n-1}/\Gamma
  \times (-\eta^{-1}, \eta^{-1})  \rightarrow N$  for some $\eta>0$. Consider  the Riemannian universal cover $\pi:\tilde{N}\rightarrow N$  of the neck $N$,  which has an induced $\eta$-neck structure given by a lift $\tilde{\psi}: \mathbb{S}^{n-1}
  \times (-\eta^{-1}, \eta^{-1})  \rightarrow \widetilde{N}$ of $\psi$. Let $\widetilde{\Gamma}$ be the deck transformation group of the Riemannian covering
  $\pi:\tilde{N}\rightarrow N$. We have $\widetilde{\Gamma} \cong \Gamma$.  Note that for any $z\in (-\eta^{-1}, \eta^{-1})$, $\tilde{\psi}( \mathbb{S}^{n-1}
  \times \{z\})$ is $\widetilde{\Gamma}$-invariant.  By  (the proof of) Theorem C2.2 and Corollary C2.3   in  \cite{H97}, there exists  a positive number  $\tilde{\varepsilon}_1=\tilde{\varepsilon}_1(\varepsilon)$ with the following property. If $0< \eta \leq  \tilde{\varepsilon}_1$, there exists a  Hamilton's canonical parametrization $\tilde{\phi}:\mathbb{S}^{n-1}\times (a, b)\rightarrow \widetilde{\mathcal{T}} \subset \widetilde{N} $ which gives  an $\varepsilon$-tube $\widetilde{\mathcal{T}}$ in $\widetilde{N}$  with $\tilde{\psi}( \mathbb{S}^{n-1}
  \times (-\eta^{-1}+1, \eta^{-1}-1))\subset \widetilde{\mathcal{T}}$, in the sense that for any $(c,d)\subset (a,b)$ with $d-c=2\varepsilon^{-1}$, $\tilde{\phi}|_{\mathbb{S}^{n-1}\times (c,d)}$ gives an $\varepsilon$-neck structure (up to a translation of the parameter) for its image, such that for each $s\in (a,b)$, $\tilde{\phi}(\mathbb{S}^{n-1}
  \times \{s\})$ is a constant mean curvature hypersurface which is close (with error controlled in terms of $\varepsilon$) to horizontal w.r.t. the $\eta$-neck structure of $\widetilde{N}$. The closeness of the constant mean curvature hypersurface $\tilde{\phi}(\mathbb{S}^{n-1}
  \times \{s\})$ to the corresponding horizontal hypersurface (w.r.t. the $\eta$-neck structure on $\widetilde{N}$) can be described more precisely as follows. Given any $s\in (a,b)$, and $\tilde{x}_0\in \tilde{\phi}(\mathbb{S}^{n-1}\times \{s\})$, let $\widetilde{\Sigma}_0$ be the cross section of the neck $\widetilde{N}$ passing through $\tilde{x}_0$, then we have $\tilde{\phi}(\mathbb{S}^{n-1}\times \{s\})=\{ \exp_{\tilde{x}}(Q^{-1/2}\tilde{f}(\tilde{x})\tilde{\nu}(\tilde{x})) \hspace*{1mm}|\hspace*{1mm}  \tilde{x}\in \widetilde{\Sigma}_0 \}$, where $\exp$ is the exponential map of the Riemannian manifold $\widetilde{N}$, $\tilde{\nu}$ is a unit normal vector field on $\widetilde{\Sigma}_0$, and $\tilde{f}$ is certain smooth function on $\widetilde{\Sigma}_0$ with $\tilde{f}(\tilde{x}_0)=0$ and $||\tilde{f}||_{C^{4,\frac{1}{2}}}< \frac{1}{2}\varepsilon$. Here the $C^{4,\frac{1}{2}}$-norm is taken w.r.t. the round metric on $\widetilde{\Sigma}_0$ with scalar curvature 1. See also the proof of Proposition D.1 in version 1 of \cite{BN2} on arXiv. We take $\tilde{\phi}$ to be maximal in the sense that $\widetilde{\mathcal{T}}$ cannot be properly included  in the image of any  Hamilton's canonical parametrization which gives an $\varepsilon$-tube in  $\widetilde{N}$.

  Given any $ \tilde{\gamma}\in \widetilde{\Gamma}$, which is  an isometry of $\widetilde{N}$, $\tilde{\gamma} \tilde{\phi}$ is a Hamilton's canonical parametrization for
  $\tilde{\gamma}(\widetilde{\mathcal{T}})\subset \widetilde{N}$. If $\varepsilon \leq \varepsilon_2=\varepsilon_2(n)$,  where $\varepsilon_2(n)< 0.1$ is a positive constant depending only on the dimension $n$, by our choice of  $\tilde{\phi}$ and Theorem C2.4 in \cite{H97} we have $\tilde{\gamma}(\widetilde{\mathcal{T}})\subset \widetilde{\mathcal{T}}$. If $\varepsilon_2=\varepsilon_2(n)$ is sufficiently small and $\varepsilon \leq \varepsilon_2$, by Theorem C2.4 and Lemma C2.1 in \cite{H97} there  exists an isometry $\hat{\gamma}$ of the standard round cylinder such that  $\hat{\gamma}(\mathbb{S}^{n-1}\times (a, b)) \subset  \mathbb{S}^{n-1}\times (a, b)$    and  $\tilde{\gamma} \tilde{\phi}=\tilde{\phi} \hat{\gamma}|_{\mathbb{S}^{n-1}\times (a, b)}$.

   As in our current situation both $a$ and  $b$ are finite numbers,   the $\mathbb{R}$-factor of $\hat{\gamma}$  cannot be a nontrivial translation.
We would like to point out that there is another argument for this fact without using the finiteness of $a$ and $b$.  We argue by contradiction. Suppose otherwise, that is  the $\mathbb{R}$-factor of $\hat{\gamma}$ is a translation $s\mapsto s+L$ with $L\neq 0$. Fix a point $\tilde{x}\in  \widetilde{\mathcal{T}}$. Then $\tilde{\gamma}^k\tilde{x} \in \widetilde{\mathcal{T}}$, $k=1,2,\cdot\cdot\cdot$. Note that  the $\mathbb{R}$-coordinate of the point $\tilde{\gamma}^k\tilde{x}$  w.r.t. the parametrization $\tilde{\phi} $ is equal to the $\mathbb{R}$-coordinate of the point $\tilde{x}$ w.r.t.  $\tilde{\phi} $ plus $kL$, $k=1,2,\cdot\cdot\cdot$.  In particular, $\tilde{\gamma}^k\tilde{x}\neq \tilde{x}$ for any positive integer $k$. But on the other hand, $\widetilde{\Gamma}$ is a finite group, so there exists a positive integer $m$ such that $\tilde{\gamma}^m=Id$. It follows that $\tilde{\gamma}^m\tilde{x}=\tilde{x}$. A contradiction.

   The $\mathbb{R}$-factor of $\hat{\gamma}$  cannot be a reflection either. We argue by contradiction. Otherwise, as
   $\hat{\gamma}(\mathbb{S}^{n-1}\times (a, b)) \subset  \mathbb{S}^{n-1}\times (a, b)$, and both $a$ and $b$ are finite numbers, the   $\mathbb{R}$-factor of $\hat{\gamma}$ must be the reflection  w.r.t. the point $\frac{a+b}{2}$. Fix  any point $p \in \mathbb{S}^{n-1}\times (a, b)$ whose $\mathbb{R}$-component  is close to $a$, the $\mathbb{R}$-component of the point $\hat{\gamma}p$ will be close to   $b$.   The  distance of  the points $\tilde{\phi}(p)$ and $\tilde{\phi}(\hat{\gamma}p)$ (computed in $\widetilde{\mathcal{T}}$) must be much larger than the diameter of the  cross section, say $\widetilde{\Sigma}$, of the neck $\widetilde{N}$ passing through the point $\tilde{\phi}(p)$. On the other hand, the distance of the points $\tilde{\phi}(p)$ and $\tilde{\gamma}  \tilde{\phi}(p)$ (computed in $\widetilde{\mathcal{T}}$) is not larger than the diameter of $\widetilde{\Sigma}$.  But $\tilde{\phi}(\hat{\gamma}p)=\tilde{\gamma} \tilde{\phi}(p)$. Thus we get  a contradiction. (By the way, if in some other situation  exactly one of  $a$ and $b$ is
    finite, since $\hat{\gamma}(\mathbb{S}^{n-1}\times (a, b)) \subset  \mathbb{S}^{n-1}\times (a, b)$, the $\mathbb{R}$-factor of $\hat{\gamma}$  cannot be a reflection.)

   Thus  we see that in effect, the group $\hat{\Gamma}:=\{\hat{\gamma}|_{\mathbb{S}^{n-1}\times (a, b)}  \hspace{1mm} | \hspace{1mm} \tilde{\gamma} \in \widetilde{\Gamma}\}=\{{\tilde{\phi}}^{-1}\tilde{\gamma}\tilde{\phi}  \hspace{1mm}  | \hspace{1mm} \tilde{\gamma} \in \widetilde{\Gamma}\} $   only acts on the $\mathbb{S}^{n-1}$ factor of $\mathbb{S}^{n-1}   \times (a, b)$, and the
  parametrization $\tilde{\phi}$, which is equivariant w.r.t. the actions of $\hat{\Gamma}$ and $\widetilde{\Gamma}$,  can
be pushed down  to get a Hamilton's canonical parametrization   $\phi:\mathbb{S}^{n-1}/\widehat{\Gamma} \times (a, b)\rightarrow \mathcal{T}\subset N$ which  gives  an $\varepsilon$-tube $\mathcal{T}$ in $N$  with $\psi( \mathbb{S}^{n-1}/\Gamma
  \times (-\eta^{-1}+1, \eta^{-1}-1))\subset \mathcal{T}$, in the sense that for any $(c,d)\subset (a,b)$ with $d-c=2\varepsilon^{-1}$, $\phi|_{\mathbb{S}^{n-1}/\widehat{\Gamma}\times (c,d)}$ gives an $\varepsilon$-neck structure (up to a translation of the parameter) for its image. (Compare the first paragraph of Section D.1 in \cite{H97}.)

Let $s\in (a,b)$,  $\tilde{x}_0\in \tilde{\phi}(\mathbb{S}^{n-1}\times \{s\})$, and $\widetilde{\Sigma}_0$  be the cross section of the neck $\widetilde{N}$ passing through $\tilde{x}_0$ as above. Let $x_0=\pi(\tilde{x}_0)$,   and  $\Sigma_0=\pi(\widetilde{\Sigma}_0)$. As both the ambient Riemannian manifold $\widetilde{N}$ and the hypersurfaces $\widetilde{\Sigma}_0$ and $\tilde{\phi}(\mathbb{S}^{n-1} \times \{s\})$ are $\widetilde{\Gamma}$-invariant, the function $\tilde{f}$ above must also be $\widetilde{\Gamma}$-invariant, so it descends to a smooth function $f$ on $\Sigma_0$ with   $f(x_0)=0$ and $||f||_{C^{4,\frac{1}{2}}}< \frac{1}{2}\varepsilon$  such that $\phi(\mathbb{S}^{n-1}/\widehat{\Gamma}\times \{s\})=\{ \exp_x(Q^{-1/2}f(x)\nu(x)) \hspace*{1mm}|\hspace*{1mm}  x\in \Sigma_0 \}$. Here the $C^{4,\frac{1}{2}}$-norm is taken w.r.t. the round metric on $\Sigma_0$ with scalar curvature 1, $\exp$ is the exponential map of the Riemannian manifold $(N, g)$, and $\nu$ is a unit normal vector field on $\Sigma_0$.  Now we see that $\phi(\mathbb{S}^{n-1}/\widehat{\Gamma}\times \{s\})$ is diffeomorphic to $\Sigma_0$.  Note that $\Sigma_0$ is diffeomorphic to $\mathbb{S}^{n-1}/\Gamma$. It follows that $\widehat{\Gamma}$ is conjugate to $\Gamma$ in $O(n)$, and we can view  the Hamilton's canonical parametrization $\phi:\mathbb{S}^{n-1}/\widehat{\Gamma} \times (a, b)\rightarrow \mathcal{T}$  as a map from
$\mathbb{S}^{n-1}/\Gamma \times (a, b)$ to $ \mathcal{T}$  which is still a Hamilton's canonical parametrization.
 \hfill{$\Box$}

\vspace *{0.2cm}

The following result  is based on Lemma \ref{gluing2} above and Theorem C2.4 in \cite{H97}. For relevant discussion in dimension $n=4$ see Appendix B in \cite{Hu13} and Proposition 2.3 in \cite{Hu15}.

\begin{prop} \label{gluing3}  (cf.  Theorem C2.4 in \cite{H97})
 Let $n\geq 4$. There exists a positive constant $\varepsilon_3=\varepsilon_3(n) \leq \varepsilon_2$, where $\varepsilon_2$ is as in Lemma \ref{gluing2}, with the following property. Let $0<\varepsilon \leq \varepsilon_3$,  and $\tilde{\varepsilon}_1=\tilde{\varepsilon}_1(\varepsilon)$, where $\tilde{\varepsilon}_1(\cdot)$ is as in Lemma \ref{gluing2}. Let $(\mathcal{O},g)$ be a Riemannian orbifold of dimension $n$ with at most isolated singularities, $0< \eta \leq \tilde{\varepsilon}_1$, and $N_i$ be an $\eta$-neck centered at a point $x_i$ in $\mathcal{O}$ given by a diffeomorphism $\psi_i: \mathbb{S}^{n-1}/\Gamma_i
  \times (-\eta^{-1}, \eta^{-1})  \rightarrow N_i$, $i=1,2$, such that  $\psi_1(\mathbb{S}^{n-1}/\Gamma_1
  \times (-\eta^{-1}+2, \eta^{-1}-2))  \cap \psi_2(\mathbb{S}^{n-1}/\Gamma_2
  \times (-\eta^{-1}+1, \eta^{-1}-1))  \neq \emptyset$.
   Then $\Gamma_1$ and $\Gamma_2$ are conjugate in $O(n)$.  Moreover, if   $\psi_1(\mathbb{S}^{n-1}/\Gamma_1
  \times (-\eta^{-1}, -\eta^{-1}+0.01))  \cap  N_2 = \emptyset$  and   $\psi_2(\mathbb{S}^{n-1}/\Gamma_2
  \times (\eta^{-1}-0.01, \eta^{-1}))  \cap  N_1 = \emptyset$, there exists a Hamilton's canonical parametrization $\phi: \mathbb{S}^{n-1}/\Gamma_1 \times (a,b) \rightarrow \mathcal{T}\subset N_1\cup N_2$  which gives  an $\varepsilon$-tube $\mathcal{T}$ in $N_1\cup N_2$  with
  $\psi_1( \mathbb{S}^{n-1}/\Gamma_1   \times (-\eta^{-1}+1, \eta^{-1}))  \cup
    \psi_2( \mathbb{S}^{n-1}/\Gamma_2   \times (-\eta^{-1}, \eta^{-1}-1))      \subset \mathcal{T}$.
  \end{prop}

\noindent {\bf Proof}.\ \  Given $0< \varepsilon \leq \varepsilon_2$ and $\eta$-neck $N_i$ ($i=1,2$) in $(\mathcal{O},g)$ with $\eta\leq \tilde{\varepsilon}_1(\varepsilon)$,  by Lemma \ref{gluing2}, there exists a
  Hamilton's canonical parametrization $\phi_i:\mathbb{S}^{n-1}/\Gamma_i \times (a_i, b_i)\rightarrow \mathcal{T}_i\subset N_i$ which gives  an $\varepsilon$-tube $\mathcal{T}_i$ in $N_i$  with $\psi_i( \mathbb{S}^{n-1}/\Gamma_i
  \times (-\eta^{-1}+1, \eta^{-1}-1))\subset \mathcal{T}_i$. We can take $\phi_i$ to be maximal in the sense that $\mathcal{T}_i$ cannot be properly included  in the image of any  Hamilton's canonical parametrization which gives an $\varepsilon$-tube in  $N_i$, $i=1,2$. Let $\Phi_i= \phi_i \circ \pi_i:\mathbb{S}^{n-1} \times (a_i, b_i)\rightarrow \mathcal{T}_i\subset N_i$ be Hamilton's canonical uniformization of $\mathcal{T}_i$, where $\pi_i:\mathbb{S}^{n-1} \times (a_i, b_i)  \rightarrow
  \mathbb{S}^{n-1}/ \Gamma_i \times (a_i, b_i)$ is the natural projection. Choose a point $x_0 \in \psi_1(\mathbb{S}^{n-1}/\Gamma_1
  \times (-\eta^{-1}+2, \eta^{-1}-2))  \cap \psi_2(\mathbb{S}^{n-1}/\Gamma_2
  \times (-\eta^{-1}+1, \eta^{-1}-1))$. Then there exists $s_1 \in  (a_1,b_1)$ away from $a_1$ and $b_1$, and $s_2 \in  (a_2,b_2)$, such that $x_0 \in \Phi_1(\mathbb{S}^{n-1} \times \{s_1\})
 \cap \Phi_2(\mathbb{S}^{n-1} \times \{s_2\})$.  If $\varepsilon_3=\varepsilon_3(n)\leq \varepsilon_2$ is sufficiently small and $0< \varepsilon \leq \varepsilon_3$,
by  the proof of   Theorem C2.4 in \cite{H97}, which uses Lemma C2.1 in \cite{H97},   there  exists an isometry of the round cylinder  $\mathbb{S}^{n-1}\times \mathbb{R}$, say $G$,  such that $G(\mathbb{S}^{n-1} \times \{s_2\})= \mathbb{S}^{n-1} \times \{s_1\}$ and $\Phi_1G|_{\mathbb{S}^{n-1} \times \{s_2\}}=\Phi_2|_{\mathbb{S}^{n-1} \times \{s_2\}}$. So $\phi_1(\mathbb{S}^{n-1}/\Gamma_1 \times \{s_1\})= \phi_2(\mathbb{S}^{n-1}/\Gamma_2 \times \{s_2\})$.  It follows that $\Gamma_1$ and $\Gamma_2$ are conjugate in $O(n)$.

If   $\psi_1(\mathbb{S}^{n-1}/\Gamma_1
  \times (-\eta^{-1}, -\eta^{-1}+0.01))  \cap  N_2 = \emptyset$  and   $\psi_2(\mathbb{S}^{n-1}/\Gamma_2
  \times (\eta^{-1}-0.01, \eta^{-1}))  \cap  N_1 = \emptyset$, as in   the proof of   Theorem C2.4 in \cite{H97}, we can glue the domain of definition of $\Phi_1$ and that of $\Phi_2$ together using the isometry $G$ above, and get a Hamilton's canonical uniformization $\Phi:\mathbb{S}^{n-1} \times (a, b)\rightarrow \mathcal{T}_1 \cup \mathcal{T}_2 \subset N_1 \cup N_2$. From $\Phi$ we get a Hamilton's canonical parametrization $\phi$ as desired.
\hfill{$\Box$}

\vspace *{0.2cm}

For $n\geq 3$, let $(\mathbb{R}^n,\hat{g}(t))$, $t\in [0,\frac{n-1}{2})$, be
 a smooth standard solution constructed in for example, \cite{P2},  \cite{KL08}, Chapter 7 of  \cite{BBB+}, Chapter 12 of \cite{MT}, Appendix A in \cite{CZ} (but beware that the formula for the function $f$ near the point $z=4$ on p. 260 of \cite{CZ} is not correct), and Theorem 9.1 in \cite{B19} (the $\hat{g}(t)$ here is denoted by $g(t)$ there), which is a rotationally symmetric (i.e. $O(n)$-invariant) solution to the Ricci flow with nonnegative curvature operator on $\mathbb{R}^n$, and let $p_0$ be its tip, which is the (unique) point fixed by $O(n)$. Recall that the initial metric $\hat{g}(0)$ of a smooth standard solution is obtained by gluing a rotationally symmetric cap, which is diffeomorphic to $D^n$ and called the surgery cap, to the half-infinite round cylinder $\mathbb{S}^{n-1}\times (-\infty, 0]$ with scalar curvature 1 along their boundaries.  More precisely, we proceed as in Section 72 of  \cite{KL08}. (Compare Chapter 7 of  \cite{BBB+}.) First note that for a smooth function $\varphi$ on $\mathbb{R}$, the  sectional curvatures of the metric  $e^{-2\varphi} g_{\text{cly}}$ are $e^{2\varphi}\varphi''$ and $e^{2\varphi}(\frac{1}{(n-1)(n-2)}-(\varphi')^2)$, where $g_{\text{cly}}$ denotes the standard product metric  on the  cylinder  $\mathbb{S}^{n-1}\times \mathbb{R}$ with scalar curvature 1.  For any positive number $k$, let
  \begin{equation*}
\varphi_k(z)=\log (1+e^{c_nz})-\frac{c_nz}{2}-\log \frac{c_n}{k},   \hspace*{8mm}  z\in \mathbb{R},
\end{equation*}
where $c_n=\frac{2}{\sqrt{(n-1)(n-2)}}$.  One can check that the metric $e^{-2\varphi_k} g_{\text{cly}}$ has constant sectional curvature $k^2$, and extends to a smooth metric on $\mathbb{S}^n$ of constant sectional curvature $k^2$ by adding one point  for each of the two ends.

  As in Section 72 of  \cite{KL08}, there exist numbers $b_0$ and $b'$ with  $\frac{1}{10} < b' < b_0$, and a smooth convex function $\varphi$ defined on $(-\infty, b_0]$ with $\varphi|_{(-\infty, 0]}=0$  and  $\varphi(z)= ce^{-1/z}$ for $z \in (0, \frac{1}{10}]$, where $c< \frac{c_n}{2}$ is a small positive constant, such that $\varphi$ coincides with $\varphi_k$ on $[b',b_0]$ for some $k>0$, and the metric $e^{-2\varphi} g_{\text{cly}}|_{\mathbb{S}^{n-1}\times (-\infty,b_0]}$ has nonnegative curvature operator. Here, we use the fact that if $c< \frac{c_n}{2}$,
  \begin{equation*}
  \frac{d}{dz} ce^{-1/z}|_{z=1/10} < \lim_{z\rightarrow \infty} \varphi_k'(z)=\frac{1}{\sqrt{(n-1)(n-2)}}.
   \end{equation*}
   In particular,  the metric $e^{-2\varphi} g_{\text{cly}}|_{\mathbb{S}^{n-1}\times (-\infty,b_0]}$  has constant  sectional curvature $k^2$ on the region $\mathbb{S}^{n-1}\times [b',b_0]$. Then we smoothly glue a closed metric $n$-ball of the same constant  sectional curvature $k^2$  to $(\mathbb{S}^{n-1}\times (-\infty,b_0], e^{-2\varphi} g_{\text{cly}}|_{\mathbb{S}^{n-1}\times (-\infty,b_0]})$ along their boundaries, and the resulting Riemannian manifold is a possible choice for $(\mathbb{R}^n,\hat{g}(0))$. Let  $\Gamma$ be a finite subgroup of $O(n)$ acting freely on $\mathbb{S}^{n-1}$,  and $(\mathbb{R}^n//\Gamma,\hat{g}_\Gamma(t))$ be the quotient of a smooth standard solution  by $\Gamma$, called an orbifold standard solution. So the initial metric $\hat{g}_\Gamma(0)$ of an orbifold standard solution is obtained by gluing a (quotient) cap which is diffeomorphic to $D^n//\Gamma$ and also called the surgery cap, to the half-infinite quotient cylinder $\mathbb{S}^{n-1}/\Gamma \times (-\infty, 0]$ with scalar curvature 1 along their boundaries.  Let  $p_\Gamma$ be the tip
of  an orbifold standard solution $(\mathbb{R}^n//\Gamma, \hat{g}_\Gamma(t))$, that is the image of the point $p_0$ under the natural projection  $\pi_{\Gamma}: \mathbb{R}^n \rightarrow
\mathbb{R}^n//\Gamma$.

By the maximum principle $R(x,t)= R_{\hat{g}(t)}(x)\geq \frac{1}{1-\frac{2}{n}t}$ for any $x\in \mathbb{R}^n$.  So if $t> \frac{n}{n+2}$, $tR(x,t)\geq t\frac{1}{1-\frac{2}{n}t} >1$. We'll fix an $n$-dimensional smooth standard solution. Let  $K_{\text{st}}$ be the superemum of the norm of the Riemannian curvature tensor of the
$n$-dimensional smooth standard solution on the time interval $[0,\frac{n+1}{n+2}]$.

The
following lemma on strengthening necks extends and improves Lemma
4.3.5 in \cite{BBB+} and Lemma 4.11 in \cite{BBM}; for a version in the 4-dimensional case see also Lemma A.2 in \cite{Hu13} and Lemma 2.5 in \cite{Hu15}.
It is used several times  in this paper.

\begin{lem} \label{lem 2.5} (cf.  Lemma
4.3.5 in \cite{BBB+})  Let $n\geq 3$, and $K=\max \{K_{\text{st}}, 1\}$.
  For any $\varepsilon \in (0, 10^{-3})$
 there exists a number $\beta=\beta (\varepsilon)(=\beta (\varepsilon, n))\in (0,10^{-2})$ with the
 following property.
  Let $a, b$ be real numbers with $a< b \leq 0$, $|b| \leq
  \frac{n}{n+1}$,  and $(\mathcal{O}(\cdot), g(\cdot))$ be an $n$-dimensional complete surgical solution to the Ricci flow
  defined on $(a,0]$  with at most isolated singularities, and $x$ be a point in $\mathcal{O}(b)$ such that:

  (i) $R(x,b)=1$;

  (ii) $(x,b)$ is a center of a  $\beta \varepsilon$-neck with a scaling factor $Q$;

  (iii) the space-time subset $B(x, b, (\beta \varepsilon)^{-1}) \times  [b-Q^{-1}, 0]$ is unscathed with $|Rm|\leq 2K$.

 \noindent  Then $(x,0)$ is a center of a strong $\varepsilon$-neck.
\end{lem}

\noindent {\bf Proof}.  The proof is adapted from that of   Lemma
4.3.5 in \cite{BBB+}. We  argue by contradiction. Otherwise there exists
$\varepsilon \in (0, 10^{-3})$, a sequence of positive numbers $\beta_k \rightarrow 0$,
sequences $a_k < b_k$, $b_k \in [-\frac{n}{n+1},0]$,   and a sequence of $n$-dimensional
surgical solutions $(\mathcal{O}_k(t), g_k(t))$ ($t \in (a_k,0]$) with a point
$x_k \in \mathcal{O}_k(b_k)$ such that

(i) $R(x_k,b_k)=1$;

  (ii) $(x_k,b_k)$ is a center of a  $\beta_k \varepsilon$-neck $N_k$ given by a diffeomorphism $\psi_k: \mathbb{S}^{n-1}/\Gamma_k \times
(-(\beta_k \varepsilon)^{-1},(\beta_k \varepsilon)^{-1})\rightarrow  N_k$ with a scaling factor $Q_k$;

  (iii) the space-time subset $B(x_k, b_k, (\beta_k \varepsilon)^{-1}) \times  [b_k-Q_k^{-1}, 0]$ is unscathed with $|\text{Rm}|\leq 2K$, but

  (iv) $(x_k,0)$ is not a center of any strong $\varepsilon$-neck.

Let $N_k'= \psi_k(\mathbb{S}^{n-1}/\Gamma_k \times
(-0.99 (\beta_k \varepsilon)^{-1}, 0.99 (\beta_k \varepsilon)^{-1}))$. From (i) (ii) (iii) above we have $Q_k\rightarrow 1$, and $|\text{Rm}|\leq 2K$ on $N_k' \times [b_k-Q_k^{-1}, b_k]$ for any large $k$.  By (iii) and the distance distortion inequality (see Theorem 17.1 in \cite{H95b}) there exists a positive constant $c$ such that $B(x_k,0,c(\beta_k\varepsilon)^{-1})\subset
B(x_k,b_k,(\beta_k\varepsilon)^{-1})$ for any $k$. So $|\text{Rm}|\leq 2K$ on  $P(x_k, 0, c(\beta_k \varepsilon)^{-1}, b_k-Q_k^{-1})$.

For each $k$, let
\begin{equation*}
\Psi_k=\psi_k \circ \pi_k: \mathbb{S}^{n-1} \times (-0.99 (\beta_k \varepsilon)^{-1}, 0.99 (\beta_k \varepsilon)^{-1})
\rightarrow  N_k',
\end{equation*}
where
\begin{equation*}
\pi_k: \mathbb{S}^{n-1}\times
(-0.99 (\beta_k \varepsilon)^{-1},0.99 (\beta_k \varepsilon)^{-1})\rightarrow \mathbb{S}^{n-1}/\Gamma_k \times (-0.99 (\beta_k \varepsilon)^{-1}, 0.99(\beta_k \varepsilon)^{-1})
\end{equation*}
  is the natural projection. Now we  restrict the solutions $g_k(\cdot)$ to the space-time regions $N_k' \times [b_k-Q_k^{-1}, 0] $, and pull them back  to  $\mathbb{S}^{n-1} \times
(-0.99(\beta_k \varepsilon)^{-1}, 0.99 (\beta_k \varepsilon)^{-1})$  via $\Psi_k$. Choose a point  $*\in \mathbb{S}^{n-1}\times \{0\}$.  Using Theorem B.1.2 in \cite{BBB+} we get a positive lower bound (independent of $k$) for the injectivity radius of the pullback solutions at the space-time point $(*, 0)$. Up to an extraction, the sequence $b_k$ converges to some number $b\in [-\frac{n}{n+1},0]$. By using a local  version of Hamilton's compactness theorem for the Ricci flow \cite{H95}  (cf. Theorem 4.1.5 in \cite{CaZ}, Theorem 16.1 in \cite{H95b}, and Appendix E in \cite{KL08}) and the assumption   $\beta_k \rightarrow 0$,  we see that up to  a further extraction,  the sequence of pointed Ricci flow solutions $(B_{\Psi_k^*g_k(0)}(*, c(\beta_k \varepsilon)^{-1}), \Psi_k^*g_k(\cdot), (*,0)) $  converges smoothly to a complete pointed Ricci flow $(M_\infty, g_\infty(\cdot), (x_\infty,0))$ with bounded curvature, defined on the time interval $(b-1,0]$. In particular, the sequence of pointed Riemannian manifolds $(B_{\Psi_k^*g_k(0)}(*, c(\beta_k \varepsilon)^{-1}), \Psi_k^*g_k(b_k), *) $ subconverges smoothly to $(M_\infty, g_\infty(b), x_\infty)$.

 On the other hand, from the assumption (ii) we see that the sequence of  pointed Riemannian manifolds $(B_{\Psi_k^*g_k(0)}(*, c(\beta_k \varepsilon)^{-1}), Q_k\Psi_k^*g_k(b_k), *) $ converges smoothly to $(\mathbb{S}^{n-1}\times \mathbb{R}, g_{\text{cly}}(0), *)$, where $g_{\text{cly}}(\cdot)$ denotes the round cylinder solution to the Ricci flow on $\mathbb{S}^{n-1}\times \mathbb{R}$ with scalar curvature 1 at time 0, which is called the cylindrical flow. As $Q_k\rightarrow 1$, it follows that $(M_\infty, g_\infty(b))$ is isometric to  $(\mathbb{S}^{n-1}\times \mathbb{R}, g_{\text{cly}}(0))$. By the uniqueness theorem (see \cite{CZ06}) and the backwards uniqueness theorem (see \cite{Ko}) for the Ricci flow we see that
 $(M_\infty, \{g_\infty(t)\}_{t\in (b-1,0]})$ is isometric to  the  flow  $(\mathbb{S}^{n-1}\times \mathbb{R}, \{g_{\text{cly}}(t-b)\}_{t\in (b-1,0]})$.

 Now we see that the sequence $(B_{\Psi_k^*g_k(0)}(*, c(\beta_k \varepsilon)^{-1}), \{\Psi_k^*g_k(t)\}_{t\in (b_k-Q_k^{-1},0]}, (*,0)) $  subconverges  smoothly to $(\mathbb{S}^{n-1}\times \mathbb{R}, \{\tilde{g}(t)\}_{t\in (b-1,0]}, (*,0))$, where   $\tilde{g}(t)=g_{\text{cly}}(t-b)$.
 Using  assumption (ii), $Q_k\rightarrow 1$,  and $b_k\rightarrow b$  we deduce  that the  tensor $\Psi_k^*g_k(b)-\tilde{g}(b)$ together with  its the covariant derivatives are tending to 0 on arbitrarily large balls centered at $(*,b)$.
Then by the argument in Section 2 of \cite{H95}, we see that the diffeomorphisms involved in the subconvergence here can be chosen to be  the identity maps.
 So if $k$ is sufficiently large, the evolving tensor $\Psi_k^*g_k(\cdot)-\tilde{g}(\cdot)$ together with  all its covariant derivatives are arbitrarily close to 0 on  arbitrary compact subsets of $\mathbb{S}^{n-1}\times \mathbb{R} \times (b-1,0]$ containing $(*,0)$.   As $\Psi_k$ is $\Gamma_k$-invariant, it follows from Remark 4.2.7(i) in \cite{BBB+}
    and Lemma \ref{lem 2.7}  that (a suitable restriction of) the diffeomorphism $\psi_k$ will give a strong $\varepsilon$-neck centered at $(x_k,0)$ for sufficiently large $k$. (By the way, Remark 4.2.7(i) in \cite{BBB+} follows from the fact that given  $n>1$, $\varepsilon >0$ and $Q>0$, there exists
 $\hat{Q} > Q^{-1}$ such that
  \begin{equation*}
\sup_{t\in [-1,0]} |\hat{Q}Qg_{\text{cyl}}((\hat{Q}Q)^{-1}t)-g_{\text{cyl}}(t)|_{[\varepsilon^{-1}],\mathbb{S}^{n-1}\times (-\frac{1}{\varepsilon}, \frac{1}{\varepsilon}), g_{\text{cyl}}(t)} < \varepsilon,
\end{equation*}
 using the notation in Definition  2.2.4 in \cite{BBB+}.)  This contradicts assumption (iv).  \hfill{$\Box$}

\vspace *{0.2cm}
Note that the assumption of our Lemma \ref{lem 2.5} is slightly weaker than that of   Lemma
4.3.5 in \cite{BBB+}. That is, in assumption (ii) we do not assume that the $\beta \varepsilon$-neck is strong, which is assumed in Lemma
4.3.5(ii) in \cite{BBB+}. (Of course, we need to adjust assumption (iii) compared to Lemma
4.3.5(iii) in \cite{BBB+}.)


\vspace *{0.4cm}

Laboratory of Mathematics and Complex Systems (Ministry of Education),

School of Mathematical Sciences, Beijing Normal University,

Beijing 100875,  People's Republic of China

 E-mail address: hhuang@bnu.edu.cn


\begin{thebibliography}{99}


\bibitem {Ba} W. Ballmann,  Vector bundles and connections, unpublished lecture notes available at Prof. Ballmann's homepage.

\bibitem {B07}  R. Bamler, Ricci flow with surgery, Diploma thesis,   Ludwig-Maximilians-University, Munich 2007.



\bibitem {B} R. Benedetti, Lectures on differential topology.  Graduate Studies in Mathematics, 218. American Mathematical Society, Providence, RI, 2021.


\bibitem {BCO}  J. Berndt, S. Console, C. E. Olmos, Submanifolds and holonomy, Second edition. CRC Press, 2016.

\bibitem {Bs} A. L. Besse,  Einstein manifolds. Reprint of the 1987 edition. Classics in Mathematics. Springer-Verlag, Berlin, 2008.

\bibitem {BBB+} L. Bessi\`{e}res, G. Besson, M. Boileau, S.
Maillot and J. Porti, Geometrisation of 3-manifolds, Europ. Math.
Soc. 2010.

\bibitem {BBM} L. Bessi\`{e}res, G. Besson and S. Maillot, Ricci flow on
 open 3-manifolds and positive scalar curvature,  Geometry and Topology  15 (2011), 927-975.


\bibitem {BW} C. B\"{o}hm, B. Wilking, Manifolds with positive curvature operators are space forms, Ann. Math. 167 (2008), 1079-1097.


\bibitem {BMP} M. Boileau, S. Maillot, and J. Porti, Three-dimensional orbifolds and their geometric
structures, Soci\'{e}t\'{e} Math\'{e}matique de
France, Paris, 2003.



\bibitem {B02} F. Bonahon,  Geometric structures on 3-manifolds. Handbook of geometric topology, 93-164, North-Holland, Amsterdam, 2002.


\bibitem {BSi85} F. Bonahon, L. Siebenmann,  The classification of Seifert fibred 3-orbifolds. Low-dimensional topology (Chelwood Gate, 1982), 19-85. London Math. Soc. Lecture Note Ser., 95, Cambridge University Press, Cambridge, 1985.


\bibitem {BSi} F. Bonahon, L. C. Siebenmann, The characteristic toric splitting of irreducible compact 3-orbifolds,  Math. Ann. 278 (1987), 441-479.

\bibitem {B93} J. Borzellino, Orbifolds of maximal diameter, Indiana Univ. Math. J. 42 (1993),  37-53.

\bibitem {BB12} J. Borzellino, V. Brunsden,  Elementary orbifold differential topology, Topology Appl. 159 (2012), no.17, 3583-3589.

\bibitem {BB13} J. Borzellino, V. Brunsden, The stratified structure of spaces of smooth orbifold mappings, Commun. Contemp. Math. 15 (2013), no.5, 1350018, 37 pp.


\bibitem {B72} G. E. Bredon,  Introduction to compact transformation groups. Pure and Applied Mathematics, Vol. 46. Academic Press, New York-London, 1972.

\bibitem {B09} S. Brendle,  A generalization of Hamilton's differential Harnack inequality for the Ricci flow. J. Differential Geom. 82 (2009), no. 1, 207-227.

\bibitem {B10a} S. Brendle,  Einstein manifolds with nonnegative isotropic curvature are locally symmetric. Duke Math. J. 151 (2010), no. 1, 1-21.

\bibitem {B10} S. Brendle, Ricci flow and the sphere theorem. Graduate Studies in Mathematics, 111. American Mathematical Society, Providence, RI, 2010.

\bibitem {B18} S. Brendle, Ricci flow with surgery in higher dimensions, Ann. of Math. (2)
187 (2018) no. 1, 263-299.

\bibitem {B19} S. Brendle, Ricci flow with surgery on manifolds with positive isotropic curvature, Ann. of Math. 190 (2019), 465-559.

\bibitem {BN1} S. Brendle, P. Daskalopoulos, K. Naff, N. Sesum, Uniqueness of compact ancient solutions to the higher dimensional Ricci flow, J. reine angew. Math. 795 (2023), 85-138.

\bibitem {BN2} S. Brendle, K. Naff, Rotational symmetry of ancient solutions to the Ricci flow in higher dimensions, arXiv:2005.05830;  Geom. Topol. 27(2023), no.1, 153-226.



\bibitem {BS} S. Brendle, R. Schoen,  Manifolds with 1/4-pinched curvature are space forms. J. Amer. Math. Soc. 22 (2009), no. 1, 287-307.

\bibitem {BS08} S. Brendle, R. Schoen,  Classification of manifolds with weakly 1/4-pinched curvatures. Acta Math. 200 (2008), no. 1, 1-13.

\bibitem {BJ}  Th. Br\"{o}cker, K. J\"{a}nich, Introduction to Differential Topology, Cambridge University Press, 1982.

\bibitem {BBI}  D. Burago, Y. Burago, S. Ivanov, A course in metric geometry, Grad. Stud. Math., 33, American Mathematical Society, Providence, RI, 2001.

\bibitem {BHZ}  T. Buttsworth, M. Hallgren, Y. Zhang, Canonical surgeries in rotationally invariant Ricci flow, arXiv:2201.09387;  Trans. Amer. Math. Soc. 377 (2024), no. 11,  7877-7944.


\bibitem {CaZ} H.-D. Cao, X.-P. Zhu, A complete proof of the
Poincar\'{e} and geometrization conjectures---application of
the Hamilton-Perelman theory of the Ricci flow, Asian J. Math. 10
(2006), 165-492; arXiv:math/0612069.



\bibitem {CH} B.-L. Chen, X.-T. Huang, Path-connectedness of the moduli spaces of metrics with positive isotropic curvature on four-manifolds,  Math. Ann.  366 (2016), 819-851.


\bibitem {CTZ} B.-L. Chen, S.-H. Tang and X.-P. Zhu,  Complete classification of
compact four-manifolds with positive isotropic curvature,  J. Diff.
Geom. 91 (2012), 41-80.

\bibitem {CZ06} B.-L. Chen, X.-P. Zhu, Uniqueness of the Ricci flow on complete noncompact manifolds, J. Diff. Geom.
74 (2006), 119-154.

\bibitem {CZ} B.-L. Chen, X.-P. Zhu, Ricci flow with surgery
on  four-manifolds with positive isotropic curvature, J. Diff. Geom.
74 (2006), 177-264.


\bibitem {CR} W. Chen, Y. Ruan,  Orbifold Gromov-Witten theory. Orbifolds in mathematics and physics (Madison, WI, 2001), 25-85, Contemp. Math., 310, Amer. Math. Soc., Providence, RI, 2002.


\bibitem {Ch} Z. Chen, Manifolds with positive isotropic curvature of dimension at least nine, arXiv:2410.21078.

\bibitem {CL} J. H. Cho, Y. Li,  Ancient solutions to the Ricci flow with isotropic curvature conditions, arXiv:2005.11866;   Math. Ann.  387 (2023), 1009-1041.


\bibitem {CC} B. Chow, S.-C. Chu, D. Glickenstein, C. Guenther, J. Isenberg, T. Ivey, D. Knopf, P. Lu, F. Luo, L. Ni, The Ricci flow: techniques and applications, Part I: Geometric aspects, American Mathematical Society, 2007.



\bibitem {C+} B. Chow, S.-C. Chu, D. Glickenstein, C. Guenther, J. Isenberg, T. Ivey, D. Knopf, P. Lu, F. Luo, L. Ni, The Ricci flow: techniques and applications, Part II: Analytic aspects, American Mathematical Society, 2008.

\bibitem {CK} B. Chow, D. Knopf,  The Ricci flow: an introduction,  American Mathematical Society,  2004.


\bibitem {CHK} D. Cooper, C. D. Hodgson, S. P. Kerckhoff,
Three-dimensional orbifolds and cone-manifolds,  Mathematical Society of Japan, 2000.

\bibitem {CGGK} D. Corro, K. Garcia, M. G\"{u}nther, J.-B. Korda{\ss}, Bundles with even-dimensional spherical space form as fibers and fiberwise quarter pinched Riemannian metrics, Proc. Amer. Math. Soc. 149 (2021), 5407-5416.



\bibitem {dR} G. de Rham,  Complexes \`{a} automorphismes et hom\'{e}omorphie diff\'{e}rentiable.  Ann. Inst. Fourier (Grenoble) 2 (1950), 51-67 (1951).


\bibitem {D08}  J. Dinkelbach, Equivariant Ricci flow with surgery,  PhD thesis, Ludwig-Maximilians-University, Munich 2008.


\bibitem {DL} J. Dinkelbach and B. Leeb, Equivariant Ricci flow with surgery and application to finite group
actions on geometric 3-manifolds, Geom.  Topol. 13 (2009), p. 1129-1173.


\bibitem {FGKO} T. Farrell, Z. Gang, D. Knopf, P. Ontaneda, Sphere bundles with 1/4-pinched fiberwise metrics, Trans. Amer. Math. Soc. 369 (2017), no. 9, 6613-6630.


\bibitem {F} A. M. Fraser, Fundamental groups of manifolds with positive isotropic curvature,   Ann. of Math. 158 (2003), 345-354.


\bibitem {G} M. Gromov, Positive curvature, macroscopic dimension, spectral gaps and higher
signatures, Functional analysis on the eve of the 21st century, Vol. II (New
Brunswick 1993), 1-213, Progr. Math., 132, Birkh\"{a}user.

\bibitem {H93} R. Hamilton, The Harnack estimate for the Ricci flow, J. Diff. Geom. 37 (1993), 225-243.

\bibitem {H95} R. Hamilton, A compactness property for solutions of the Ricci flow, Amer. J.
Math. 117 (1995), no.3, 545-572.

\bibitem {H95b}  R. Hamilton,  The formation of singularities in the Ricci flow. Surveys in differential geometry, Vol. II, 7-136, International Press, Cambridge, MA, 1995.



\bibitem {H97} R. Hamilton, Four-manifolds with positive isotropic curvature, Comm. Anal. Geom. 5 (1997), 1-92; also in Collected Papers on Ricci flow, 342-407, edited by H. D. Cao, B. Chow, S. C. Chu and S. T. Yau, International Press 2003.

\bibitem {Hi} M. W. Hirsch,  Differential topology. Corrected reprint of the 1976 original. Graduate Texts in Mathematics, 33. Springer-Verlag, New York, 1994.

\bibitem {Hu13} H. Huang, Ricci flow on open 4-manifolds with positive isotropic
curvature, J. Geom. Anal. 23 (2013), no.3, 1213-1235.

\bibitem {Hu15} H. Huang, Four-orbifolds with positive isotropic curvature,  Comm. Anal. Geom. 23 (2015), no. 5, 951-991; arXiv:1107.1469.


\bibitem {Hu23} H. Huang, Classification of compact manifolds with positive isotropic curvature, arXiv:2305.18154; submitted to a journal.

\bibitem {Hu23b} H. Huang, Open manifolds with uniformly positive isotropic curvature, arXiv:2311.15825.

\bibitem {HS}  G. Huisken, C. Sinestrari, Mean curvature flow with surgeries of two-convex hypersurfaces, Invent. Math. 175 (2009), no.1, 137-221.


\bibitem {KL08}  B. Kleiner, J. Lott, Notes on Perelman's papers, Geom. Topol. 12 (2008), 2587-2855;  arXiv:0605667v5.


\bibitem {KL}  B. Kleiner, J. Lott, Geometrization of three-dimensional orbifolds via Ricci flow,  Ast\'{e}risque No. 365 (2014), 101-177;  arXiv:1101.3733v3.

\bibitem {KN} S. Kobayashi, K. Nomizu,  Foundations of differential geometry. Vol I. Interscience Publishers,  1963.


\bibitem {K} A. Kosinski, Differential Manifolds, Academic Press, 1993.

\bibitem {Ko} B. Kotschwar, Backwards uniqueness for the Ricci flow, Intern. Math.  Res. Not.  2010, Issue 21, 4064-4097.

\bibitem {Le} J. M. Lee,  Introduction to smooth manifolds. Second edition. Graduate Texts in Mathematics, 218. Springer, New York, 2013.

\bibitem {LZ} X. Li, Y. Zhang, Ancient solutions to the Ricci flow in higher dimensions, arXiv:1812.04156; Commun. Anal. Geom.  30 (2022), no.9, 2011-2048.

\bibitem {L} P. Lu, A compactness property for solutions of
the Ricci flow on orbifolds, Amer. J. Math.  123 (2001), 1103-1134.

\bibitem {MM} M. Micallef, J.D. Moore, Minimal two-spheres
and the topology of manifolds with positive curvature on totally
isotropic two-planes, Ann. Math. (2) 127 (1988), 199-227.



\bibitem {MW} M. J. Micallef, M. Y. Wang, Metrics with nonnegative isotropic curvature. Duke Math. J. 72 (1993), no. 3, 649-672.

\bibitem {M69} J. Milnor, Morse theory, Princeton University Press, 1969.


\bibitem {M07} J. Milnor,  Collected papers of John Milnor. III. Differential topology. American Mathematical Society, Providence, RI, 2007.

 \bibitem {MT} J. Morgan, G. Tian, Ricci flow and the Poincar\'{e} conjecture. Clay Mathematics Monographs, 3. American Mathematical Society, Providence, RI; Clay Mathematics Institute, Cambridge, MA, 2007.


 \bibitem {Mu} A. Mukherjee, Differential topology. Second edition. Hindustan Book Agency, New Delhi; Birkh\"{a}user/Springer, Cham, 2015.


 \bibitem {N} H. Nguyen, Isotropic curvature and the Ricci flow, Internat. Math. Res. Notices,
2010, no. 3, 536-558.


 \bibitem {NW} L. Ni, B. Wilking,  Manifolds with 1/4-pinched flag curvature, Geom. Funct. Anal. 20 (2010), 571-591.


\bibitem {P0} G. Perelman, A proof of the soul conjecture of Cheeger and Gromoll, J. Diff. Geom. 40 (1994), 209-212.


\bibitem {P1} G. Perelman, The entropy formula for the Ricci flow and its geometric applications,
arXiv:math.DG/0211159.


\bibitem {P2} G. Perelman, Ricci flow with surgery on three-manifolds, arXiv:math.DG/0303109.



\bibitem {Pe} P. Petersen,  Riemannian geometry, Third edition,  Grad. Texts in Math., 171,  Springer,  2016.


\bibitem {R} J. Ratcliffe, Foundations of hyperbolic manifolds, Third edition, GTM 149, Springer, 2019.


\bibitem {Sa} S. Salamon, Riemannian geometry and holonomy groups. Pitman Research Notes in Mathematics Series, 201. Longman Scientific and Technical, 1989.

\bibitem {S} R. Schoen, Minimal submanifolds in higher codimension, Mat. Contemp. 30 (2006), 169-199.

\bibitem {Sc} P. Scott,  The geometries of 3-manifolds. Bull. London Math. Soc. 15 (1983), no. 5, 401-487.


\bibitem {Si} J. Simons, On transitivity of holonomy systems, Ann. Math. 76 (1962), 213-234.


\bibitem {So} I. Solonenko,  answer to the question``How does one prove that $\text{Isom}(\mathbb{S}^2 \times \mathbb{R}) = \text{Isom}(\mathbb{S}^2)  \times \text{Isom}(\mathbb{R})$?"  in math.stackexchange.com



\bibitem {T} W. Thurston, The geometry and topology of three-manifolds, American Mathematical Society, Providence,  RI, 2022.

\bibitem {To} P. M. Topping, Ricci flow and PIC1, Surv. Differ. Geom., 27, 189-211,  International Press, Somerville, MA, 2024.

\bibitem {W} C. T. C. Wall,  Differential topology. Cambridge Studies in Advanced Mathematics, 156. Cambridge University Press, Cambridge, 2016.

\bibitem {Wi} B. Wilking, A Lie algebraic approach to Ricci flow invariant curvature conditions and
Harnack inequalities, J. Reine Angew. Math. 679 (2013), 223-247.


\bibitem {Wo} J. A. Wolf, Spaces of constant curvature. Sixth edition. AMS Chelsea Publishing, Providence, RI, 2011.


\bibitem {Y} D. Yeroshkin,  Riemannian orbifolds with non-negative curvature, Ph.D. Thesis, University of Pennsylvania, 2014.

\bibitem {Z11} Q. S. Zhang,  Sobolev inequalities, heat kernels under Ricci flow, and the Poincar\'{e} conjecture. CRC Press, Boca Raton, FL, 2011.


\end{thebibliography}
\end{document}